\documentclass[11pt]{amsart}%
\usepackage{ae,aecompl}
\usepackage{chngcntr}
\usepackage{graphicx}
\usepackage[all,cmtip]{xy}
\usepackage{amsmath, amscd}
\usepackage{amsthm}
\usepackage{amssymb}
\usepackage{amsfonts}
\usepackage{qsymbols}
\usepackage{latexsym}
\usepackage{mathrsfs}
\usepackage{cite}
\usepackage{color}
\usepackage{url}
\usepackage{enumerate}
\usepackage[utf8]{inputenc}
\usepackage[T1]{fontenc}
\usepackage{esint}
\usepackage{mathtools}
\usepackage{bbm}
\usepackage{verbatim}
\usepackage[colorlinks, citecolor = blue, pagebackref]{hyperref}
\usepackage[colorinlistoftodos,prependcaption,textsize=small]{todonotes}
\usepackage{accents}
\usepackage{amsmath}%
\providecommand{\U}[1]{\protect\rule{.1in}{.1in}}
\allowdisplaybreaks
\newcommand {\R}{{\mathbb R}}

\newcommand {\vanish}[1]{\relax}
\DeclarePairedDelimiter\floor{\lfloor}{\rfloor}

\newtheorem{theorem}{Theorem}[section]

\newtheorem{proposition}[theorem]{Proposition}
\newtheorem{corollary}[theorem]{Corollary}

\theoremstyle{definition}

\newtheorem{remark}[theorem]{Remark}
\newtheorem{example}[theorem]{Example}

\renewcommand{\asymp}{\eqsim}
\numberwithin{equation}{section}
\protected
\def\ignorethis#1\endignorethis{}
\let\endignorethis\relax

\begin{document}
\title[New Sobolev inequalities]{New Brezis--Van Schaftingen--Yung--Sobolev type inequalities connected with
maximal inequalities and one parameter families of operators}
\author{Oscar Dom\'inguez}
\address{O. Dom\'inguez, Departamento de An\'alisis Matem\'atico y Matem\'atica
Aplicada, Facultad de Matem\'aticas\\
Universidad Complutense de Madrid\\
Plaza de Ciencias 3\\
28040 Madrid\\
Spain}
\email{oscar.dominguez@ucm.es}
\author{Mario Milman}
\address{M. Milman, Instituto Argentino de Matematica\\
Buenos Aires\\
Argentina}
\email{mario.milman@icloud.com}
\urladdr{https://sites.google.com/site/mariomilman/}
\thanks{The first named author is supported in part by MTM2017-84058-P (AEI/FEDER, UE)}

\begin{abstract}
Motivated by the recent characterization of Sobolev spaces due to Brezis--Van
Schaftingen--Yung we prove new weak-type inequalities for one parameter
families of operators connected with mixed norm inequalities. The novelty of
our approach comes from the fact that the underlying measure space
incorporates the parameter as a variable. We also show that our framework can
be adapted to treat related characterizations of Sobolev spaces obtained
earlier by Bourgain--Nguyen. The connection to classical and fractional order
Sobolev spaces is shown through the use of generalized Riesz potential spaces
and the Caffarelli--Silvestre extension principle. Higher order inequalities
are also considered. We indicate many examples and applications to PDE's and
different areas of Analysis, suggesting a vast potential for future research.
In a different direction, and inspired by methods originally due to Gagliardo
and Garsia, we obtain new maximal inequalities which combined with mixed norm
inequalities are applied to obtain Brezis--Van Schaftingen--Yung type
inequalities in the context of Calder\'{o}n--Campanato spaces. In particular,
Log versions of the Gagliardo--Brezis--Van Schaftingen--Yung spaces are
introduced and \ compared with corresponding limiting versions of
Calder\'{o}n--Campanato spaces, resulting in a sharpening of recent
inequalities due to Crippa--De Lellis and Bru\'{e}--Nguyen.

\end{abstract}
\maketitle
\tableofcontents

\section{Introduction\label{sec:introduction}}

In their recent work \cite{Brezis, BrezisSY}, Brezis--Van Schaftingen--Yung obtained a
striking novel way to recover the $L^{p}$ norm of the gradient of a Sobolev
function using a weak-type version of the classical Gagliardo seminorms. Let
$N\geq1,1\leq p<\infty$, then (cf. \cite[Theorem 1.1]{Brezis}) \footnote{Here the symbol
$f\asymp g$ indicates the existence of a universal constant $c>0$ (independent
of all parameters involved) such that $(1/c)f\leq g\leq cf$. Likewise the
symbol $f\lesssim g$ will mean that there exists a universal constant $c>0$
(independent of all parameters involved) such that $f\leq cg$.}\footnote{It
was recently shown by Poliakovsky \cite[Theorem 1.3]{poliakovsky} that \eqref{I1} holds, more generally, for functions $f\in W_{p}^{1}({\mathbb{R}%
}^{N})$ if $p>1$ or $f\in \text{BV}({\mathbb{R}}^{N})$ if $p=1.$
Conversely, if $f \in L^p(\mathbb{R}^N)$ then the finiteness of the left-hand side of \eqref{I1} implies that
$f\in W_{p}^{1}({\mathbb{R}}^{N})$ if $p>1$ or $f\in \text{BV}({\mathbb{R}%
}^{N})$ if $p=1.$ Similar results are valid in smooth domains.}
\begin{equation}
\bigg\|\frac{f(x)-f(y)}{|x-y|^{\frac{N}{p}+1}}\bigg\|_{L(p,\infty
)({\mathbb{R}}^{N}\times{\mathbb{R}}^{N})}\asymp\Vert\nabla f\Vert
_{L^{p}({\mathbb{R}}^{N})},\quad f\in C_{c}^{\infty}({\mathbb{R}}^{N}).
\label{I1}%
\end{equation}
Here $L(p,\infty)(X,m)$ is the weak $L^{p}$ space on the measure space
$(X,m),$ defined by the condition
\begin{equation}
\left\Vert f\right\Vert _{L(p,\infty)(X,m)}^{p}=\sup_{\lambda>0}\lambda
^{p}m(\{x\in X:\left\vert f(x)\right\vert >\lambda\})<\infty.
\label{defweakLp}%
\end{equation}
Then \eqref{I1} reads as\footnote{Here and in what follows, $\mathcal{L}^{N}$
denotes the Lebesgue measure on ${\mathbb{R}}^{N}$ and $\mathcal{L}%
^{1}=\mathcal{L}$.}%

\begin{equation}
\sup_{\lambda>0}\lambda^{p}\mathcal{L}^{2N}(E_{\lambda})\asymp\Vert\nabla
f\Vert_{L^{p}({\mathbb{R}}^{N})}^{p} \label{l1*}%
\end{equation}
where
\[
E_{\lambda}=\bigg\{(x,y)\in{\mathbb{R}}^{N}\times{\mathbb{R}}^{N}%
:\frac{\left\vert f(x)-f(y)\right\vert }{\left\vert x-y\right\vert ^{\frac
{N}{p}+1}}>\lambda\bigg\}.
\]
As an application (cf. \cite{BrezisSY}) these authors consider some endpoint inequalities that fail
for the classical fractional Gagliardo seminorms, and prove alternative
inequalities using the \textquotedblleft weak-type Gagliardo
seminorms\textquotedblright\ that appear on the left-hand sides of (\ref{I1})
and \eqref{l1*}. Furthermore, it is also shown in \cite[Theorem 1.2]{Brezis} that the lower
bound in \eqref{l1*} can be considerably sharpened\footnote{The formula \eqref{I1'} was recently extended to $f \in W^1_p(\mathbb{R}^N), \, p \geq 1,$ in \cite[Corollary 1.2]{poliakovsky}.}
\begin{equation}
\lim_{\lambda\rightarrow+\infty}\lambda^{p}\mathcal{L}^{2N}(E_{\lambda}%
)=\frac{1}{N}k(p,N)\Vert\nabla f\Vert_{L^{p}({\mathbb{R}}^{N})}^{p},\quad f\in
C_{c}^{\infty}({\mathbb{R}}^{N}), \label{I1'}
\end{equation}
where $k(p,N)$ is an explicit computable constant\footnote{$k(p,N)=\int%
_{\mathbb{S}^{N-1}}\left\vert \langle e,w\rangle\right\vert ^{p}%
\mathrm{d}\sigma^{N-1}(w),$ where $\mathbb{S}^{N-1}$ is the unit sphere in
$\mathbb{R}^{N}$ and $e$ is any unit vector in $\mathbb{R}^{N}.$} that depends
only on $p$ and $N$. 

In what follows we shall use the classical notation for the Gagliardo spaces:
for $s\in(0,1], \, 1 \leq p < \infty,$ we let
\begin{equation}
\left\Vert f\right\Vert _{W^{s,p}(\mathbb{R}^{N})}:=\bigg\|\frac
{f(x)-f(y)}{|x-y|^{\frac{N}{p}+s}}\bigg\|_{L^{p}({\mathbb{R}}^{N}%
\times{\mathbb{R}}^{N})}. \label{DefG}%
\end{equation}
Moreover, due to their central role in our development, we single out the
$BSY_{p}^{s}$ spaces (Brezis--Van Schaftingen--Yung spaces), defined through
the functional%
\begin{equation}
\text{ }\left\Vert f\right\Vert _{BSY_{p}^{s}(\mathbb{R}^{N})}:=\bigg\|\frac
{f(x)-f(y)}{|x-y|^{\frac{N}{p}+s}}\bigg\|_{L(p,\infty)({\mathbb{R}}^{N}%
\times{\mathbb{R}}^{N})},\quad1\leq p<\infty,\quad s\in(0,1]. \label{pppp}%
\end{equation}

The use of Gagliardo weak-type functionals on product domains in this context
is apparently new. An earlier closely related result due to
Bourgain--Brezis--Mironescu \cite{Bourgain00}, provided a different method to
compute the $L^{p}$ norm of the gradient of a function, via limits of the
classical Gagliardo seminorms: Let $1\leq p<\infty,f\in\mathring{W}_{p}%
^{1}({\mathbb{R}}^{N}),$ then%
\begin{equation}
\lim_{s\rightarrow1-}(1-s)\left\Vert f\right\Vert _{W^{s,p}(\mathbb{R}^{N}%
)}^{p}=\frac{k(p,N)}{p}\Vert\nabla f\Vert_{L^{p}({\mathbb{R}}^{N})}^{p}.
\label{I2}%
\end{equation}
The presence of the factor $(1-s)$ in (\ref{I2}) serves to mitigate the
behavior of the Gagliardo seminorms as $s\rightarrow1.$ In fact, when $s=1$ it
is easy to see that the corresponding Gagliardo space $W^{1,p}(\mathbb{R}%
^{N})$, does not coincide with the homogeneous Sobolev space $\mathring{W}%
_{p}^{1}(\mathbb{R}^{N}),$ and indeed (cf. \cite{Bourgain00}),%
\[
\left\Vert f\right\Vert _{W^{1,p}(\mathbb{R}^{N})}<\infty\text{ implies
}f=\text{constant.}%
\]

In particular, the result has been extended to more general domains, it has
been applied to formulate conditions that imply constancy of functions, it has
been used to investigate endpoint Sobolev-type inequalities, etc. (cf.
\cite{ponce} and the references therein). Furthermore, one can naturally place
(\ref{I2}) in the more general framework of scales of interpolation spaces
(cf. \cite{Milman}). In particular, this makes it possible to give a unified
approach\footnote{In particular, the interpolation point of view unifies the
Bourgain--Brezis--Mironescu formula ($s=1$) with the Maz'ya--Shaposhnikova (cf.
\cite{Maz}) formula ($s=0$). In this concern we note parentetically that the
mixed norm inequalities of this paper easily yield%
\[
\left\Vert \frac{f(x)-f(y)}{\left\vert x-y\right\vert ^{\frac{N}{p}}%
}\right\Vert _{L(p,\infty)({\mathbb{R}}^{N}\times{\mathbb{R}}^{N})}%
\lesssim \left\Vert f\right\Vert _{L^{p}({\mathbb{R}}^{N})},\quad1\leq
p<\infty.
\]
} to related formulae for generalized Sobolev-type spaces in different
contexts, including Carnot groups, Besov spaces associated with semigroups, or
even non commutative versions of Sobolev spaces (cf. \cite{Bus}, \cite{Maa},
\cite{XiXio}, and the references therein).

These characterizations of Sobolev seminorms, and many of its variants, have
been the subject of intense research interest. Indeed, during the review
process of our work, the referees called our attention to the characterization
of Sobolev spaces obtained by Nguyen \cite{ngu06} and Bourgain--Nguyen
\cite{bn06} (cf. also the references therein). In these papers a somewhat
related set of functionals was used to characterize the $L^{p}$ norms of
gradients. Let $1\leq p<\infty,\delta>0,f\in L^{p}(\mathbb{R}^{N}),$ and
define
\begin{equation}\label{Idelta}
I_{\delta}(f)=\int\int_{\{\left\vert f(x)-f(y)\right\vert >\delta\}}%
\frac{\delta^{p}}{\left\vert x-y\right\vert ^{N+p}} \, \mathrm{d}x \, \mathrm{d}y.
\end{equation}
Suppose that $1<p<\infty,$ and $f\in W^{1,p}(\mathbb{R}^{N}),$ then (cf.
\cite{ngu06}, \cite{bn06}) there exist positive constants\footnote{The
precise value of the constants is known. In particular, $c(p, N) = \frac{k(p, N)}{p}$.} $C(p,N),$ $c(p,N),$ such that%
\begin{equation}
\sup_{\delta>0}I_{\delta}(f)\leq C(p,N)\int_{\mathbb{R}^{N}}\left\vert \nabla
f(x)\right\vert ^{p} \, \mathrm{d}x, \label{nb1}%
\end{equation}
and%
\begin{equation}
\lim_{\delta\rightarrow0^{+}}I_{\delta}(f)=c(p,N)\int_{\mathbb{R}^{N}%
}\left\vert \nabla f(x)\right\vert ^{p} \, \mathrm{d}x. \label{nb2}%
\end{equation}
Moreover (cf. \cite{ngu06}, \cite{bn06}, \cite{bn18}), if $p\geq1,$ $f\in
L^{p}(\mathbb{R}^{N}),$ and
\[
\liminf_{\delta\rightarrow0^{+}}I_{\delta}(f)<\infty,
\]
then $f\in \text{BV}(\mathbb{R}^{N})$ if $p=1$ and $f\in W^{1,p}(\mathbb{R}^{N})$ if
$p>1,$ and
\begin{equation}
\liminf_{\delta\rightarrow0^{+}}I_{\delta}(f)\geq c(p,N)\int_{\mathbb{R}^{N}%
}\left\vert \nabla f(x)\right\vert ^{p} \mathrm{d}x. \label{nb3}%
\end{equation}
Following questions and suggestions of the reviewers, we have added an
Appendix, which can be read independently from the rest of the paper, where we
show, using a variant of Theorem \ref{ThmMaximal} below (cf. Theorem
\ref{ThmMaximal2} in the Appendix), that our general framework can be used to
provide a unified abstract treatment of the Brezis--Van Schaftingen--Yung and
the Bourgain--Nguyen formulae.

Based on all these considerations, one can expect that (\ref{I1}) will have a
considerable impact\footnote{We should also note here that another striking
feature of these results is that one could apriori use the left hand sides of
(\ref{I1}) or (\ref{I2}) to \textit{define} Sobolev spaces in contexts where
there is no apriori notion of a gradient.}.

It is natural to ask if results similar to \eqref{I1} can be obtained for
other intermediate spaces. In this paper we shall focus on two different type
of spaces. The Riesz potential spaces, which are defined using fractional
derivatives, and the Calder\'{o}n--Campanato spaces, defined using suitable
maximal operators instead. The results that we obtain for each of these
classes of spaces are very different in nature and specific toolboxes had to
be tailored to deal with each of them.

When dealing with Riesz potential spaces we have to contend with the fact that
$(-\Delta)^{s},s\in(0,1),$ is a non-local operator and thus the localization
techniques, that were successfully applied in \cite{Brezis} to establish
(\ref{I1}), may fail. To avoid this problem we apply the celebrated
Caffarelli--Silvestre extension theorem. We shall now explain our point of
view in more detail.

We consider the Riesz potential spaces $H^{2s,p}({\mathbb{R}}^{N})$ endowed
with%
\begin{equation}
\label{IntroRiesz}\left\Vert f\right\Vert _{H^{2s,p}({\mathbb{R}}^{N})}%
:=\Vert(-\Delta)^{s}f\Vert_{L^{p}({\mathbb{R}}^{N})},\quad s\in(0,1),\quad1 <
p<\infty.
\end{equation}
We then ask: What is an appropriate Brezis--Van Schaftingen--Yung condition
(cf. \eqref{pppp}) controlling $\|f\|_{H^{2 s,p}(\mathbb{R}^{N})}$? We first
observe that the non-local character of $(-\Delta)^{s}$ can be overcome by
using the Caffarelli--Silvestre extension theorem \cite{Caffarelli}. Recall
that this extension allows for problems involving $(-\Delta)^{s}$ on
$\mathbb{R}^{N},$ to be transformed into local PDE's of degenerate type in the
upper half-space $\mathbb{R}_{+}^{N+1}=\mathbb{R}^{N}\times(0,\infty).$ More
precisely, for a given $f\in\mathcal{S}(\mathbb{R}^{N}),$ we consider
$u_{(-\Delta)^{s}}$, the unique solution on $\mathbb{R}_{+}^{N+1}$ to%

\begin{equation}
\left\{
\begin{array}
[c]{ll}%
\Delta u+\frac{1-2s}{t}u_{t}+u_{tt}=0 & \quad\text{in}\quad\mathbb{R}%
_{+}^{N+1},\\
u(x,0)=f(x) & \quad\text{on}\quad\mathbb{R}^{N}.
\end{array}
\right.  \label{11111}%
\end{equation}
This solution can be expressed in terms of the Poisson kernels,
\[
P_{(-\Delta)^{s}}(x,t)= C_{N,s} \, \frac{t^{2s}}{(|x|^{2}+t^{2})^{\frac{N}%
{2}+s}},\quad x\in{\mathbb{R}}^{N},\quad t>0,
\]
that is,
\begin{equation}
u_{(-\Delta)^{s}}(x,t)=P_{(-\Delta)^{s}}[f](x,t)=(P_{(-\Delta)^{s}}%
(\cdot,t)\ast f)(x)=\int_{\mathbb{R}^{N}}P_{(-\Delta)^{s}}%
(x-y,t)f(y)\,\mathrm{d}y. \label{CSExplicit}%
\end{equation}
In particular, if $s=1/2$ then $u_{(-\Delta)^{1/2}}=P_{(-\Delta)^{1/2}}[f]$
coincides with $P[f]$, the classical Poisson extension of $f$ to the upper
half-space (cf. \cite[p. 61]{Stein70b}):
\begin{equation}
P[f](x,t)= P_{t}[f](x)= C_{N} \, \int_{\mathbb{R}^{N}}\frac{t}{(|x-y|^{2}%
+t^{2})^{\frac{N+1}{2}}}f(y)\,\mathrm{d}y. \label{ooo}%
\end{equation}
Then there exists a positive constant $\mu_{s}$ such that,
\begin{equation}
(-\Delta)^{s}f(x)=-\mu_{s}\lim_{t\rightarrow0+}\frac{P_{(-\Delta)^{s}%
}[f](x,t)-f(x)}{t^{2s}}\quad(\text{in the sense of }L^{\infty}(\mathbb{R}%
^{N})). \label{combinada}%
\end{equation}
In this setting our version of the Brezis--Van Schaftingen--Yung type
conditions corresponds naturally to controlling the weak $L(p,\infty
)(\mathbb{R}_{+}^{N+1})$ norm of the expression $\frac{P_{(-\Delta)^{s}%
}[f](x,t)-f(x)}{t^{2s}}.$ Indeed, combining (\ref{combinada}) and Theorem
\ref{ThmMaximal}(ii) we show (cf. Theorem \ref{ThmCSBSY} below)%
\begin{equation}
\left\Vert f\right\Vert _{H^{2s,p}(\mathbb{R}^{N})}\leq\mu_{s}\bigg\|\frac
{P_{(-\Delta)^{s}}[f](x,t)-f(x)}{t^{2s+\frac{1}{p}}}\bigg\|_{L(p,\infty
)(\mathbb{R}_{+}^{N+1})}. \label{I5}%
\end{equation}
In particular, an upper estimate for the $L^{p}$ norm of the gradient in terms
of the Poisson extension \eqref{ooo} can be recovered by letting $s=1/2,$ and
using the well-known formula \textquotedblleft$\left\vert (-\Delta
)^{1/2}\right\vert =\left\vert \nabla\right\vert $\textquotedblright%
\thinspace\ for Riesz transforms,
\begin{equation}
\left\Vert \nabla f\right\Vert _{L^{p}(\mathbb{R}^{N})}\lesssim\left\Vert
\frac{P[f](x,t)-f(x)}{t^{1+\frac{1}{p}}}\right\Vert _{L(p,\infty
)(\mathbb{R}_{+}^{N+1})},\quad1<p<\infty. \label{tx}%
\end{equation}
Concerning \eqref{tx}, characterizations of function spaces of smooth
functions via Poisson kernels owe much to the pioneering work of Taibleson
\cite{Taibleson1, Taibleson2}; cf. also \cite{Dominguez} for a more recent
account on this topic. As a prototype of these characterizations, we find
\begin{equation}
\Vert\Delta f\Vert_{L^{p}(\mathbb{R}^{N})}\asymp\Big\|\Big\|\frac{(P_{t}
-\text{id})^{2}f}{t^{2}}\Big\|_{L_{x}^{p}({\mathbb{R}}^{N})}\Big\|_{L_{t}%
^{\infty}(0,\infty)},\quad1<p<\infty. \label{tx1}%
\end{equation}
Note that the right-hand sides of \eqref{tx} and \eqref{tx1} are of a
completely different nature and, in particular, the function norms in
\eqref{tx} are weak $L^{p}$ norms on the product space ${\mathbb{R}}_{+}%
^{N+1}$, while the function norms in \eqref{tx1} define the mixed norm spaces
$L_{t}^{\infty}((0,\infty),L_{x}^{p}({\mathbb{R}}^{N}))$.

The estimates \eqref{I5} and \eqref{tx} can be examined within the more
general framework of non local hypoelliptic operators in the sense of
H\"{o}rmander \cite{Hormander}. Consider
\begin{equation}
\mathscr{K}u:=\mathscr{A}u-u_{t}:=\text{tr}(Q\nabla^{2}u)+\langle Bx,\nabla
u\rangle-u_{t}=0\quad\text{in}\quad{\mathbb{R}}_{+}^{N+1} \label{Hor}%
\end{equation}
where ${\small Q=Q}^{\star}{\small \geq0}$ \ and{\small \ }${\small B}$ \ are
real $N\times N$ matrices with constant coefficients, and the hypoelliptic
condition
\begin{equation}
\text{det}\,K(t)>0\quad\text{for all}\quad t>0 \label{Hypoelliptic}%
\end{equation}
is satisfied, where
\[
K(t)=\frac{1}{t}\int_{0}^{t}e^{uB}Qe^{uB^{\ast}}\,\mathrm{d}u.
\]
As is well known, for suitable choices of $Q$ and $B,$ the equation
(\ref{Hor}) contains many models of interest in mathematical physics such as
the standard heat equation ($Q=I_{N}$ and $B=0_{N}$), the classical
Ornstein--Uhlenbeck equation ($Q=I_{N}$ and $B=-I_{N}$), Kolmogorov equation,
Kolmogorov with friction, Kramers equation, etc. Recently, Garofalo--Tralli
\cite{Garofalo20} extended the Caffarelli--Silvestre theorem for fractional
powers of the Laplacian to fractional powers of operators associated with the
H\"{o}rmander equation, $(-\mathscr{K})^{s}$, and its diffusive part
$(-\mathscr{A})^{s}$ (cf. Section \ref{sec:Caff-Silv} below). For example, for
$(-\mathscr{A})^{s},$ the Caffarelli--Silvestre--Garofalo--Tralli extension
reads as follows: If $f\in\mathcal{S}({\mathbb{R}}^{N})$ then there exists a
function $u_{(-\mathscr{A})^{s}}\equiv u\in C^{\infty}({\mathbb{R}}_{+}%
^{N+1})$ such that%
\[
\left\{
\begin{array}
[c]{cc}%
t^{1-2s}(\mathscr{A}u+\frac{1-2s}{t}u_{t}+u_{tt})=0 & \text{in}\quad
\mathbb{R}_{+}^{N+1},\\
u(x,0)=f(x) & \text{on}\quad\mathbb{R}^{N}.
\end{array}
\right.
\]
The solution of this problem can be expressed in terms of $f$ via the
associated Poisson kernel $P_{(-\mathscr{A})^{s}}$. Then corresponding to
(\ref{combinada}) we have%
\[
(-\mathscr{A})^{s}f(x)=-\mu_{s}\lim_{t\rightarrow0+}\frac
{P_{(-\mathscr{A})^{s}}[f](x,t)-f(x)}{t^{2s}},\quad\text{in}\quad L^{\infty
}(\mathbb{R}^{N}).
\]
To formulate the corresponding Brezis--Van Schaftingen--Yung result we
consider the spaces $H_{\mathscr{A}}^{2s,p}({\mathbb{R}}^{N})$ (resp.
$H_{\mathscr{K}}^{2s,p}({\mathbb{R}}^{N})$) defined by means of replacing
$(-\Delta)^{s}$ in \eqref{IntroRiesz} by $(-\mathscr{A})^{s}$ (resp.
$(-\mathscr{K})^{s}$) and prove, \textit{mutatis mutandis, }the corresponding
versions of (\ref{I5}). We refer to Section \ref{sec:Caff-Silv} (cf. Theorems
\ref{ThmCSBSYGarofalo}, \ref{ThmCSBSYGarofalo2}, for the details). These ideas
can be also applied to other PDE's where Caffarelli--Silvestre type extensions
are available, e.g., using the extension theorem of \cite{Stinga} we can deal
in a similar fashion with the harmonic oscillator
\[
H=-\Delta+\left\vert x\right\vert ^{2}.
\]

Underlying these results is a general principle, see Theorem \ref{ThmMaximal}
below, valid for nonlinear semigroups of operators, $\{T_{t}:t>0\},$ which can
be formulated in terms of weak-type spaces defined on $\mathbb{R}_{+}^{N+1}$.
In fact, as we shall now explain, our results do not require \textit{apriori}
the semigroup property and work in abstract measure spaces.

Let $(X,m)$ be a $\sigma$-finite measure space, and let $\{T_{t}:t>0\}$ be a
one-parameter family of (not necessarily linear) operators on $L^{p}(X,m)$. A
central issue often associated to such families is to show, under suitable
assumptions (cf. \cite{Stein, Stein61, Stein76, Stein70b, Garsia, Torchinsky},
and the references therein), that the associated maximal operator $T^{\ast
}f:=\sup_{t>0}|T_{t}f|$ is bounded on $L^{p}(X,m),$
\begin{equation}
\left\Vert T^{\ast}f\right\Vert _{L^{p}(X,m)}\leq C_{p}\left\Vert f\right\Vert
_{L^{p}(X,m)},\quad1<p<\infty, \label{I3}%
\end{equation}
and is of weak-type for $p=1.$ Accepting temporarily the validity of
(\ref{I3}) the problem we tackle is how to proceed to prove weak-type
estimates, and limit theorems, on the product domain $(X\times\mathbb{R}%
_{+},m\times\mathcal{L}).$

More generally, let $\gamma>0,$ we shall consider the measure on
$\mathbb{R}_{+}$ defined by
\[
w_{\gamma}(A)=\int_{A}t^{\gamma-1} \, \mathrm{d} t
\]
for a measurable set $A$ in ${\mathbb{R}}_{+}$. In particular, setting
$\gamma=1$ one recovers $\mathcal{L}$. Then we have (cf. Section
\ref{sec:maximal})

\begin{theorem}
\label{ThmMaximal} Let $(X,m)$ be a $\sigma$-finite measure space, and let
$\{T_{t}:t>0\}$ be a one-parameter family of (not necessarily linear)
operators on $L^{p}(X,m),\,1\leq p<\infty$, furthermore, suppose that
$\gamma>0.$

\begin{enumerate}
[\upshape(i)]

\item Assume\footnote{In what follows we only consider operators such that
$T_{t}f(x)$ is measurable as a function of $(x,t).$}
\begin{equation}
T^{\ast}f\in L^{p}(X,m). \label{MaxAssump1+new}%
\end{equation}
Then
\begin{equation}
\sup_{\lambda>0}\lambda^{p}(m\times w_{\gamma})\bigg(\bigg\{(x,t)\in
X\times(0,\infty):\frac{|T_{t}f(x)|}{t^{\gamma/p}}>\lambda\bigg\}\bigg)\leq
\frac{1}{\gamma}\Vert T^{\ast}f\Vert_{L^{p}(X,m)}^{p}. \label{IX}%
\end{equation}
In particular, if $p>1$ and the maximal operator $T^{\ast}$ is bounded on
$L^{p}(X,m)$ (i.e., \eqref{I3} holds) then
\[
\sup_{\lambda>0}\lambda^{p}(m\times w_{\gamma})\bigg(\bigg\{(x,t)\in
X\times(0,\infty):\frac{|T_{t}f(x)|}{t^{\gamma/p}}>\lambda\bigg\}\bigg)\leq
\frac{C_{p}^{p}}{\gamma}\Vert f\Vert_{L^{p}(X,m)}^{p}.
\]

\item Suppose that
\begin{equation}
\lim_{t\rightarrow0+}T_{t}f(x)<\infty\quad\text{$m$-a.e.}\quad x\in X.
\label{MaxAssump1}%
\end{equation}
Then
\[
\frac{1}{\gamma}\,\Big\|\lim_{t\rightarrow0+}T_{t}f\Big\|_{L^{p}(X,m)}^{p}%
\leq\liminf_{\lambda\rightarrow\infty}\lambda^{p}(m\times w_{\gamma
})\bigg(\bigg\{(x,t)\in X\times(0,\infty):\frac{|T_{t}f(x)|}{t^{\gamma/p}%
}>\lambda\bigg\}\bigg).
\]

\end{enumerate}
\end{theorem}

In his Corona seminar presentation (cf. \cite{van}) Jean Van Schaftingen asked
if, like its counterpart (\ref{I2}), it is possible to place (\ref{I1}) and
\eqref{I1'} in a more general framework. One could consider Theorem
\ref{ThmMaximal} as part of the program\footnote{In this connection, the
extrapolation formula that results combining (\ref{I1}) and (\ref{I2})%
\[
\lim_{s\rightarrow1-}(1-s)^{1/p}\bigg\|\frac{f(x)-f(y)}{|x-y|^{\frac{N}{p}+s}%
}\bigg\|_{L^{p}({\mathbb{R}}^{N}\times{\mathbb{R}}^{N})}\asymp\bigg\|\frac
{f(x)-f(y)}{|x-y|^{\frac{N}{p}+1}}\bigg\|_{L(p,\infty)({\mathbb{R}}^{N}%
\times{\mathbb{R}}^{N})},
\]
is a surprising result that also demands explanations from the interpolation
theory community.
\par
In the same vein, but in a different direction, M. Cwikel (private
communication) pointed out a possible connection of (\ref{I1'}) and the exact
formulae for the canonical maximal seminorm which is dominated by the weak
$L^{1}$ quasinorm (cf. \cite{CF}).} to understand (\ref{I1}) and \eqref{I1'}
from a general point of view\footnote{In the Appendix we prove a closely
related result (cf. Theorem \ref{ThmMaximal2}) which is used to incorporate
(\ref{nb1}), (\ref{nb2}), (\ref{nb3}) to our framework.}. Apart from the
applications that we have outlined above, Theorem \ref{ThmMaximal} is
connected with many different topics in Analysis. We have collected many
examples and applications in Section \ref{SubsectionAppl}, and we invite the
reader to use the table of contents as a road map to the many different
examples it contains, which can be read independently.

It will be instructive to go over the simple proof of Theorem \ref{ThmMaximal}%
(i) now, because it illustrates a mechanism that we frequently use to go from
maximal norm inequalities to weak-type inequalities on product
domains\footnote{As pointed out by the referees, an even easier proof can be
obtained using Fubini's theorem directly. Indeed, in the Appendix we use this
approach to deal with Theorem \ref{ThmMaximal2}. We prefer the method
presented here since not only it is useful to deal with many related results
throughout the paper, but also it applies to $L(p,q)$ spaces, and, indeed, informs
the type of spaces one needs to use to extend the results beyond Lebesgue
spaces.}. A basic consequence of Fubini's theorem that will be useful in the
sequel, can be stated as follows. Consider a product of $\sigma$-finite
measure spaces $(X_{1}\times X_{2},m_{1}\times m_{2}),$ then we have (cf.
Section \ref{Subsection:Mixed}, Proposition \ref{PropMixed}):%
\begin{equation}
\left\Vert f\right\Vert _{L(p,\infty)(X_{1}\times X_{2},m_{1}\times m_{2}%
)}\leq\min\bigg\{\left\Vert \left\Vert f\right\Vert _{L(p,\infty)(X_{i}%
,m_{i})}\right\Vert _{L^{p}(X_{j},m_{j})}:\text{ }i\neq j\in
\{1,2\}\bigg\},\quad1\leq p<\infty. \label{I4}%
\end{equation}
To prove (\ref{IX}) under the assumptions of Theorem \ref{ThmMaximal}, we
consider the product space $X\times\mathbb{R}_{+}$ with measure $m\times
w_{\gamma}.$ According to \eqref{defweakLp}, one has that
\[
\sup_{\lambda>0}\lambda^{p}(m\times w_{\gamma})\bigg(\bigg\{(x,t)\in
X\times(0,\infty):\frac{|T_{t}f(x)|}{t^{\gamma/p}}>\lambda\bigg\}\bigg)=\Vert
t^{-\gamma/p}T_{t}f(x)\Vert_{L(p,\infty)(X\times{\mathbb{R}}_{+},m\times
w_{\gamma})}^{p}.
\]
Now, by (\ref{I4}), the lattice property of the weak $L^{p}$ space and the
fact that
\[
\Vert t^{-\gamma/p}\Vert_{L(p,\infty)({\mathbb{R}}_{+},w_{\gamma})}=\frac
{1}{\gamma^{1/p}},
\]
we successively have
\begin{align*}
\Vert t^{-\gamma/p}T_{t}f(x)\Vert_{L(p,\infty)(X\times{\mathbb{R}}_{+},m\times
w_{\gamma})}  &  \leq\left\Vert \Vert t^{-\gamma/p}T_{t}f(x)\Vert
_{L(p,\infty)({\mathbb{R}}_{+},w_{\gamma})}\right\Vert _{L^{p}(X,m)}\\
&  \leq\left\Vert \Vert t^{-\gamma/p}\Vert_{L(p,\infty)({\mathbb{R}}%
_{+},w_{\gamma})}T^{\ast}f(x)\right\Vert _{L^{p}(X,m)}\\
&  =\frac{1}{\gamma^{1/p}}\left\Vert T^{\ast}f\right\Vert _{L^{p}(X,m)}.
\end{align*}
For the proof of Theorem \ref{ThmMaximal}(ii) we refer to Section
\ref{sec:maximal}.

An easy consequence of Theorem \ref{ThmMaximal} is that, under its
assumptions, we can recover the $L^{p}$ norms from weak $L^{p}$ norms as
follows (cf. Corollary \ref{CorMaximal}): If $1<p<\infty$ then
\begin{align}
\Vert f\Vert_{L^{p}(X,m)}^{p}  &  \leq\lim_{\lambda\rightarrow\infty}%
\lambda^{p}(m\times\mathcal{L})\bigg(\bigg\{(x,t)\in X\times(0,\infty
):\frac{|T_{t}f(x)|}{t^{1/p}}>\lambda\bigg\}\bigg)\nonumber\\
&  \leq\sup_{\lambda> 0}\lambda^{p}(m\times\mathcal{L})\bigg(\bigg\{(x,t)\in
X\times(0,\infty):\frac{|T_{t}f(x)|}{t^{1/p}}>\lambda\bigg\}\bigg) \leq
C_{p}^{p}\,\Vert f\Vert_{L^{p}(X,m)}^{p}. \label{lim1}%
\end{align}

The boundedness of the maximal operator $T^{\ast}$ is usually not available to
us when $p=1,$ and our arguments do not work directly when only weak-type
$(1,1)$ holds. However, for special families of operators $\{T_{t}\}$, we can
still obtain results like (\ref{lim1}) for $p=1$ with explicit values of the
equivalence constants. In Section \ref{sec:p}, this is achieved under the
assumption that for a dense class of functions, $\{T_{t} : t > 0\}$
satisfies\footnote{For interpolation theory aficionados we mention that
conditions of this type are connected with the computation of suitable $K$ and
$E$ functionals; cf. \cite{Jawerth}.} (cf. Theorem \ref{ThmPh}),
\[
\Vert T_{t}f-f\Vert_{L^{\infty}(X,m)}\leq C_{f}\,t^{1/p}\text{ for all }t>0
\]
where $C_{f}$ is a positive constant which may depend on $f$. Note that this
assumption means that $f$ is a H\"older--Lipschitz continuous function of
order $1/p$ with norm $C_{f}$. In particular, we recover the formula
\eqref{I1'} with $N=1$ (cf. Subsection \ref{sec:p}: Examples \ref{exbsyp=1},
\ref{exbysp=2} and also Subsection \ref{ExampleBSY}), that is,
\begin{equation}
\label{fl1}\lim_{\lambda\to\infty}\lambda^{p} \mathcal{L}^{2}
\bigg( \bigg\{(x,t)\in{\mathbb{R}} \times(0,\infty):\frac{\left\vert
f(x+t)-f(x)\right\vert }{t^{\frac{1}{p} + 1}} >\lambda\bigg\} \bigg) = \Vert
f^{\prime}\Vert_{L^{p}({\mathbb{R}})}^{p}, \quad p \geq1.
\end{equation}

Still dealing with $p=1$, we also address the counterpart of \eqref{fl1} which
is obtained by replacing $\lim_{\lambda\rightarrow\infty}$ by $\sup
_{\lambda>0}$, i.e., we investigate the equivalence provided in \eqref{l1*}.
In our notation, the approach of \cite{Brezis} related to \eqref{l1*} relies
on establishing first a weak-type $(1,1)$ inequality for the one dimensional
case and the special family of operators given by $T_{t}f(x)=\frac
{f(x+t)-f(x)}{t^{2}},$ more precisely,
\begin{equation}
\sup_{\lambda>0}\lambda\,\mathcal{L}^{2}\bigg(\bigg\{(x,t)\in\mathbb{R}%
\times(0,\infty):\frac{|f(x+t)-f(x)|}{t^{2}}>\lambda\bigg\}\bigg)\lesssim\Vert
f^{\prime}\Vert_{L^{1}(\mathbb{R})}. \label{BSYLim}%
\end{equation}
The extension to $\mathbb{R}^{N}\times(0,\infty)$ is then achieved by the
method of rotations. In Section \ref{sec:p} we consider a version of
(\ref{BSYLim}) in the setting of metric doubling measure spaces which thus
avoids the use of the method of rotations.

Let $(X,d,m)$ be a doubling metric measure space and consider the sequence of
integral averages $T_{t}f(x)=\frac{1}{m(B(x,t))}\int_{B(x,t)}f(y)\,\mathrm{d}
m(y).$ Using ideas of Carleson \cite{Carleson}, we extend Vitali coverings
from $X$ to $X\times(0,\infty)$ and show that (cf. Theorem \ref{ThmDav}
below): For $1\leq p<\infty,$
\[
(m\times\mathcal{L})\bigg(\bigg\{(x,t)\in X\times(0,\infty):\frac{|T_{t}%
f(x)|}{t^{1/p}}>\lambda\bigg\}\bigg)\leq C\frac{\Vert f\Vert_{L^{p}(X,m )}%
^{p}}{\lambda^{p}} \quad\text{for all} \quad\lambda> 0.
\]
We refer to Section \ref{sec:p} for more on the applications to Sobolev inequalities.

In Section \ref{sec:Ca-Ca-section} we study embeddings of
Calder\'{o}n--Campanato spaces (cf. \cite{Calderon72}, \cite{CalderonScott},
\cite{DeVore}, \cite{see}) defined by%
\begin{equation}
\left\Vert f\right\Vert _{C_{p}^{s}({\mathbb{R}}^{N})}:=\left\Vert f_{s}%
^{\#}\right\Vert _{L^{p}({\mathbb{R}}^{N})},\quad1<p<\infty,\quad s\in
\lbrack0,1], \label{DefCC}%
\end{equation}
where $f_{s}^{\#}$ is the sharp fractional maximal operator given by
\begin{equation}
f_{s}^{\#}(x)=\sup_{r>0}\frac{1}{r^{s+N}}\int_{B(x,r)}|f(y)-(f)_{B(x,r)}%
|\,\mathrm{d}y,\quad x\in{\mathbb{R}}^{N}, \label{DefCC2}%
\end{equation}
and $(f)_{B(x,r)}=\frac{1}{\mathcal{L}^{N}(B(x,r))}\int_{B(x,r)}%
f(y)\,\mathrm{d}y.$ Note that, when $s=0,$ we recover the Fefferman--Stein
maximal operator \footnote{It is well known this operator plays a key role in
the theory of the space $\text{BMO}({\mathbb{R}}^{N})$ formed by functions of
bounded mean oscillation; cf. \cite{Fefferman}.} $f^{\#}$. We also remark that
Seeger \cite{see} has shown that%
\[
C_{p}^{s}({\mathbb{R}}^{N})=F_{p,\infty}^{s}({\mathbb{R}}^{N}).
\]

Using the improvement of Bojarski's inequality obtained by DeVore--Sharpley
(cf. \eqref{Calderon}), we show that\footnote{Recall that the space
$BSY_{p}^{s}({\mathbb{R}}^{N})$ was introduced in \eqref{pppp}.} (cf. Example
\ref{ExampleC})%

\begin{equation}
C_{p}^{s}({\mathbb{R}}^{N})\subset BSY_{p}^{s}({\mathbb{R}}^{N}),\quad
1<p<\infty,\quad s\in(0,1]. \label{zas}%
\end{equation}
The sharpness of (\ref{zas}) will be discussed in Section \ref{CC}.

To incorporate to our theory generalized Calder\'{o}n--Campanato spaces
defined on metric spaces we shall require new maximal inequalities. Our
results in this direction are presented in Section \ref{sec:max-ca-ca}, where
we prove new maximal inequalities that were inspired by the beautiful
works\footnote{We believe that the strategy of proof, which in particular
avoids the assumption of doubling conditions, is of independent interest.} of
Gagliardo \cite{Gagliardo61, Gagliardo63} and Garsia \cite{Garsia72}. We
believe that the methods and results presented in this section are of
independent interest. Here we will just mention a typical embedding result.

Suppose that the metric measure space $(X,d,m)$ is $N$-regular, and let
$\rho:(0,\infty)\rightarrow(0,\infty)$ be continuous and increasing, with
$\lim_{r \rightarrow0+}\rho(r)=0$. Then (cf. Theorem \ref{ThmCCBSY})
\[
\left\Vert \frac{ f(x)-f(y) }{d(x,y)^{\frac{N}{p}}\Big(\int_{0}^{2d(x,y)}%
\frac{\mathrm{d}\bar{\rho}(\lambda)}{\lambda^{2N}}\Big)}\right\Vert
_{L(p,\infty)(X\times X,m\times m)}\lesssim\Vert f\Vert_{C_{p}^{\bar{\rho}%
}(X,m)}, \quad1 < p < \infty,
\]
where $\bar{\rho}(r)=\rho(2 r),$ and the generalized Calder\'{o}n--Campanato
space $C_{p}^{\bar{\rho}}(X,m)$ is defined by%
\[
\Vert f\Vert_{C_{p}^{\bar{\rho}}(X,m)}:=\left\Vert f_{\bar{\rho}}^{\#}
\right\Vert _{L^{p}(X,m)}
\]
and
\[
f_{\bar{\rho}}^{\#}(x)=\sup_{0<r<\text{diam}(X)}\frac{m(B(x,r))}{\bar{\rho
}(r)}\int_{B(x,r)}|f(y)-(f)_{B(x,r)}|\,\mathrm{d}m(y)
\]
where $(f)_{B(x,r)}=\frac{1}{m(B(x,r)) }\int_{B(x,r)}f(y)\,\mathrm{d} m( y).$

In Section \ref{sec:Crippa} we consider the limiting Calder\'{o}n--Campanato
spaces $\mathcal{N}^{s,p}(\mathbb{R}^{N})$ endowed with
\begin{equation}
\label{qqq}\Vert f\Vert_{\mathcal{N}^{s,p}(\mathbb{R}^{N})}:=\left\Vert
\Phi^{s}(f)\right\Vert _{L^{p}(\mathbb{R}^{N})},\quad s\in(0,1], \quad1 < p <
\infty,
\end{equation}
where the maximal function $\Phi^{s}$ is defined by%
\begin{equation}
\label{qqq56}\Phi^{s}(f)(x)=\sup_{r>0} \frac{1}{\mathcal{L}^{N}(B(x,r))}%
\int_{B(x,r)}\log\bigg(\frac{|f(x)-f(y)|}{r^{s}}+1\bigg)\,\mathrm{d}y,\quad
x\in{\mathbb{R}}^{N}.
\end{equation}
These spaces play an important role in the study of Lagrangian flows to
Sobolev fields (cf. \cite{Crippa}, \cite{Brue} and \cite{BrueSemola}). We
compare the functionals \eqref{qqq} to logarithmic versions of the
Gagliardo--Brezis--Van Schaftingen--Yung functionals. In fact, using the
recent estimates of Crippa--De Lellis \cite{Crippa}, and mixed norm
inequalities as above, we show (cf. Theorem \ref{BrezisExp})%

\[
\bigg\|\log\bigg(1+\frac{|f(x)-f(y)|}{|x-y|^{s}}\bigg)\frac{1}{|x-y|^{\frac
{N}{p}}}\bigg\|_{L(p,\infty)(\mathbb{R}^{N}\times\mathbb{R}^{N})}\lesssim\Vert
f\Vert_{\mathcal{N}^{s,p}(\mathbb{R}^{N})}.
\]
Conversely (cf. Proposition \ref{ThmBN}), when $s=1$ we have that, for $f\in
C^{1}({\mathbb{R}}^{N})$,%
\begin{align*}
\int_{{\mathbb{R}}^{N}}\log{\small \,(1+|\nabla f(x)|)}^{p}%
{\small \,\mathrm{d}x}  & \\
&  \hspace{-3cm}\lesssim\liminf_{\lambda\rightarrow\infty}\lambda
^{p}\mathcal{L}^{2N}\bigg (\bigg\{(x,y)\in{\mathbb{R}}^{N}\times{\mathbb{R}%
}^{N}:\log\bigg(1+\frac{|f(x)-f(y)|}{|x-y|}\bigg)\frac{1}{|x-y|^{\frac{N}{p}}%
}>\lambda\bigg\}\bigg).
\end{align*}
In particular, this result sharpens the recent inequality due to
Bru\'{e}--Nguyen \cite[Proposition 2.6]{Brue},
\[
\int_{{\mathbb{R}}^{N}}\log(1+|\nabla f(x)|)^{p}\,\mathrm{d}x\lesssim\Vert
f\Vert_{\mathcal{N}^{1,p}({\mathbb{R}}^{N})}^{p}.
\]

The paper is naturally divided in three parts: The first part, comprising
Sections \ref{sec:maximal} and \ref{SubsectionAppl}, deals with Theorem
\ref{ThmMaximal} and Applications to Analysis, the second part (Sections
\ref{sec:Ca-Ca-section}, \ref{sec:max-ca-ca}, \ref{sec:Crippa}) devoted to
Calder\'{o}n--Campanato spaces, and the third part (Sections
\ref{sec:Caff-Silv}, \ref{sec:further}) that concerns with generalized Riesz
potential type spaces, extension theorems and applications. Finally in the
Appendix we show how to incorporate the Bourgain--Nguyen characterizations of
Sobolev spaces to our framework, by means of a suitable variant of Theorem
\ref{ThmMaximal} (cf. Theorem \ref{ThmMaximal2}).

\vspace{2mm}
\textbf{Acknowledgements}. \textit{We are very grateful to the referees for
providing precious information, references and many interesting questions that
helped us improve the quality of the paper. In particular, the Appendix was
developed to address their questions.}

\textit{We are also grateful to our colleagues Michael Cwikel, Carlos P\'erez and Andreas Seeger for useful comments and references to the literature.}

\section{Weak-type estimates related to sequences of operators: Proof of
Theorem \ref{ThmMaximal}\label{sec:maximal}}

Before we proceed with the proof of Theorem \ref{ThmMaximal}, we highlight an
immediate consequence of the result that enables us to characterize Lebesgue
norms, in terms of the distribution function of both $(x,t)$ variables.
Namely, we have the following

\begin{corollary}
\label{CorMaximal} Let $(X,m)$ be a $\sigma$-finite measure space, and let
$\{T_{t}:t>0\}$ be a one-parameter family of sublinear operators on $L^{p}%
(X,m),\,1<p<\infty$, furthermore suppose that $\gamma>0.$ Assume that there
exists a positive constant $C_{p}$ such that
\begin{equation}
\Vert T^{\ast}f\Vert_{L^{p}(X,m)}\leq C_{p}\,\Vert f\Vert_{L^{p}(X,m)}%
\quad\text{for all}\quad f\in L^{p}(X,m) \label{MaxAssump2}%
\end{equation}
and
\begin{equation}
\lim_{t\rightarrow0+}T_{t}f(x)=f(x)\quad\text{for a dense subclass of}\quad
L^{p}(X,m). \label{MaxAssump3}%
\end{equation}
Then
\begin{equation}
\frac{1}{\gamma}\,\Vert f\Vert_{L^{p}(X,m)}^{p}\leq\lim_{\lambda
\rightarrow\infty}\lambda^{p}\iint_{E_{\lambda,\gamma/p}}t^{\gamma
-1}\,\mathrm{d}t\,\mathrm{d}m(x)\leq\sup_{\lambda>0}\lambda^{p}\iint%
_{E_{\lambda,\gamma/p}}t^{\gamma-1}\,\mathrm{d}t\,\mathrm{d}m(x)\leq
\frac{C_{p}^{p}}{\gamma}\,\Vert f\Vert_{L^{p}(X,m)}^{p} \label{weight}%
\end{equation}
where
\[
E_{\lambda,\gamma/p}=\{(x,t)\in X\times(0,\infty):|T_{t}f(x)|>\lambda
\,t^{\gamma/p}\}.
\]
In particular, if $\gamma=1$ then
\begin{align*}
\Vert f\Vert_{L^{p}(X,m)}^{p}  &  \leq\lim_{\lambda\rightarrow\infty}%
\lambda^{p}(m\times\mathcal{L})\bigg(\bigg\{(x,t)\in X\times(0,\infty
):\frac{|T_{t}f(x)|}{t^{1/p}}>\lambda\bigg\}\bigg)\\
&  \leq\sup_{\lambda>0}\lambda^{p}(m\times\mathcal{L})\bigg(\bigg\{(x,t)\in
X\times(0,\infty):\frac{|T_{t}f(x)|}{t^{1/p}}>\lambda\bigg\}\bigg)\leq
C_{p}^{p}\,\Vert f\Vert_{L^{p}(X,m)}^{p}.
\end{align*}

\end{corollary}

The corollary follows readily from Theorem \ref{ThmMaximal} and the well-known
fact that under the assumptions \eqref{MaxAssump2} and \eqref{MaxAssump3} one
has that, for each $f\in L^{p}(X,m)$, $\lim_{t\rightarrow0+}T_{t}f(x)=f(x)$
$m$-a.e. $x\in X$.

\begin{remark}\label{RemarkSharpAssertion}
	Under assumptions of Corollary \ref{CorMaximal}, the first inequality in \eqref{weight} becomes in fact an equality, i.e.,
	$$
		\frac{1}{\gamma}\,\Vert f\Vert_{L^{p}(X,m)}^{p} = \lim_{\lambda
\rightarrow\infty}\lambda^{p}\iint_{E_{\lambda,\gamma/p}}t^{\gamma
-1}\,\mathrm{d}t\,\mathrm{d}m(x).
	$$
	This exact formula is a special case of a more general phenomenon explained in Section \ref{SectionBN} below (cf. \eqref{Statement1.2}).   
\end{remark}


\begin{proof}
[Proof of Theorem \ref{ThmMaximal}]The first part of the theorem was proved in
the Introduction. It remains to prove (ii). \ According to \eqref{MaxAssump1}
there exists a measurable set $N\subset X$ such that $m(X\backslash N)=0$ and
\[
g(x):=\lim_{t\rightarrow0+}T_{t}f(x)\quad\text{for all}\quad x\in N.
\]
Let $x\in N$ be fixed such that $g(x)\neq0$. Given any $n\in{{\mathbb{N}}%
},\,n>1,$ there exists $t_{0}=t_{0}(x,n)>0$ for which
\[
\Big(1-\frac{1}{n}\Big)|g(x)|\leq|T_{t}f(x)|\quad\text{for all}\quad
t\in(0,t_{0}).
\]
Thus, for each $\lambda>0$ we have,
\[
\bigg\{t\in(0,t_{0}):\Big(1-\frac{1}{n}\Big)\frac{|g(x)|}{t^{\gamma/p}%
}>\lambda\bigg\}\subseteq\bigg\{t\in(0,\infty):\frac{|T_{t}f(x)|}{t^{\gamma
/p}}>\lambda\bigg\}.
\]
This implies
\[
\bigg\{t\in(0,\infty):t<\min\bigg\{\bigg(\frac{|g(x)|}{\lambda}\Big(1-\frac
{1}{n}\Big)\bigg)^{\frac{p}{\gamma}},t_{0}(x,n)\bigg\}\bigg\}\subseteq
\bigg\{t\in(0,\infty):\frac{|T_{t}f(x)|}{t^{\gamma/p}}>\lambda\bigg\}
\]
and consequently
\[
w_{\gamma}\bigg(\bigg(0,\min\bigg\{\bigg(\frac{|g(x)|}{\lambda}\Big(1-\frac
{1}{n}\Big)\bigg)^{\frac{p}{\gamma}},t_{0}(x,n)\bigg\}\bigg)\bigg)\leq
w_{\gamma}\bigg(\bigg\{t\in(0,\infty):\frac{|T_{t}f(x)|}{t^{\gamma/p}}%
>\lambda\bigg\}\bigg).
\]
The previous inequality can be rewritten as
\[
\frac{1}{\gamma}\,\min\bigg\{\bigg(\frac{|g(x)|}{\lambda}\Big(1-\frac{1}%
{n}\Big)\bigg)^{p},(t_{0}(x,n))^{\gamma}\bigg\}\leq w_{\gamma}%
\bigg(\bigg\{t\in(0,\infty):\frac{|T_{t}f(x)|}{t^{\gamma/p}}>\lambda
\bigg\}\bigg).
\]
Integrating both sides we arrive at,
\begin{align*}
\frac{1}{\gamma}\,\int_{\{x\in N:g(x)\neq0\}}\min\bigg\{\bigg(\frac
{|g(x)|}{\lambda}\Big(1-\frac{1}{n}\Big)\bigg)^{p},(t_{0}(x,n))^{\gamma
}\bigg\}\,\mathrm{d}m(x)  &  \leq\\
&  \hspace{-7cm}\int_{\{x\in N:g(x)\neq0\}}w_{\gamma}\bigg(\bigg\{t\in
(0,\infty):\frac{|T_{t}f(x)|}{t^{\gamma/p}}>\lambda\bigg\}\bigg)\,\mathrm{d}%
m(x)\\
&  \hspace{-7cm}\leq(m\times w_{\gamma})\bigg(\bigg\{(x,t)\in X\times
(0,\infty):\frac{|T_{t}f(x)|}{t^{\gamma/p}}>\lambda\bigg\}\bigg).
\end{align*}
Since $\lambda$ is independent of $x$ we can factor it out $\lambda^{p}$, and
obtain
\begin{align*}
\frac{1}{\gamma}\,\int_{\{x\in N:g(x)\neq0\}}\min
\bigg\{\bigg(|g(x)|\Big(1-\frac{1}{n}\Big)\bigg)^{p},\lambda^{p}%
(t_{0}(x,n))^{\gamma}\bigg\}\,\mathrm{d}m(x)  &  \leq\\
&  \hspace{-7cm}\lambda^{p}(m\times w_{\gamma})\bigg(\bigg\{(x,t)\in
X\times(0,\infty):\frac{|T_{t}f(x)|}{t^{\gamma/p}}>\lambda\bigg\}\bigg).
\end{align*}
We now let $\lambda\rightarrow\infty,$ then by Fatou's Lemma, (note that $n$
is independent of $x$), we see that
\begin{align*}
\frac{1}{\gamma}\,\Big(1-\frac{1}{n}\Big)^{p}\int_{\{x\in N:g(x)\neq
0\}}|g(x)|^{p}\,\mathrm{d}m(x)  &  \leq\\
&  \hspace{-5cm}\frac{1}{\gamma}\,\liminf_{\lambda\rightarrow\infty}%
\int_{\{x\in N:g(x)\neq0\}}\min\bigg\{\bigg(|g(x)|\Big(1-\frac{1}%
{n}\Big)\bigg)^{p},\lambda^{p}(t_{0}(x,n))^{\gamma}\bigg\}\,\mathrm{d}m(x)\\
&  \hspace{-5cm}\leq\liminf_{\lambda\rightarrow\infty}\lambda^{p}(m\times
w_{\gamma})\bigg(\bigg\{(x,t)\in X\times(0,\infty):\frac{|T_{t}f(x)|}%
{t^{\gamma/p}}>\lambda\bigg\}\bigg).
\end{align*}
Finally, we let $n\rightarrow\infty$ to obtain
\begin{align*}
\frac{1}{\gamma}\,\int_{\{x\in N:g(x)\neq0\}}|g(x)|^{p}\,\mathrm{d}m(x)  &
\leq\\
&  \hspace{-3cm}\liminf_{\lambda\rightarrow\infty}\lambda^{p}(m\times
w_{\gamma})\bigg(\bigg\{(x,t)\in X\times(0,\infty):\frac{|T_{t}f(x)|}%
{t^{\gamma/p}}>\lambda\bigg\}\bigg).
\end{align*}
This concludes the proof of (ii).
\end{proof}

\begin{remark}
\label{RemMulti} The unweighted multiparameter version of Theorem
\ref{ThmMaximal} reads as follows: Let $(X,m)$ be a $\sigma$-finite measure
space, and let $\{T_{t_{1},\ldots,t_{l}}\}$ be a multiparameter family of (not
necessarily linear) operators on $L^{p}(X,m), \, 1 \leq p < \infty$.

\begin{enumerate}
[\upshape(i)]

\item Assume
\[
\sup_{t_{1}>0,\ldots,t_{l}>0}|T_{t_{1},\ldots,t_{l}}f|\in L^{p}(X,m).
\]
Then
\begin{align*}
\sup_{\lambda>0}\lambda^{p}(m\times\mathcal{L}^{l})\bigg(\bigg\{(x,t_{1}%
,\ldots,t_{l})  &  \in X\times(0,\infty)^{l}:\frac{|T_{t_{1},\ldots,t_{l}%
}f(x)|}{\max\{t_{1},\ldots,t_{l}\}^{l/p}}>\lambda\bigg\}\bigg)\leq\\
&  \Big\|\sup_{t_{1}>0,\ldots,t_{l}>0}|T_{t_{1},\ldots,t_{l}}f|\Big\|_{L^{p}%
(X,m)}^{p}.
\end{align*}

\item Suppose that
\[
\lim_{t_{1}\rightarrow0+,\ldots,t_{l}\rightarrow0+}T_{t_{1},\ldots,t_{l}%
}f(x)<\infty\quad\text{$m$-a.e.}\quad x\in X.
\]
Then
\begin{align*}
\Big\|\lim_{t_{1}\rightarrow0+,\ldots,t_{l}\rightarrow0+}T_{t_{1},\ldots
,t_{l}}f\Big\|_{L^{p}(X,m)}^{p}  &  \leq\\
&  \hspace{-5cm}\liminf_{\lambda\rightarrow\infty}\lambda^{p}(m\times
\mathcal{L}^{l})\bigg(\bigg\{(x,t_{1},\ldots,t_{l})\in X\times(0,\infty
)^{l}:\frac{|T_{t_{1},\ldots,t_{l}}f(x)|}{\max\{t_{1},\ldots,t_{l}\}^{l/p}%
}>\lambda\bigg\}\bigg).
\end{align*}

\end{enumerate}

The weighted version of the previous result can also be obtained, but we omit
here further details.
\end{remark}

\subsection{A limiting case}

Let $\omega$ \ be a weight on ${\mathbb{R}}^{N}$, it is plain that
\[
\Vert f\Vert_{L^{p}({\mathbb{R}}^{N})}=\bigg\|\frac{f}{\omega^{1/p}%
}\bigg\|_{L^{p}({\mathbb{R}}^{N},\omega\,\mathrm{d}x)}%
\]
in other words
\[
f\in L^{p}({\mathbb{R}}^{N})\quad\text{if and only if}\quad\frac{f}%
{\omega^{1/p}}\in L^{p}({\mathbb{R}}^{N},\omega\,\mathrm{d}x).
\]
By Theorem \ref{ThmMaximal} (and Corollary \ref{CorMaximal}) a similar result
holds for families of operators, and the family of weights $\{t^{\gamma
}\}_{\gamma>0}$. As it is illustrated by \eqref{weight}, the statement does
not hold for the limiting value $\gamma=0$. In this section we consider this
case and prove the following.

\begin{theorem}
\label{ThmMaximalNew} Let $(X,m)$ be a $\sigma$-finite measure space, and let
$\{T_{t}:t>0\}$ be a one-parameter family of (not necessarily linear)
operators on $L^{p}(X,m), \, 1 \leq p < \infty$. We let $\eta>1$, and for
Lebesgue measurable sets $A\subset(0,1/2)$ define
\begin{equation}
v_{\eta}(A)=\int_{A}t^{-1}\log^{-\eta}\Big(\frac{1}{t}\Big)\,\mathrm{d}t.
\label{Weight}%
\end{equation}

\begin{enumerate}
[\upshape(i)]

\item Assume
\[
T^{\ast}f \in L^{p}(X,m).
\]
Then
\[
\sup_{\lambda>0}\lambda^{p}(m\times v_{\eta})\bigg(\bigg\{(x,t)\in
X\times(0,1/2):\frac{|T_{t}f(x)|}{\log^{\frac{-\eta+1}{p}}\Big(\frac{1}%
{t}\Big)}>\lambda\bigg\}\bigg)\leq\frac{1}{\eta-1}\,\Vert T^{\ast}%
f\Vert_{L^{p}(X,m)}^{p}.
\]
In particular, if $p > 1$ and the maximal operator $T^{\ast}$ is bounded on
$L^{p}(X,m)$:
\[
\left\Vert T^{\ast}f\right\Vert _{L^{p}(X,m)}\leq C_{p}\left\Vert f\right\Vert
_{L^{p}(X,m)}%
\]
then
\[
\sup_{\lambda>0}\lambda^{p}(m\times v_{\eta})\bigg(\bigg\{(x,t) \in
X\times(0,1/2):\frac{|T_{t}f(x)|}{\log^{\frac{-\eta+1}{p}}\Big(\frac{1}%
{t}\Big)}>\lambda\bigg\}\bigg)\leq\frac{C_{p}^{p}}{\eta- 1} \Vert
f\Vert_{L^{p}(X,m)}^{p}.
\]

\item Suppose that
\[
\lim_{t\rightarrow0+}T_{t}f(x)<\infty\quad\text{$m$-a.e.}\quad x\in X.
\]
Then
\begin{align*}
&  \frac{1}{\eta-1}\,\Big\|\lim_{t\rightarrow0+}T_{t}f\Big\|_{L^{p}(X,m)}%
^{p}\\
&  \leq\liminf_{\lambda\rightarrow\infty}\lambda^{p}(m\times v_{\eta
})\bigg(\bigg\{(x,t)\in X\times(0,1/2):\frac{|T_{t}f(x)|}{\log^{\frac{-\eta
+1}{p}}\Big(\frac{1}{t}\Big)}>\lambda\bigg\}\bigg).
\end{align*}

\end{enumerate}
\end{theorem}

The proof of this result is similar to the proof of Theorem \ref{ThmMaximal}
and we shall leave it to the reader. When comparing Theorems \ref{ThmMaximal}
and \ref{ThmMaximalNew}, note that in Theorem \ref{ThmMaximalNew} we switch
from weights of logarithmic order $\log^{-\eta}\Big(\frac{1}{t}\Big)$ (see
\eqref{Weight}) to denominators $\log^{\frac{-\eta+1}{p}}\Big(\frac{1}%
{t}\Big)$ in the Gagliardo-type quotients. The explanation of this phenomenon
hinges on elementary computations: for any $T\in(0,1/2)$
\[
w_{\gamma}([0,T])\asymp T^{\gamma}\quad\text{and}\quad v_{\eta}([0,T])\asymp
\log^{-\eta+1}\Big(\frac{1}{t}\Big).
\]
These special choices of weights seem to indicate that the denominator in the
Gagliardo-type quotients related the measure on $(0,\infty)$ induced by a
general weight $\omega$ is $\int_{0}^{t} \omega(u) \, \mathrm{d} u$. Indeed,
this conjecture will be confirmed in Remark \ref{RemWeights} below.

As an immediate consequence we have the following

\begin{corollary}
\label{CorMaximalNew}Let $(X,m)$ be a $\sigma$-finite measure space, and let
$\{T_{t}:t>0\}$ be a one-parameter family of sublinear operators on $L^{p}%
(X,m),\,1<p<\infty$, furthermore suppose that $\eta>1$. Assume that there
exists a positive constant $C_{p}$ such that
\[
\Vert T^{\ast}f\Vert_{L^{p}(X,m)}\leq C_{p}\,\Vert f\Vert_{L^{p}(X,m)}%
\quad\text{for all}\quad f\in L^{p}(X,m)
\]
and
\[
\lim_{t\rightarrow0+}T_{t}f(x)=f(x)\quad\text{for a dense subclass of}\quad
L^{p}(X,m).
\]
Then
\begin{align*}
\frac{1}{\eta-1}\,\Vert f\Vert_{L^{p}(X,m)}^{p}  &  \leq\lim_{\lambda
\rightarrow\infty}\lambda^{p}\iint_{\mathcal{E}_{\lambda,\frac{-\eta+1}{p}}%
}t^{-1}\log^{-\eta}\Big(\frac{1}{t}\Big)\,\mathrm{d}t\,\mathrm{d}m(x)\\
&  \leq\sup_{\lambda>0}\lambda^{p}\iint_{\mathcal{E}_{\lambda,\frac{-\eta
+1}{p}}}t^{-1}\log^{-\eta}\Big(\frac{1}{t}\Big)\,\mathrm{d}t\,\mathrm{d}m(x)\\
&  \leq\frac{C_{p}^{p}}{\eta-1}\,\Vert f\Vert_{L^{p}(X,m)}^{p},
\end{align*}
where
\[
\mathcal{E}_{\lambda,\frac{-\eta+1}{p}}=\Big\{(x,t)\in X\times(0,1/2):|T_{t}%
f(x)|>\lambda\,\log^{\frac{-\eta+1}{p}}\Big(\frac{1}{t}\Big)\Big\}.
\]

\end{corollary}

With Theorem \ref{ThmMaximalNew} and Corollary \ref{CorMaximalNew} at hand, we
are able to state and prove the corresponding limiting versions for all the
applications that are given in Section \ref{SubsectionAppl} below. Here we
shall only write down the limiting case of the Brezis--Van Schaftingen--Yung
formula for first order Sobolev seminorms given in Example \ref{ExampleBSY}:
If $\eta> 1$ and $1 < p < \infty$ then
\begin{align*}
\lim_{\lambda \to \infty} \iint_{\mathcal{E}%
_{\lambda,\frac{-\eta+1}{p}}}|x-y|^{-1}\log^{-\eta}\Big(\frac{1}{|x-y|}\Big)\,\mathrm{d}x\,\mathrm{d}y &\asymp \Vert \nabla f\Vert_{L^{p}%
({\mathbb{R}^N})}^{p} \\
& \hspace{-5cm}\asymp \sup_{\lambda> 0}\lambda^{p}\iint_{\mathcal{E}%
_{\lambda,\frac{-\eta+1}{p}}}|x-y|^{-1}\log^{-\eta}\Big(\frac{1}{|x-y|}\Big)\,\mathrm{d}x\,\mathrm{d}y
\end{align*}
where
\[
\mathcal{E}_{\lambda,\frac{-\eta+1}{p}}=\bigg\{(x,y)\in{\mathbb{R}^N}%
\times \mathbb{R}^N: |x-y| < \frac{1}{2}, \quad |f(x)-f(y)|>\lambda\,|x-y| \, \log^{\frac{-\eta+1}{p}}\Big(\frac{1}%
{|x-y|}\Big)\bigg\}.
\]

\subsection{The case $p=1$\label{sec:p}}

The proof of Theorem \ref{ThmMaximal} works for $p=1$ without any changes, if
we assume that $T^{\ast}$ is bounded $L^{1}(X,m);$ however this assumption is
too strong, and does not hold in the main applications. In this section we
offer a formulation of Corollary \ref{CorMaximal} that avoids this assumption
and, furthermore, provides an exact formula.

\begin{theorem}
\label{ThmPh}Let $(X,m)$ be a $\sigma$-finite measure space, and let
$\{T_{t}:t>0\}$ be a one-parameter family of operators on $L^{p}(X,m),1\leq
p<\infty$. Then
\[
\left\Vert f\right\Vert _{L^{p}(X,m)}^{p}=\lim_{\lambda\rightarrow\infty
}\lambda^{p}(m\times\mathcal{L})\bigg(\bigg\{(x,t)\in X\times(0,\infty
):\frac{|T_{t}f(x)|}{t^{1/p}}>\lambda\bigg\}\bigg)
\]
for all $f\in L^{p}(X,m)$ satisfying
\begin{equation}
\left\Vert T_{t}f-f\right\Vert _{L^{\infty}(X,m)}\leq C_{f}\,t^{1/p},\quad
t>0, \label{formerly25}%
\end{equation}
for some $C_{f}>0$ which is independent of $t$ (but depends on $f$).
\end{theorem}

\begin{proof}
Using \eqref{formerly25} and the triangle inequality we find
\begin{equation}
|T_{t}f(x)|\leq C_{f}\,t^{1/p}+|f(x)|. \label{ProofLimit1}%
\end{equation}

Given $\lambda>0$, let
\[
E(f,\lambda)=\bigg\{(x,t)\in X\times(0,\infty):\frac{|T_{t}f(x)|}{t^{1/p}%
}>\lambda\bigg\}.
\]
Since we will be only interested in large values of $\lambda$ we may assume,
without loss of generality, that $\lambda>C_{f}$. Let $(x,t)\in E(f,\lambda)$,
then combining with \eqref{ProofLimit1} and rearranging yields,
\[
t\leq\bigg(\frac{|f(x)|}{\lambda-C_{f}}\bigg)^{p}.
\]
Consequently,
\[
E(f,\lambda)\subseteq\bigg\{(x,t):f(x)\neq0,\,t\leq\bigg(\frac{|f(x)|}%
{\lambda-C_{f}}\bigg)^{p}\bigg\}.
\]
It follows that%
\begin{align*}
\limsup_{\lambda\rightarrow\infty}\lambda^{p}(m\times\mathcal{L}%
)(E(f,\lambda))  &  \leq\limsup_{\lambda\rightarrow\infty}\lambda^{p}%
\int_{\{x\in X:f(x)\neq0\}}\bigg(\frac{|f(x)|}{\lambda-C_{f}}\bigg)^{p}%
\,\mathrm{d}m(x)\\
&  =\Big(\limsup_{\lambda\rightarrow\infty}\frac{\lambda}{\lambda-C_{f}%
}\Big)^{p}\,\int_{X}|f(x)|^{p}\,\mathrm{d}m(x)=\int_{X}|f(x)|^{p}%
\,\mathrm{d}m(x).
\end{align*}

The converse inequality, i.e.,
\[
\int_{X}|f(x)|^{p}\,\mathrm{d} m(x)=\int_{X}|f(x)|^{p}\,\mathrm{d}
m(x)\leq\liminf_{\lambda\rightarrow\infty}\lambda^{p}(m\times\mathcal{L}%
)(E(f,\lambda))
\]
follows from our previous results since the condition \eqref{formerly25}
implies that
\[
\lim_{t\rightarrow0+}T_{t}f(x)=f(x),\quad x\in X.
\]

\end{proof}

\begin{example}
\label{exbsyp=1}Let $1\leq p<\infty.$ For locally integrable functions
consider the the family of operators%
\[
T_{t}f(x)=\frac{1}{m(B(x,t))}\int_{B(x,t)}f(y) \, \mathrm{d} m(y),\quad t>0,
\quad x \in X.
\]
Let us consider the H\"older--Lipschitz class
\begin{equation}
\label{HLNew}\mathcal{C}^{1/p}(X)=\{f:X\rightarrow{\mathbb{R}}:\exists M>0
\quad\text{s.t.} \quad\left\vert f(x)-f(y)\right\vert \leq M\left(
d(x,y)\right)  ^{1/p} \quad\text{for all} \quad x,y\in X\}
\end{equation}
endowed with the seminorm
\begin{equation}
\label{HLNew2}\|f\|_{\mathcal{C}^{1/p}(X)} = \inf M.
\end{equation}
Assume $f \in\mathcal{C}^{1/p}(X)$. Then it is readily verified that%
\[
\left\Vert T_{t}f-f\right\Vert _{L^{\infty}(X,m)}\leq\|f\|_{\mathcal{C}%
^{1/p}(X)} \, t^{1/p}%
\]
and consequently by Theorem \ref{ThmPh}%
\begin{equation}
\left\Vert f\right\Vert _{L^{p}(X,m)}^{p}=\lim_{\lambda\rightarrow\infty
}\lambda^{p}(m\times\mathcal{L})\bigg( \bigg\{(x,t)\in X\times(0,\infty
):\frac{ \frac{1}{m(B(x,t))} \left\vert \int_{B(x,t)} f(y) \, \mathrm{d}
m(y)\right\vert }{t^{1/p}}>\lambda\bigg\}\bigg). \label{dav}%
\end{equation}
Note that the previous reasoning works with arbitrary $\sigma$-finite measures
$m$ on $X$, not necessarily doubling measures.
\end{example}

\begin{example}
\label{exbysp=2}Let $X={\mathbb{R}}$ and $1\leq p<\infty$. For $f\in
\mathcal{S}({\mathbb{R}}),$ let%
\[
T_{t}f(x)=\frac{1}{t}\int_{x}^{x+t}f(y)\,\mathrm{d}y,\quad x\in{\mathbb{R}%
},\quad t>0,
\]
then (\ref{formerly25}) is verified with $C_{f}=\left\Vert f\right\Vert
_{\mathcal{C}^{1/p}({\mathbb{R}})}$ (cf. \eqref{HLNew}, \eqref{HLNew2}). By
assumption, $f^{\prime}$ also satisfies the $\mathcal{C}^{1/p}$ condition and
therefore (\ref{dav}) implies that
\begin{align*}
\left\Vert f^{\prime}\right\Vert _{L^{p}({\mathbb{R}})}^{p}  &  =\lim
_{\lambda\rightarrow\infty}\lambda^{p}\mathcal{L}^{2}\bigg(\bigg\{(x,t)\in
{\mathbb{R}}\times(0,\infty):\frac{\left\vert \frac{1}{t}\int_{x}%
^{x+t}f^{\prime}(y)\,\mathrm{d}y\right\vert }{t^{1/p}}>\lambda\bigg\}\bigg)\\
&  =\lim_{\lambda\rightarrow\infty}\lambda^{p}\mathcal{L}^{2}%
\bigg(\bigg\{(x,t)\in{\mathbb{R}}\times(0,\infty):\frac{\left\vert
f(x+t)-f(x)\right\vert }{t^{1+1/p}}>\lambda\bigg\}\bigg).
\end{align*}
This recovers \eqref{I1'} with $N=1$.
\end{example}

\begin{remark}
As pointed out by the referees, in the context of the Bourgain--Nguyen
functionals and Sobolev spaces, \cite{ngu06} and \cite{bn06} rely on weaker convergence conditions than the uniform convergence \eqref{formerly25}. We feel that the convergence
issues deserve a separate treatment and we hope to deal with these issues elsewhere.
\end{remark}

A natural question here is: under what conditions one can replace
\textquotedblleft$\lim_{\lambda\rightarrow\infty}$\textquotedblright\ for by
\textquotedblleft$\sup_{\lambda>0}$\textquotedblright\ in \eqref{dav}? This
was already investigated in previous sections for $p>1,$ as a by-product of
the $L^{p}$ boundedness of the maximal operator $T^{\ast}f=\sup_{t>0}|T_{t}%
f|$. However, our method does not provide a satisfactory answer when $p=1,$
due to the lack of strong type estimates for $T^{\ast}$. As we pointed out in
the Introduction, in the special case $T_{t}f(x)=\frac{f(x+t)-f(x)}{t}%
,\,x\in\mathbb{R}$, \cite{Brezis} provides a proof of (\ref{BSYLim}) using a
Vitali covering for the product space $\mathbb{R}\times\mathbb{R}$. The
extension to $\mathbb{R}^{N}\times{\mathbb{R}}^{N}$ is achieved from the case
$N=1$ via the method of rotations. However, it is not clear how to extend
these techniques to the more general setting of metric spaces.

Our aim now is to provide an alternative methodology for \eqref{BSYLim}. Our
method is very close to the method of \cite{Brezis} but it has the advantage
that it works in metric spaces. Before we state our result, let us outline
what is the basic idea behind the proof: Instead of applying the Vitali
covering theorem on the product space $\mathbb{R}\times\mathbb{R}$ as in
\cite{Brezis}, we construct a suitable Vitali covering on $X$ and then extend
it to $X\times(0,\infty)$ using Carleson boxes. In this fashion we are able to
apply Fubini on the corresponding covering for $X\times(0,\infty)$.

\begin{theorem}
\label{ThmDav}Let $(X,d)$ be a metric space, with doubling measure $m$ and let
$1\leq p<\infty.$ Then, there exists a positive constant $C$ such that for all
$\lambda>0$%
\[
(m\times\mathcal{L}) \bigg(\bigg\{(x,t)\in X\times(0,\infty):\frac{\left\vert
T_{t}f(x)\right\vert }{t^{1/p}}>\lambda\bigg\} \bigg)\leq C \, \frac
{\left\Vert f\right\Vert _{L^{p}(X,m)}^{p}}{\lambda^{p}},
\]
where $T_{t}f(x)=\frac{1}{m(B(x,t))}\int_{B(x,t)}f(y) \, \mathrm{d} m(y).$
\end{theorem}

\begin{proof}
For each $\lambda>0$, we let
\[
E_{\lambda}=\bigg\{(x,t)\in X\times(0,\infty):\frac{|T_{t}f(x)|}{t^{1/p}%
}>\lambda\bigg\}
\]
and
\[
E_{\lambda}^{\prime}=\bigg\{x\in X:\sup_{t>0}\frac{|T_{t}f(x)|}{t^{1/p}%
}>\lambda\bigg\}.
\]
Given $x\in E_{\lambda}^{\prime}$ we introduce
\[
t_{x}=\sup\bigg\{t>0:\frac{|T_{t}f(x)|}{t^{1/p}}>\lambda\bigg\}.
\]
Then $t_{x}>0,$ and
\begin{equation}
\frac{|T_{t_{x}}f(x)|}{t_{x}^{1/p}}\geq\lambda. \label{Sup3}%
\end{equation}
Without loss of generality we can assume that $t_{x}$ is everywhere finite.
Suppose first that $E_{\lambda}^{\prime}$ is bounded. By Vitali's covering
theorem there exists a sequence of disjoint balls $B(x_{i},t_{x_{i}}%
),\,x_{i}\in E_{\lambda}^{\prime}$, such that
\begin{equation}
E_{\lambda}^{\prime}\subset\bigcup_{i}B(x_{i},4\,t_{x_{i}}). \label{Sup1}%
\end{equation}
This implies that
\begin{equation}
E_{\lambda}\subset\bigcup_{i}B(x_{i},4\,t_{x_{i}})\times(0,2\,t_{x_{i}}).
\label{Sup2}%
\end{equation}
Indeed, if $(x,t)\in E_{\lambda}$ then $x\in E_{\lambda}^{\prime}$ and by
\eqref{Sup1} there exists $i$ so that $x\in B(x_{i},t_{x_{i}})$. Furthermore,
following the proof of the Vitali covering theorem given in \cite{Coif}, we
may assume that $t_{x}<2\,t_{x_{i}}$. Therefore, since $t\leq t_{x}$, we infer
that $t<2\,t_{x_{i}}$.

Let $c_{m}$ be a constant depending only on the doubling property of $m.$
Using \eqref{Sup2}, \eqref{Sup3}, and Jensen's inequality, we get
\begin{align*}
(m\times\mathcal{L})(E_{\lambda})  &  \leq2\sum_{i}t_{x_{i}}m (B(x_{i}%
,4\,t_{x_{i}}))\leq\frac{2c_{m}}{\lambda^{p}}\sum_{i}|T_{t_{x_{i}}}%
f(x_{i})|^{p}m(B(x_{i},t_{x_{i}}))\\
&  \leq\frac{2c_{m}}{\lambda^{p}}\sum_{i}\int_{B(x_{i},t_{x_{i}})}%
|f(y)|^{p}\,\mathrm{d} m(y)\\
&  \leq\frac{2c_{m}}{\lambda^{p}}\Vert f\Vert_{L^{p}(X,m)}^{p}%
\end{align*}
where the last step follows from the disjointness of the collection of balls
$\{B(x_{i},t_{x_{i}})\}$.

To remove the assumption that $E_{\lambda}^{\prime}$ is bounded we apply the
argument given above to $E_{\lambda}^{\prime}\cap B(a,R)$ with $a\in X$ and
$R>0$. This yields,
\[
(m\times\mathcal{L})\bigg(\bigg\{(x,t)\in X\times(0,\infty):\frac{|T_{t}%
f(x)|}{t^{1/p}}>\lambda\quad\text{and}\quad x\in B(a,R)\bigg\}\bigg)\leq
\frac{2c_{m}}{\lambda^{p}}\Vert f\Vert_{L^{p}(X,m )}^{p}.
\]
Taking limits on both sides of the previous inequality as $R\rightarrow\infty$
we arrive at the desired result.
\end{proof}

\begin{remark}
The proof of the previous theorem can be applied, with minor modifications, to
deal with more general operators as we now indicate. Let $(X,d,m)$ be as in
the previous theorem, and let $\nu$ be a Borel measure on the product space
$X\times(0,\infty).$ For $x\in X,t>0,$ we let $\mathcal{C}(x,t)=\frac
{\nu(\tilde{B}(x,4t))}{m(B(x,4t))},$ where $\tilde{B}(x,t)=B(x,t)\times(0,t)$
is a Carleson tent based on $B(x,t).$ For fixed $p\in\lbrack1,\infty)$ and a
locally integrable function $f$, we define
\[
T_{t}f(x)=\frac{1}{\mathcal{C}(x,t)^{1/p}}\int_{B(x,t)}f(y) \, \mathrm{d}
m(y), \quad x \in X, \quad t > 0.
\]
Then there exists a constant $C>0$ such that for all $\lambda>0$,
\[
\nu\bigg(\bigg\{(x,t)\in X\times(0,\infty):\frac{|T_{t}f(x)|}{\mathcal{C}%
(x,t)^{1/p}}>\lambda\bigg\}\bigg)\leq C\frac{\Vert f\Vert_{L^{p}(X,m)}^{p}%
}{\lambda^{p}}.
\]

\end{remark}

\begin{proof}
Following the proof of Theorem \ref{ThmDav} for each $\lambda>0$, let
\[
E_{\lambda}=\bigg\{(x,t)\in X\times(0,\infty):\frac{|T_{t}f(x)|}%
{\mathcal{C}(x,t)^{1/p}}>\lambda\bigg\}
\]
and
\[
E_{\lambda}^{\prime}=\bigg\{x\in X:\sup_{t>0}\frac{|T_{t}f(x)|}{\mathcal{C}%
(x,t)^{1/p}}>\lambda\bigg\}.
\]
Given $x\in E_{\lambda}^{\prime}$ we introduce
\[
t_{x}=\sup\bigg\{t>0:\frac{|T_{t}f(x)|}{\mathcal{C}(x,t)^{1/p}}>\lambda
\bigg\},
\]
and, as before, we infer that $t<2\,t_{x_{i}}$ (recall that the $t_{x_{i}}$'s
are the radius of the balls of the corresponding Vitali covering, see
\eqref{Sup1}). Following the argument above \textit{mutatis mutandis}, we get
\begin{align*}
\nu(E_{\lambda})  &  \leq\sum_{i}\nu(\tilde{B}(x_{i},4t_{x_{i}}))=\sum
_{i}m(B(x_{i},4t_{x_{i}}))\mathcal{C}(x_{i},t_{x_{i}})\\
&  \leq\frac{2c_{m}}{\lambda^{p}}\sum_{i}m(B(x_{i},t_{x_{i}}))|T_{t_{x_{i}}%
}f(x_{i})|^{p}\\
&  \leq\frac{2c_{m}}{\lambda^{p}}\sum_{i}\int_{B(x_{i},t_{x_{i}})}%
|f(y)|^{p}\,\mathrm{d}m(y)\\
&  \leq\frac{2c_{m}}{\lambda^{p}}\Vert f\Vert_{L^{p}(X,m)}^{p}.
\end{align*}

\end{proof}

\begin{remark}
\label{RemWeights} Note that in the special case $\nu=m\times\mathcal{L}$ one
has $\mathcal{C}(x,t)=4 t$. On the other hand, if $\omega$ is an arbitrary
weight on $(0,\infty)$ and we let $\nu=m\times\omega$ then $\mathcal{C}%
(x,t)=\int_{0}^{t}\omega(u)\,\mathrm{d} u.$
\end{remark}

\section{Examples and Applications of Theorem \ref{ThmMaximal} and Corollary
\ref{CorMaximal}}

\label{SubsectionAppl}

In this section we collect a number of applications of Theorem
\ref{ThmMaximal} to different areas of Analysis. The section and each of the
Examples can be read independently from the remaining parts of the paper.

Unless otherwise stated, in this section we assume $1 < p < \infty$.

\subsection{Brezis--Van Schaftingen--Yung formula\label{ExampleBSY}}

We show how to approach the Brezis--Van Schaftingen--Yung formula via Theorem
\ref{ThmMaximal}. Let us consider first functions defined on ${\mathbb{R}}$.

For $t>0,$ $f\in C^{1}({\mathbb{R}}),$ we let
\[
T_{t}f(x)=\frac{f(x+t)-f(x)}{t},\quad x\in{\mathbb{R}}.
\]
Then one has
\[
\lim_{t\rightarrow0+}T_{t}f(x)=f^{\prime}(x)
\]
and
\begin{equation}
\sup_{t>0}|T_{t}f(x)|\leq M(f^{\prime})(x). \label{BSY1}%
\end{equation}
Here $M$ is the (right) Hardy--Littlewood maximal operator defined for
integrable functions by
\[
M(g)(x)=\sup_{t>0}\frac{1}{t}\int_{x}^{x+t}|g(y)|\,\mathrm{d}y,\quad
x\in{\mathbb{R}}.
\]
The weighted version of the Brezis--Van Schaftingen--Yung formulae \eqref{l1*} and \eqref{I1'} that follow from Corollary \ref{CorMaximal} (cf. also Remark \ref{RemarkSharpAssertion}) and
the Hardy--Littlewood maximal theorem reads as 
\[
\lim_{\lambda\rightarrow\infty}\lambda^{p}\iint_{E_{\lambda,\gamma/p}%
}t^{\gamma-1}\,\mathrm{d}t\,\mathrm{d}x = \frac{1}{\gamma} \Vert
f^{\prime}\Vert_{L^{p}({\mathbb{R}})}^{p} \asymp\sup_{\lambda>0}\lambda^{p}%
\iint_{E_{\lambda,\gamma/p}}t^{\gamma-1}\,\mathrm{d}t\,\mathrm{d}x,\quad\gamma>0,
\]
where
\[
E_{\lambda,\gamma/p}=\{(x,t)\in{\mathbb{R}}\times(0,\infty
):|f(x+t)-f(x)|>\lambda\,t^{\frac{\gamma}{p}+1}\}.
\]
In particular, letting $\gamma=1$ we obtain
\begin{align*}
\Vert f^{\prime}\Vert_{L^{p}({\mathbb{R}})}^{p}  &  = \lim_{\lambda
\rightarrow\infty}\lambda^{p}\mathcal{L}^{2}\bigg(\bigg\{(x,t)\in{\mathbb{R}%
}\times(0,\infty):\frac{|f(x+t)-f(x)|}{t^{\frac{1}{p}+1}}>\lambda
\bigg\}\bigg)\\
&  \asymp\sup_{\lambda>0}\lambda^{p}\mathcal{L}^{2}\bigg(\bigg\{(x,t)\in
{\mathbb{R}}\times(0,\infty):\frac{|f(x+t)-f(x)|}{t^{\frac{1}{p}+1}}%
>\lambda\bigg\}\bigg).
\end{align*}

To deal with the $N$-dimensional case we use polar coordinates. Let
$\lambda>0$ and $\gamma>0$, then
\begin{equation}
\underset{|f(x)-f(y)|>\lambda|x-y|^{\frac{\gamma}{p}+1}}{\iint}\lambda
^{p}|x-y|^{\gamma-N}\,\mathrm{d}x\,\mathrm{d}y=\underset{|f(x)-f(x+t\omega
)|>\lambda t^{\frac{\gamma}{p}+1}}{\int_{\mathbb{S}^{N-1}}\int_{\mathbb{R}%
^{N}}\int_{0}^{\infty}}\lambda^{p}t^{\gamma-1}\,\mathrm{d}t\,\mathrm{d}%
x\,\mathrm{d}\sigma^{N-1}(\omega). \label{Polar}%
\end{equation}
Given $f\in W^{1, p}(\mathbb{R}^{N}),$ and $\omega\in\mathbb{S}^{N-1}$, we
claim that
\begin{equation}
\underset{|f(x)-f(x+t\omega)|>\lambda t^{\frac{\gamma}{p}+1}}{\int%
_{\mathbb{R}^{N}}\int_{0}^{\infty}}\lambda^{p}t^{\gamma-1}\,\mathrm{d}%
t\,\mathrm{d}x\leq\frac{C_{p}^{p}}{\gamma}\int_{\mathbb{R}^{N}}|\langle\nabla
f(x), \omega\rangle|^{p}\,\mathrm{d}x\qquad\text{for all}\qquad\lambda>0
\label{Polar1}%
\end{equation}
and
\begin{equation}
\frac{1}{\gamma}\int_{\mathbb{R}^{N}}|\langle\nabla f(x), \omega\rangle
|^{p}\,\mathrm{d}x = \lim_{\lambda\rightarrow\infty}%
\underset{|f(x)-f(x+t\omega)|>\lambda t^{\frac{\gamma}{p}+1}}{\int%
_{\mathbb{R}^{N}}\int_{0}^{\infty}}\lambda^{p}t^{\gamma-1}\,\mathrm{d}%
t\,\mathrm{d}x. \label{Polar2}%
\end{equation}

Assume momentarily that \eqref{Polar1}-\eqref{Polar2} hold. It follows from
\eqref{Polar}, \eqref{Polar1} and Fubini's theorem that
\[
\underset{|f(x)-f(y)|>\lambda|x-y|^{\frac{\gamma}{p}+1}}{\iint}\lambda
^{p}|x-y|^{\gamma-N}\,\mathrm{d}x\,\mathrm{d}y\leq\frac{C_{p}^{p}}{\gamma
}k(p,N)\int_{\mathbb{R}^{N}}|\nabla f(x)|^{p}\,\mathrm{d}x,\qquad\lambda>0,
\]
 where 
 \begin{equation}\label{ConstantkpN}
 k(p,N)=\int_{\mathbb{S}^{N-1}}|\langle e,\omega\rangle
|^{p}\,\mathrm{d}\sigma^{N-1}(\omega)
\end{equation}
 and $e$ is any unit vector in
$\mathbb{R}^{N}$. In light of \eqref{Polar2} and applying Lebesgue's dominated
convergence theorem, we get
\begin{equation}
\frac{k(p,N)}{\gamma}\int_{\mathbb{R}^{N}}|\nabla f(x)|^{p}\,\mathrm{d}%
x = \lim_{\lambda\rightarrow\infty}\underset{|f(x)-f(y)|>\lambda
|x-y|^{\frac{\gamma}{p}+1}}{\iint}\lambda^{p}|x-y|^{\gamma-N}\,\mathrm{d}%
x\,\mathrm{d}y. \label{BSYLimits}%
\end{equation}
Then
\begin{align*}
\lim_{\lambda\rightarrow\infty}\underset{|f(x)-f(y)|>\lambda|x-y|^{\frac
{\gamma}{p}+1}}{\iint}\lambda^{p}|x-y|^{\gamma-N}\,\mathrm{d}x\,\mathrm{d}y
&  \asymp\Vert\nabla f\Vert_{L^{p}(\mathbb{R}^{N})}^{p}\\
&  \hspace{-5cm}\asymp\sup_{\lambda>0}\underset{|f(x)-f(y)|>\lambda
|x-y|^{\frac{\gamma}{p}+1}}{\iint}\lambda^{p}|x-y|^{\gamma-N}\,\mathrm{d}%
x\,\mathrm{d}y
\end{align*}
and, in particular (letting $\gamma=N$)
\begin{align*}
\Vert\nabla f\Vert_{L^{p}({\mathbb{R}^{N}})}^{p}  &  \asymp\lim_{\lambda
\rightarrow\infty}\lambda^{p}\mathcal{L}^{2N}\bigg(\bigg\{(x,y)\in
{\mathbb{R}^{N}}\times\mathbb{R}^{N}:\frac{|f(x)-f(y)|}{|x-y|^{\frac{N}{p}+1}%
}>\lambda\bigg\}\bigg)\\
&  \asymp\sup_{\lambda>0}\lambda^{p}\mathcal{L}^{2N}\bigg(\bigg\{(x,y)\in
{\mathbb{R}^{N}}\times\mathbb{R}^{N}:\frac{|f(x)-f(y)|}{|x-y|^{\frac{N}{p}+1}%
}>\lambda\bigg\}\bigg).
\end{align*}

To prove \eqref{Polar1} and \eqref{Polar2}, we may assume without loss of
generality that $\omega=e_{N}=(0,\ldots,0,1)$. For $t>0,$ we let
\[
T_{t}f(x)=\frac{f(x+te_{N})-f(x)}{t},\qquad x\in{\mathbb{R}^{N}}.
\]
Then one has
\[
\lim_{t\rightarrow0+}T_{t}f(x)=\frac{\partial f}{\partial x_{N}}(x)
\]
and
\[
\sup_{t>0}|T_{t}f(x)|\leq M_{N}\Big(\frac{\partial f}{\partial x_{N}%
}\Big)(x).
\]
Here $M_{N}$ is the (right) Hardy--Littlewood maximal operator with respect to
the variable $x_{N}$ defined by
\[
M_{N}(g)(x^{\prime},x_{N})=\sup_{t>0}\frac{1}{t}\int_{x_{N}}^{x_{N}%
+t}|g(x^{\prime},y)|\,\mathrm{d}y,\quad x^{\prime}\in{\mathbb{R}^{N}},\quad
x_{N}\in\mathbb{R}.
\]
It follows from Corollary \ref{CorMaximal} (see also Remark \ref{RemarkSharpAssertion}) and the Hardy--Littlewood maximal
theorem that
\[
\sup_{\lambda>0}\lambda^{p}\iint_{E_{\lambda,\gamma/p}}t^{\gamma
-1}\,\mathrm{d}t\,\mathrm{d}x\leq\frac{C_{p}^{p}}{\gamma}\int_{\mathbb{R}^{N}%
}|\langle\nabla f(x),e_{N}\rangle|^{p}\,\mathrm{d}x
\]
and
\[
\frac{1}{\gamma}\int_{\mathbb{R}^{N}}|\langle\nabla f(x),e_{N}\rangle
|^{p}\,\mathrm{d}x = \lim_{\lambda\rightarrow\infty}\lambda^{p}%
\iint_{E_{\lambda,\gamma/p}}t^{\gamma-1}\,\mathrm{d}t\,\mathrm{d}x
\]
where
\[
E_{\lambda,\gamma/p}=\{(x,t)\in{\mathbb{R}^{N}}\times(0,\infty):|f(x+te_{N}%
)-f(x)|>\lambda\,t^{\frac{\gamma}{p}+1}\}.
\]

\subsection{Higher order Sobolev norms\label{ex:ver}}

The work of Alabern--Mateu--Verdera \cite{Verdera} (with the forerunners
\cite{Wheeden, Wheeden72}) will be used to construct a natural family of
operators $\{T_{t}\}$ that can be used to recover $\Vert\Delta f\Vert
_{L^{p}({\mathbb{R}}^{N})}.$ Our method combines \cite{Verdera} with Theorem
\ref{ThmMaximal}. Recall that in \cite{Verdera} the authors are able to
characterize $\Vert\Delta f\Vert_{L^{p}({\mathbb{R}}^{N})}$ in terms of square
maximal functions involving only first order difference operators. Apparently
the idea behind \cite{Verdera} is the fact that if $f$ is a smooth function
then, by Taylor's expansion,
\[
f(y)=f(x)+\sum_{0<|\alpha|\leq2}\frac{1}{\alpha!}D^{\alpha}f(x)(y-x)^{\alpha
}+O(|y-x|^{3}),
\]
therefore, taking averages over $y\in B(x,t)$ on both sides of the last
equation, one obtains
\begin{equation}
f(x)-\frac{1}{\mathcal{L}^{N}(B(x,t))}\int_{B(x,t)}f(y)\,\mathrm{d}y=-\frac
{1}{2(N+2)}\Delta f(x)t^{2}+o(t^{2}). \label{BSYHigh}%
\end{equation}
We are therefore led to consider the family of operators
\[
T_{t}f(x)=\frac{f(x)-\frac{1}{\mathcal{L}^{N}(B(x,t))}\int_{B(x,t)}%
f(y)\,\mathrm{d}y}{t^{2}},\quad x\in{\mathbb{R}}^{N},
\]
since by \eqref{BSYHigh} we have
\[
\lim_{t\rightarrow0+}T_{t}f(x)=-\frac{1}{2(N+2)}\Delta f(x)
\]
and
\[
\Big\|\sup_{t>0}|T_{t}f|\Big\|_{L^{p}({\mathbb{R}}^{N})}\lesssim\Vert\Delta
f\Vert_{L^{p}({\mathbb{R}}^{N})};
\]
see also \cite[formula (2.4), p. 91]{Dai}. Applying Corollary \ref{CorMaximal} and Remark \ref{RemarkSharpAssertion}
we infer that, for $\gamma > 0$,
\[\sup_{\lambda>0}\lambda^{p}%
\iint_{E_{\lambda,\gamma/p}}t^{\gamma-1}\mathrm{d}t\,\mathrm{d}x\asymp
\Vert\Delta f\Vert_{L^{p}({\mathbb{R}}^{N})}^{p}
\]
and
$$
\lim_{\lambda\rightarrow\infty}\lambda^{p}\iint_{E_{\lambda,\gamma/p}%
}t^{\gamma-1}\mathrm{d}t\,\mathrm{d}x = \frac{1}{2  (N+2) \gamma} \Vert\Delta f\Vert_{L^{p}({\mathbb{R}}^{N})}^{p}
$$
where
\[
E_{\lambda,\gamma/p}=\bigg\{(x,t)\in{\mathbb{R}}^{N}\times(0,\infty
):\Big|f(x)-\frac{1}{\mathcal{L}^{N}(B(x,t))}\int_{B(x,t)}f(y)\,\mathrm{d}%
y\Big|>\lambda\,t^{\frac{\gamma}{p}+2}\bigg\}.
\]
In particular, letting $\gamma=1$ we have
\begin{equation*}
\sup_{\lambda>0}\lambda^{p}\mathcal{L}^{N+1}\bigg(\bigg\{(x,t)\in
{\mathbb{R}}^{N}\times(0,\infty):\frac{|f(x)-\frac{1}{\mathcal{L}^{N}%
(B(x,t))}\int_{B(x,t)}f(y)\,\mathrm{d}y|}{t^{2+1/p}}>\lambda\bigg\}\bigg)  \asymp \Vert\Delta f\Vert_{L^{p}({\mathbb{R}}^{N})}^{p} 
\end{equation*}
and
$$
\lim_{\lambda
\rightarrow\infty}\lambda^{p}\mathcal{L}^{N+1}\bigg(\bigg\{(x,t)\in
{\mathbb{R}}^{N}\times(0,\infty):\frac{|f(x)-\frac{1}{\mathcal{L}^{N}%
(B(x,t))}\int_{B(x,t)}f(y)\,\mathrm{d}y|}{t^{2+1/p}}>\lambda\bigg\}\bigg) = \frac{1}{2 (N+2)} \Vert\Delta f\Vert_{L^{p}({\mathbb{R}}^{N})}^{p} .
$$

The method can be extended to deal with the Sobolev spaces $W^{2k,p}%
({\mathbb{R}}^{N}),\,k\in\mathbb{N}$. For further characterizations of smooth
function spaces in terms of ball averages, we refer to \cite{Dai} and
\cite[Section 10]{Dominguez}.

\subsection{Brezis--Van Schaftingen--Yung formula: the anisotropic case}

For the sake of simplicity, we restrict ourselves to functions on
${\mathbb{R}}^{2}$. Let
\[
T_{t,s}f(x,y)=\frac{1}{ts}\int_{x}^{x+t}\int_{y}^{y+s}f(u,v)\,\mathrm{d}%
u\,\mathrm{d}v,\quad(x,y)\in{\mathbb{R}}^{2},\quad(t,s)\in(0,\infty)^{2}.
\]
The associated maximal operator, usually called the strong maximal operator,
is bounded from $L^{p}({\mathbb{R}}^{2})$ to $L^{p}({\mathbb{R}}^{2})$, i.e.,
\[
\Big\|\sup_{t,s>0}|T_{t,s}f|\Big\|_{L^{p}({\mathbb{R}}^{2})}\lesssim\Vert
f\Vert_{L^{p}({\mathbb{R}}^{2})}.
\]
This follows from the fact that $T_{t,s}$ can be expressed as the composition
of one-dimensional Hardy--Littlewood maximal functions together with the
Fubini property on $L^{p}({\mathbb{R}}^{2})$. By Remark \ref{RemMulti}, we
get
\begin{equation}
\lim_{\lambda\rightarrow\infty}\lambda^{p}\mathcal{L}^{4}(E_{\lambda}%
)\asymp\sup_{\lambda>0}\lambda^{p}\mathcal{L}^{4}(E_{\lambda})\asymp
\Vert\partial_{x}\partial_{y}f\Vert_{L^{p}({\mathbb{R}}^{2})}^{p} \label{cccc}%
\end{equation}
where
\[
E_{\lambda}=\bigg\{(x,y,t,s):\frac{|f(x+t,y+s)-f(x,y+s)-f(x+t,y)+f(x,y)|}%
{ts\max\{t,s\}^{2/p}}>\lambda\bigg\}.
\]
Comparing this result with its isotropic counterpart given by (cf. \eqref{l1*}
and \eqref{I1'})
\[
\lim_{\lambda\rightarrow\infty}\lambda^{p}\mathcal{L}^{4}%
\bigg(\bigg\{(x,y,t,s):\frac{|f(x+t,y+s)-f(x,y)|}{\max\{t,s\}\max
\{t,s\}^{2/p}}>\lambda\bigg\}\bigg)\asymp\Vert\nabla f\Vert_{L^{p}%
({\mathbb{R}}^{2})}^{p},
\]
we observe that the factor $\max\{t,s\}$ related to the classical difference
operator (i.e., $f(x+t,y+s)-f(x,y)$) is replaced in \eqref{cccc} by the factor
$ts$ related to the mixed difference operator (i.e.,
$f(x+t,y+s)-f(x,y+s)-f(x+t,y)+f(x,y)$).

Further characterizations of Sobolev spaces in terms of anisotropic and
magnetic fields may be found in \cite{ns19}, \cite{npsv18} and the references therein.

\subsection{Symmetric Markov Semigroups}

Let $(X,m)$ be a $\sigma$-finite measure space, and let $\{T_{t}:t>0\}$ be a
symmetric Markov semigroup acting on $L^{2}(X,m)$. This class of semigroups
has a canonical extension to a contraction semigroup in all $L^{p}(X,m)$
spaces; cf. \cite{Stein} and \cite{Ledoux}. Then the Rota's maximal theorem
asserts that
\begin{equation}
\Big\|\sup_{t>0}|T_{t}f|\Big\|_{L^{p}(X,m)}\lesssim\Vert f\Vert_{L^{p}(X,m)}.
\label{Rota}%
\end{equation}
In addition, for all $f\in L^{p}(X,m)$, one has $\lim_{t\rightarrow0+}%
T_{t}f(x)=f(x)$ $m$-a.e. $x\in X$ (cf. \cite[pg. 73]{Stein}, or \cite[Lemma
1.6.2]{Ledoux}). According to Corollary \ref{CorMaximal} and Remark \ref{RemarkSharpAssertion} we obtain, for each $\gamma > 0$,
\[
\sup_{\lambda> 0}\lambda
^{p}\iint_{E_{\lambda,\gamma/p}}t^{\gamma-1}\,\mathrm{d}t\,\mathrm{d}m(x)
\asymp\Vert f\Vert_{L^{p}(X,m)}^{p}
\]
and
$$
 \lim_{\lambda\rightarrow\infty}\lambda^{p}\iint_{E_{\lambda,\gamma/p}%
}t^{\gamma-1}\,\mathrm{d}t\,\mathrm{d}m(x) = \frac{1}{\gamma} \Vert f\Vert_{L^{p}(X,m)}^{p}
$$
where
\[
E_{\lambda,\gamma/p}=\{(x,t)\in X\times(0,\infty):|T_{t}f(x)|>\lambda
\,t^{\gamma/p}\}.
\]
In particular
$$
 \sup_{\lambda> 0}\lambda^{p}(m\times\mathcal{L})\bigg(\bigg\{(x,t)\in
X\times(0,\infty):\frac{|T_{t}f(x)|}{t^{1/p}}>\lambda\bigg\}\bigg) \asymp \Vert f\Vert_{L^{p}(X,m)}^{p} 
$$
and
$$
	 \lim_{\lambda\rightarrow\infty}%
\lambda^{p}(m\times\mathcal{L})\bigg(\bigg\{(x,t)\in X\times(0,\infty
):\frac{|T_{t}f(x)|}{t^{1/p}}>\lambda\bigg\}\bigg) = \Vert f\Vert_{L^{p}(X,m)}^{p}.
$$
We point out that the proof of the maximal inequality \eqref{Rota} given in
\cite{Stein} and \cite{Ledoux} strongly relies on the symmetry properties of
the semigroup. The next example concerns with non-symmetric semigroups related
to an important class of PDE's.

\subsection{H\"{o}rmander semigroups\label{ExampleHormander}}

We consider an example connected with the class of second order PDE's
introduced by H\"{o}rmander in his seminal paper \cite{Hormander} (cf.
\eqref{Hor}). Let $\mathscr{P}_{t}$ be the Poisson semigroup associated to
\eqref{Hor} given by
\begin{equation}
\mathscr{P}_{t}f(x)=\frac{1}{\sqrt{4\pi}}\int_{0}^{\infty}\frac{t}{\xi^{3/2}%
}e^{-\frac{t^{2}}{4\xi}}P_{\xi}f(x)\,\mathrm{d}\xi,\quad f\in\mathcal{S}%
({\mathbb{R}}^{N}), \label{1*p}%
\end{equation}
where $P_{\xi}$ is the H\"{o}rmander semigroup, that is, $P_{\xi}f$ is the
unique solution of the Cauchy problem
\[
\left\{
\begin{array}
[c]{ll}%
u(x,0)=f(x) & \text{on}\quad\mathbb{R}^{N},\\
\mathscr{K}u=0 & \text{in}\quad\mathbb{R}_{+}^{N+1}.
\end{array}
\right.
\]
Note that $\mathscr{P}_{t}$ is non-doubling and non-symmetric semigroup, so
that one can not apply directly Rota's theorem \eqref{Rota}. However, it was
recently shown by Garofalo and Tralli \cite[Theorem 5.5]{Garofalo} that
\eqref{Rota} still holds true for this special class of semigroups on
$L^{p}({\mathbb{R}}^{N})$. Namely, they show that if $\text{tr}\,B\geq0$ then
\[
\Big\|\sup_{t>0}|\mathscr{P}_{t}f|\Big\|_{L^{p}({\mathbb{R}}^{N})}%
\lesssim\Vert f\Vert_{L^{p}({\mathbb{R}}^{N})}.
\]
Hence, by Corollary \ref{CorMaximal} and Remark \ref{RemarkSharpAssertion}, we get, for $\gamma > 0$,
\[
\sup_{\lambda> 0}\lambda^{p}%
\iint_{E_{\lambda,\gamma/p}}t^{\gamma-1}\,\mathrm{d}t\,\mathrm{d}x\asymp\Vert
f\Vert_{L^{p}({\mathbb{R}}^{N})}^{p}
\]
and
$$
\lim_{\lambda\rightarrow\infty}\lambda^{p}\iint_{E_{\lambda,\gamma/p}%
}t^{\gamma-1}\,\mathrm{d}t\,\mathrm{d}x = \frac{1}{\gamma} \Vert
f\Vert_{L^{p}({\mathbb{R}}^{N})}^{p}
$$
where
\[
E_{\lambda,\gamma/p}=\Big\{(x,t)\in{\mathbb{R}}^{N}\times(0,\infty
):|\mathscr{P}_{t}f(x)|>\lambda\,t^{\gamma/p}\Big\}.
\]
In particular,
\begin{align}
  \sup_{\lambda> 0}\lambda^{p}\mathcal{L}^{N+1}\bigg(\bigg\{(x,t)\in
{\mathbb{R}}^{N}\times(0,\infty):\frac{|\mathscr{P}_{t}f(x)|}{t^{1/p}}%
>\lambda\bigg\}\bigg) \asymp \Vert f\Vert_{L^{p}({\mathbb{R}}^{N})}^{p} \label{GT'}%
\end{align}
and
\begin{equation}\label{GT'6626362}
	\lim_{\lambda
\rightarrow\infty}\lambda^{p}\mathcal{L}^{N+1}\bigg(\bigg\{(x,t)\in
{\mathbb{R}}^{N}\times(0,\infty):\frac{|\mathscr{P}_{t}f(x)|}{t^{1/p}}%
>\lambda\bigg\}\bigg) =  \Vert f\Vert_{L^{p}({\mathbb{R}}^{N})}^{p}.
\end{equation}
Since
\begin{equation}
u(x,t)=\mathscr{P}_{t}f(x), \quad x \in{\mathbb{R}}^{N}, \quad t > 0,
\label{1*}%
\end{equation}
satisfies
\begin{equation}
\left\{
\begin{array}
[c]{ll}%
u(x,0)=f(x) & \text{on}\quad\mathbb{R}^{N},\\
u_{tt}+\mathscr{A}u=0 & \text{in}\quad\mathbb{R}_{+}^{N+1},
\end{array}
\right.  \label{GT'3}%
\end{equation}
we can rewrite \eqref{GT'} and \eqref{GT'6626362} as
\begin{equation}
\sup_{\lambda> 0}\lambda
^{p}\mathcal{L}^{N+1}\bigg(\bigg\{(x,t)\in{\mathbb{R}}^{N}\times
(0,\infty):\frac{|u(x,t)|}{t^{1/p}}>\lambda\bigg\}\bigg) \asymp \Vert f\Vert_{L^{p}({\mathbb{R}}^{N})}^{p}\label{GT'2}%
\end{equation}
or
\begin{equation}
\label{GT'2New}
\lim
_{\lambda\rightarrow\infty}\lambda^{p}\mathcal{L}^{N+1}\bigg(\bigg\{(x,t)\in
{\mathbb{R}}^{N}\times(0,\infty):\frac{|u(x,t)|}{t^{1/p}}>\lambda\bigg\}\bigg) = \Vert f\Vert_{L^{p}({\mathbb{R}}^{N})}^{p},
\end{equation}
respectively. This establishes an interesting link between the Cauchy problem related to an
important class of PDE's (cf. \eqref{GT'3}) and the Brezis--Van
Schaftingen--Yung condition given by the left-hand sides of \eqref{GT'2}
and  \eqref{GT'2New}. In the special case $\mathscr{A}=\Delta$ (i.e., taking
$Q=I_{N}$ and $B=0_{N}$ in \eqref{Hor}), one has that $u$ defined by
\eqref{1*} and \eqref{1*p} is the classical Poisson integral for the
half-space ${\mathbb{R}}_{+}^{N+1}$ (cf. \eqref{ooo})
\begin{equation}
u(x,t) = \frac{\Gamma\Big(\frac{N+1}{2}\Big)}{\pi^{\frac{N+1}{2}}}
P[f](x,t)=\frac{\Gamma\Big(\frac{N+1}{2}\Big)}{\pi^{\frac{N+1}{2}}}%
\int_{{\mathbb{R}}^{N}}\frac{t}{(|x-y|^{2} + t^{2})^{\frac{N+1}{2}}%
}f(y)\,\mathrm{d}y \label{Poisson}%
\end{equation}
and thus, e.g., \eqref{GT'2New} reads as
\[
\lim_{\lambda\rightarrow
\infty}\lambda^{p}\mathcal{L}^{N+1}\bigg(\bigg\{(x,t)\in{\mathbb{R}}^{N}%
\times(0,\infty):\frac{|P[f](x,t)|}{t^{1/p}}>\lambda\bigg\}\bigg) = \bigg(\frac{\pi^{\frac{N+1}{2}}} {\Gamma\Big(\frac{N+1}{2}\Big)} \bigg)^p \Vert f\Vert_{L^{p}({\mathbb{R}}^{N})}^{p}.
\]

\subsection{Spherical means and wave equation}

Let $N\geq2$. We consider the spherical maximal function
\[
M_{S}f(x)=\sup_{t>0}|A_{t}f(x)|,\quad x\in{\mathbb{R}}^{N},
\]
where $A_{t}$ are the spherical averaging operators
\[
A_{t}f(x)=\int_{\mathbb{S}^{N-1}}f(x-t\omega)\,\mathrm{d}\sigma^{N-1}%
(\omega),
\]
and $\mathrm{d}\sigma^{N-1}$ is normalized surface measure on the sphere
$\mathbb{S}^{N-1}$. Here, to avoid technical issues related to the
measurability of $M_{S}f$, we restrict our attention to functions $f$ in the
Schwartz class $\mathcal{S}({\mathbb{R}}^{N})$. The maximal theorem for
spherical means, due to Stein \cite{Stein76} if $N\geq3,$ and Bourgain
\cite{Bourgain} if $N=2$ (cf. also \cite{Mockenhaupt}), claims that
\begin{equation}
\Vert M_{S}f\Vert_{L^{p}({\mathbb{R}}^{N})}\lesssim\Vert f\Vert_{L^{p}%
({\mathbb{R}}^{N})},\quad p>\frac{N}{N-1}. \label{BS}%
\end{equation}
If we apply Corollary \ref{CorMaximal} together with Lemma \ref{RemarkSharpAssertion}, we have, for $\gamma > 0$,
\begin{equation}
\sup_{\lambda> 0 }\lambda
^{p}\iint_{E_{\lambda,\gamma/p}}t^{\gamma-1}\,\mathrm{d}t\,\mathrm{d}%
x\asymp\Vert f\Vert_{L^{p}({\mathbb{R}}^{N})}^{p} \label{BS1+}%
\end{equation}
and 
\begin{equation}
\lim_{\lambda\rightarrow\infty}\lambda^{p}\iint_{E_{\lambda,\gamma/p}%
}t^{\gamma-1}\,\mathrm{d}t\,\mathrm{d}x = \frac{1}{\gamma} \Vert f\Vert_{L^{p}({\mathbb{R}}^{N})}^{p}  \label{BS1+++++++}
\end{equation}
where
\[
E_{\lambda,\gamma/p}=\bigg\{(x,t)\in{\mathbb{R}}^{N}\times(0,\infty
):\Big|\int_{\mathbb{S}^{N-1}}f(x-t\omega)\,\mathrm{d} \sigma^{N-1}%
(\omega)\Big|>\lambda\,t^{\gamma/p}\bigg\}.
\]
In particular, if $\gamma=1$ then
\begin{equation}
\Vert f\Vert_{L^{p}({\mathbb{R}}^{N})}^{p}    \asymp\sup_{\lambda> 0}\lambda^{p}\mathcal{L}^{N+1}\bigg(\bigg\{(x,t)\in
{\mathbb{R}}^{N}\times(0,\infty):\frac{|\int_{\mathbb{S}^{N-1}}f(x-t\omega
)\,\mathrm{d} \sigma^{N-1}(\omega)|}{t^{1/p}}>\lambda\bigg\}\bigg)
\label{BS1}%
\end{equation}
and
\begin{equation}
	\Vert f\Vert_{L^{p}({\mathbb{R}}^{N})}^{p} = \lim_{\lambda
\rightarrow\infty}\lambda^{p}\mathcal{L}^{N+1}\bigg(\bigg\{(x,t)\in
{\mathbb{R}}^{N}\times(0,\infty):\frac{|\int_{\mathbb{S}^{N-1}}f(x-t\omega
)\,\mathrm{d} \sigma^{N-1}(\omega)|}{t^{1/p}}>\lambda\bigg\}\bigg). \label{BS1nenennwnw}%
\end{equation}

The previous analysis can be extended to deal with the more general class of
spherical means given by
\[
M_{S}^{\alpha}f(x)=\sup_{t>0}|A_{t}^{\alpha}f(x)|, \quad x \in{\mathbb{R}}%
^{N},
\]
where
\[
\widehat{A_{t}^{\alpha}f}(\xi)=\widehat{m}_{\alpha}(\xi t)\widehat{f}(\xi)
\]
and $\widehat{m}_{\alpha}(\xi)=\pi^{-\alpha+1}|\xi|^{-N/2-\alpha
+1}J_{N/2+\alpha-1}(2\pi|\xi|)$; see \cite{Stein76}. Here, as usual,
$\widehat{f}$ denotes the Fourier transform of $f$. In particular, $A_{t}%
^{0}f$ is a constant multiple of $A_{t}f$. The extension of \eqref{BS}, when
$M_{S}$ is replaced by $M_{S}^{\alpha},$ was obtained in \cite{Stein76} under
the assumption that $N\geq3$ and one of the following conditions holds

\begin{enumerate}
\item $1 < p \leq2$, when $\alpha> 1-N+N/p$,

\item $2 \leq p \leq\infty$, when $\alpha> (2-N)/p$.
\end{enumerate}

Accordingly, it follows from Corollary \ref{CorMaximal} that
\begin{align}
\Vert f\Vert_{L^{p}({\mathbb{R}}^{N})}^{p}  &  \asymp\lim_{\lambda
\rightarrow\infty}\lambda^{p}\mathcal{L}^{N+1}\bigg(\bigg\{(x,t)\in
{\mathbb{R}}^{N}\times(0,\infty):\frac{|A_{t}^{\alpha}f(x)|}{t^{1/p}}%
>\lambda\bigg\}\bigg)\nonumber\\
&  \asymp\sup_{\lambda> 0}\lambda^{p}\mathcal{L}^{N+1}\bigg(\bigg\{(x,t)\in
{\mathbb{R}}^{N}\times(0,\infty):\frac{|A_{t}^{\alpha}f(x)|}{t^{1/p}}%
>\lambda\bigg\}\bigg). \label{BS2}%
\end{align}
Of course, in the special case $\alpha=0$ one recovers the equivalences given
in \eqref{BS1} and \eqref{BS1nenennwnw}.

There is an interesting connection between \eqref{BS2} and wave equation
\begin{equation}
\label{wave}\Delta u = u_{tt} \quad\text{on} \quad{\mathbb{R}}^{N}
\times(0,\infty)
\end{equation}
with $u(x,0)= 0$ and $u_{t}(x,0) = f(x)$. Indeed, a weak solution of this
problem is given by $u(x,t) = t A^{\alpha}_{t} f(x)$ with $\alpha= (3-N)/2$
(cf. \cite{Stein76}). Therefore, if $N \geq3$ and $\frac{2 N}{N+1} < p <
\frac{2(N-2)}{N-3}$ then we derive
\begin{equation}
\label{wave2}\|f\|_{L^{p}({\mathbb{R}}^{N})}^{p} \asymp\lim_{\lambda\to\infty}
\lambda^{p} \mathcal{L}^{N+1}\bigg( \bigg\{(x,t) \in{\mathbb{R}}^{N}
\times(0,\infty) : \frac{|u(x,t)|}{t^{1 + 1/p}} > \lambda\bigg\} \bigg)
\end{equation}
and
\begin{equation}
\label{wave2'}\|f\|_{L^{p}({\mathbb{R}}^{N})}^{p} \asymp\sup_{\lambda> 0}
\lambda^{p} \mathcal{L}^{N+1}\bigg( \bigg\{(x,t) \in{\mathbb{R}}^{N}
\times(0,\infty) : \frac{|u(x,t)|}{t^{1 + 1/p}} > \lambda\bigg\} \bigg).
\end{equation}
Therefore the datum $f$ related to \eqref{wave} belongs to $L^{p}({\mathbb{R}%
}^{N})$ if and only if the corresponding solution $u$ satisfies the
Brezis--Van Schaftingen--Yung-type condition given in the right-hand sides of
\eqref{wave2} and \eqref{wave2'}.

We also mention that the characterizations \eqref{BS1+}--\eqref{BS1nenennwnw} have
an analogue for spherical means on the Heisenberg group $\mathbb{H}^{N}$.
Indeed, we can follow the above methodology but now applying the corresponding
maximal theorems obtained in \cite{Muller} and \cite{Thangavelu}. We must
leave further details to the reader.

\subsection{Maximal operators of convolution type associate to PDE's}

Let $a,b\geq0$ with $(a,b)\neq(0,0)$. Consider the class of PDE's given by
\begin{equation}
au_{tt}-bu_{t}+\Delta u=0\quad\text{in}\quad{\mathbb{R}}^{N}\times(0,\infty).
\label{Car1}%
\end{equation}
This class of equations includes as distinguished examples the Laplace
equation ($a=1$ and $b=0$) and heat equation ($a=0$ and $b=1$).

According to \cite{Carneiro}, the solution of the equation \eqref{Car1} with
initial datum $u_{a,b}(x,0)=f(x)$ is given by
\[
u_{a,b}(x,t)=\varphi_{a,b}(\cdot,t)\ast|f|(x), \quad x \in{\mathbb{R}}^{N},
\quad t > 0,
\]
where
\[
\widehat{\varphi}_{a,b}(\xi,t)=e^{-t\Big(\frac{-b+\sqrt{b^{2}+16a\pi^{2}%
|\xi|^{2}}}{2a}\Big)}.
\]
Moreover, arguing as in \cite[Theorem 2, Chapter III]{Stein70b} shows that the
corresponding maximal function
\[
u_{a,b}^{\ast}(x)=\sup_{t>0}\varphi_{a,b}(\cdot,t)\ast|f|(x)
\]
satisfies $u_{a,b}^{\ast}(x)\leq Mf(x)$, where $M$ denotes the
Hardy--Littlewood maximal operator on ${\mathbb{R}}^{N}$. Thus
\[
\Vert u_{a,b}^{\ast}\Vert_{L^{p}({\mathbb{R}^{N}})}\lesssim\Vert f\Vert
_{L^{p}({\mathbb{R}}^{N})}.
\]
Corollary \ref{CorMaximal} and Remark \ref{RemarkSharpAssertion} imply
\begin{equation}
\lim_{\lambda\rightarrow\infty}\lambda^{p}\mathcal{L}^{N+1}%
\bigg(\bigg\{(x,t)\in{\mathbb{R}}^{N}\times(0,\infty):\frac{|u_{a,b}%
(x,t)|}{t^{1/p}}>\lambda\bigg\}\bigg) = \Vert f\Vert_{L^{p}({\mathbb{R}%
}^{N})}^{p} \label{C}%
\end{equation}
and a similar result holds true when we replace $\lim_{\lambda\to\infty}$ by
$\sup_{\lambda> 0}$ (and changing $=$ by $\asymp$). For Laplace's equation, the previous characterization
coincides with \eqref{GT'2New} and \eqref{Poisson}. On the other hand, if
$a=0$ and $b=1$ in \eqref{C} we obtain a new characterization of Lebesgue
norms in terms of the solution of the heat equation
\[
u(x,t)=\frac{1}{(4\pi t)^{\frac{N}{2}}}\int_{{\mathbb{R}}^{N}}e^{-\frac
{|x-y|^{2}}{4t}}f(y)\,\mathrm{d}y.
\]

\subsection{Poisson integral on the sphere}

Weak-type estimates for functions defined on the circle and their harmonic
extensions on the disk were recently investigated by Greco and Schiattarella
\cite{Greco}. We can now complement their results by establishing the
counterparts of \eqref{GT'2New}--\eqref{Poisson} in terms of the Poisson
integrals $u(\theta,\rho)=u(\theta\rho),\,\theta\in\mathbb{S}^{N-1},\,\rho
\in\lbrack0,1)$,
\[
u(\theta,\rho)=\int_{\mathbb{S}^{N-1}}\frac{1-\rho^{2}}{|\rho\theta-\eta|^{N}%
}|f(\eta)|\,\mathrm{d} \sigma^{N-1}(\eta).
\]
It is well known that given $f\in L^{p}(\mathbb{S}^{N-1},\mathrm{d}%
\sigma^{N-1})$ the function $u(\theta,\rho)$ provides the solution of the
Laplace's equation $\Delta u=0$ on the unit $N$-dimensional open ball $B(0,1)$
with $\lim_{\rho\rightarrow1-}u(\theta,\rho)=|f(\theta)|$ for a.e. $\theta
\in\mathbb{S}^{N-1}$. Furthermore, the following pointwise inequality holds
(cf. \cite[Theorem 2.3.6, p. 39]{DaiBook})
\[
\sup_{\rho\in\lbrack0,1)}|u(\theta,\rho)|\leq\mathcal{M}f(\theta)
\]
where $\mathcal{M}$ is the Hardy--Littlewood maximal function (taken with
respect to geodesic balls). Since $\mathcal{M}$ acts boundedly on
$L^{p}(\mathbb{S}^{N-1},\mathrm{d}\sigma^{N-1})$ (cf. \cite[Corollary 2.3.4,
p. 38]{DaiBook}), we immediately get that
\[
\Big\|\sup_{\rho\in\lbrack0,1)}|u(\theta,\rho)|\Big\|_{L^{p}(\mathbb{S}%
^{N-1},\mathrm{d}\sigma^{N-1})}\lesssim\Vert f\Vert_{L^{p}(\mathbb{S}%
^{N-1},\mathrm{d}\sigma^{N-1})}.
\]
Therefore, after a simple change of variables, we can invoke Corollary
\ref{CorMaximal} to arrive at
\begin{align*}
\lim_{\lambda\rightarrow\infty}\lambda^{p}(\mathrm{d}\sigma^{N-1}%
\times\mathcal{L})\bigg(\bigg\{(\theta,\rho)\in\mathbb{S}^{N-1}\times
(0,1):\frac{|u(\theta,1-\rho)|}{(1-\rho)^{1/p}}>\lambda\bigg\}\bigg)& = \Vert f\Vert_{L^{p}(\mathbb{S}^{N-1},\mathrm{d}\sigma^{N-1})}^{p}  \\
& \hspace{-10cm} \asymp\sup_{\lambda> 0}\lambda^{p}(\mathrm{d}\sigma^{N-1}\times
\mathcal{L})\bigg(\bigg\{(\theta,\rho)\in\mathbb{S}^{N-1}\times(0,1):\frac
{|u(\theta,1-\rho)|}{(1-\rho)^{1/p}}>\lambda\bigg\}\bigg).
\end{align*}

On the other hand, by a similar argument as above, one can show that the
solution $u(\theta,t)$ of the heat equation
\[
\left\{
\begin{array}
[c]{ll}%
u(\theta,0)=|f(\theta)| & \text{on}\quad\mathbb{S}^{N-1},\\
u_{t}-\Delta u=0 & \text{in}\quad\mathbb{S}^{N-1}\times(0,\infty),
\end{array}
\right.
\]
enables us to characterize $\Vert f\Vert_{L^{p}(\mathbb{S}^{N-1}, \mathrm{d}
\sigma^{N-1})}$. Specifically, we have
\begin{align*}
\lim_{\lambda\rightarrow\infty}\lambda^{p}(\mathrm{d}\sigma^{N-1}%
\times\mathcal{L})\bigg(\bigg\{(\theta,t)\in\mathbb{S}^{N-1}\times
(0,\infty):\frac{|u(\theta,t)|}{t^{1/p}}>\lambda\bigg\}\bigg) &= \Vert f\Vert_{L^{p}(\mathbb{S}^{N-1},\mathrm{d}\sigma^{N-1})}^{p}  \\
& \hspace{-10cm} \asymp\sup_{\lambda> 0}\lambda^{p}(\mathrm{d}\sigma^{N-1}\times
\mathcal{L})\bigg(\bigg\{(\theta,t)\in\mathbb{S}^{N-1}\times(0,\infty
):\frac{|u(\theta,t)|}{t^{1/p}}>\lambda\bigg\}\bigg).
\end{align*}

\subsection{Ergodic Theory}

We start by recalling the (local) ergodic theorem by Wiener \cite{Wiener}. Let
$\{U_{t}:t>0\}$ be a measure-preserving flow on a measure space $(X,m)$ and
suppose that $f$ is a locally integrable function. Then
\[
\lim_{t\rightarrow0+}\frac{1}{2t}\int_{-t}^{t}f(U_{s}x)\,\mathrm{d}s=f(x)\quad
m-\text{a.e.}\quad x\in X.
\]
Assume further that $(X,d,m)$ is doubling. In this case it is well known that
the Hardy--Littlewood maximal inequality for $L^{p}(X,m)$ holds and applying
the Calder\'{o}n transference principle \cite{Calderon} we see that the
maximal ergodic operator
\[
M^{\ast}f(x)=\sup_{t>0}\bigg|\frac{1}{t}\int_{0}^{t}f(U_{s}x)\,\mathrm{d}%
s\bigg|
\]
is bounded on $L^{p}(X,m)$, that is,
\[
\Vert M^{\ast}f\Vert_{L^{p}(X,m)}\lesssim\Vert f\Vert_{L^{p}(X,m)}.
\]
According to Corollary \ref{CorMaximal}, one has
\[
\lim_{\lambda\rightarrow\infty}\lambda^{p}\iint_{E_{\lambda,\gamma/p}%
}t^{\gamma-1}\,\mathrm{d}t\,\mathrm{d}m(x) = \frac{1}{\gamma} \Vert f\Vert_{L^{p}(X,m)}^{p} \asymp\sup_{\lambda> 0}\lambda
^{p}\iint_{E_{\lambda,\gamma/p}}t^{\gamma-1}\,\mathrm{d}t\,\mathrm{d}m(x),\quad\gamma>0,
\]
where
\[
E_{\lambda,\gamma/p}=\bigg\{(x,t)\in X\times(0,\infty):\Big|\frac{1}{t}%
\int_{0}^{t}f(U_{s}x)\,\mathrm{d}s\Big|>\lambda\,t^{\gamma/p}\bigg\}.
\]
As a special case, we derive
\begin{align*}
\Vert f\Vert_{L^{p}(X,m)}^{p}  &  = \lim_{\lambda\rightarrow\infty}%
\lambda^{p}(m\times\mathcal{L})\bigg(\bigg\{(x,t)\in X\times(0,\infty
):\frac{|\frac{1}{t}\int_{0}^{t}f(U_{s}x)\,\mathrm{d}s|}{t^{1/p}}%
>\lambda\bigg\}\bigg)\\
&  \asymp\sup_{\lambda> 0}\lambda^{p}(m\times\mathcal{L})\bigg(\bigg\{(x,t)\in
X\times(0,\infty):\frac{|\frac{1}{t}\int_{0}^{t}f(U_{s}x)\,\mathrm{d}%
s|}{t^{1/p}}>\lambda\bigg\}\bigg).
\end{align*}

\subsection{Martingale Differences}

For the sake of simplicity, we restrict ourselves to the Haar system
$\{H_{n}\}$ in $L^{p}(0,1)$, but the results given below can be extended to
any complete system of martingale differences. For every $n\in{{\mathbb{N}}}$,
we write
\[
S_{n}f(x)=\sum_{\nu=0}^{n}a_{\nu}H_{\nu}(x),\quad a_{\nu}=\int_{0}^{1}H_{\nu
}(x)f(x)\,\mathrm{d}x.
\]
If $f\in L^{p}(0,1)$ then
\[
\lim_{n\rightarrow\infty}S_{n}f(x)=f(x)\quad\text{a.e.}\quad x\in(0,1)
\]
and the maximal function $S^{\ast}f(x)=\sup_{n\in{{\mathbb{N}}}}|S_{n}f(x)|$
maps $L^{p}(0,1)$ into itself (cf. \cite[Theorem 3.4.2, p. 72]{Garsia}),
\[
\Vert S^{\ast}f\Vert_{L^{p}(0,1)}\lesssim\Vert f\Vert_{L^{p}(0,1)}.
\]
As a by-product of Corollary \ref{CorMaximal} we infer that
\[
\lim_{\lambda\rightarrow\infty}\lambda^{p}\iint_{E_{\lambda,\gamma/p}%
}t^{\gamma-1}\,\mathrm{d}t\,\mathrm{d}x  = \frac{1}{\gamma}\Vert f\Vert_{L^{p}(0,1)}^{p}\asymp\sup_{\lambda> 0}\lambda
^{p}\iint_{E_{\lambda,\gamma/p}}t^{\gamma-1}\,\mathrm{d}t\,\mathrm{d}x,\quad\gamma>0,
\]
where
\[
E_{\lambda,\gamma/p}=\{(x,t)\in(0,1)^{2}:|S_{\floor{1/t}}f(x)|>\lambda
\,t^{\gamma/p}\}.
\]
As usual, $\floor{x}$ denotes the integer part of the real number $x$. In
particular,
\begin{align*}
\Vert f\Vert_{L^{p}(0,1)}^{p}  & =\lim_{\lambda\rightarrow\infty}%
\lambda^{p}\mathcal{L}^{2}\bigg(\bigg\{(x,t)\in(0,1)^{2}:\frac
{|S_{\floor{1/t}}f(x)|}{t^{1/p}}>\lambda\bigg\}\bigg)\\
&  \asymp\sup_{\lambda> 0}\lambda^{p}\mathcal{L}^{2}\bigg(\bigg\{(x,t)\in
(0,1)^{2}:\frac{|S_{\floor{1/t}}f(x)|}{t^{1/p}}>\lambda\bigg\}\bigg).
\end{align*}

Corollary \ref{CorMaximal} relies on the well-known fact that maximal
inequalities imply pointwise almost everywhere convergence of sequences of
operators. The converse statement is in general not true. However, a
fundamental result due to Stein \cite{Stein61} asserts that, for some special
classes of ambient spaces and sequences of operators, maximal inequalities are
indeed necessary to establish pointwise almost everywhere convergence. As a
consequence, we can apply Theorem \ref{ThmMaximal} from \textit{apriori}
pointwise convergence statements to obtain $L^{p}$-characterizations related
to sequences of operators. This is the content of the following application.

\subsection{Stein maximal principle and partial Fourier series}

Let $G$ be a compact group, let $X$ be a homogeneous space of $G$ with a
finite Haar measure $m,$ and let $1<p\leq2$. Let $T_{n}:L^{p}(X,m)\rightarrow
L^{p}(X,m)$ be a sequence of bounded linear operators commuting with
translations, such that for each $f\in L^{p}(X, m),$ $T_{n}f$ converges almost
everywhere to $f$. Then, by Stein's maximal principle \cite{Stein61},
\begin{equation}
\Big\|\sup_{n\in{{\mathbb{N}}}}|T_{n}f|\Big\|_{L^{p}(X,m)}\lesssim\Vert
f\Vert_{L^{p}(X,m)}. \label{SteinMax}%
\end{equation}
According to Theorem \ref{ThmMaximal} and \eqref{SteinMax}, we derive
\begin{equation}
\lim_{\lambda\rightarrow\infty}\lambda^{p}\iint_{E_{\lambda,\gamma/p}%
}t^{\gamma-1}\,\mathrm{d}t\,\mathrm{d}m(x) = \frac{1}{\gamma} \Vert f\Vert_{L^{p}(X,m)}^{p} \asymp\sup_{\lambda> 0}\lambda
^{p}\iint_{E_{\lambda,\gamma/p}}t^{\gamma-1}\,\mathrm{d}t\,\mathrm{d}%
m(x),\quad\gamma>0, \label{SteinMax2*}%
\end{equation}
where
\[
E_{\lambda,\gamma/p}=\{(x,t)\in X\times(0,1):|T_{\floor{1/t}}f(x)|>\lambda
\,t^{\gamma/p}\}.
\]
Choosing $\gamma=1$,
\begin{align}
\Vert f\Vert_{L^{p}(X,m)}^{p}  &  =\lim_{\lambda\rightarrow\infty}%
\lambda^{p}(m\times\mathcal{L})\bigg(\bigg\{(x,t)\in X\times(0,1):\frac
{|T_{\floor{1/t}}f(x)|}{t^{1/p}}>\lambda\bigg\}\bigg)\nonumber\\
&  \asymp\sup_{\lambda> 0}\lambda^{p}(m\times\mathcal{L})\bigg(\bigg\{(x,t)\in
X\times(0,1):\frac{|T_{\floor{1/t}}f(x)|}{t^{1/p}}>\lambda\bigg\}\bigg).
\label{SteinMax2}%
\end{align}

In view of the celebrated theorem of Carleson on pointwise a.e. convergence of
Fourier series, the partial sums $S_{n}f(x)=\sum_{|k|\leq n}\widehat{f}%
(k)e^{ikx}$ on $L^{p}({\mathbb{T}}),$ can be seen to satisfy the previous
assumptions. In fact, working with $S_{n}$, the characterizations
\eqref{SteinMax2*} and \eqref{SteinMax2} can be extended to cover all
$p\in(1,\infty)$ via the maximal inequality related to the Carleson--Hunt
theorem, that is,
\[
\Big\|\sup_{n\in{{\mathbb{N}}}}|S_{n}f|\Big\|_{L^{p}({\mathbb{T}})}%
\lesssim\Vert f\Vert_{L^{p}({\mathbb{T}})},\quad p\in(1,\infty).
\]
In particular, the following holds
\begin{align*}
\Vert f\Vert_{L^{p}({\mathbb{T}})}^{p}  &  =\lim_{\lambda\rightarrow
\infty}\lambda^{p}\mathcal{L}^{2}\bigg(\bigg\{(x,t)\in{\mathbb{T}}%
\times(0,1):\frac{|S_{\floor{1/t}}f(x)|}{t^{1/p}}>\lambda\bigg\}\bigg)\\
&  \asymp\sup_{\lambda> 0}\lambda^{p}\mathcal{L}^{2}\bigg(\bigg\{(x,t)\in
{\mathbb{T}}\times(0,1):\frac{|S_{\floor{1/t}}f(x)|}{t^{1/p}}>\lambda
\bigg\}\bigg)
\end{align*}
for $p\in(1,\infty)$.

\subsection{Hilbert Transform}

Let $H$ be the Hilbert transform, defined for $f\in\mathcal{S}({\mathbb{R}})$
by
\[
Hf(x)=\text{p.v.}\int_{{\mathbb{R}}}\frac{f(y)}{x-y}\,\mathrm{d} y.
\]
It is well known that $H$ extends to a bounded map from $L^{p}({\mathbb{R}})$
into itself. For $\varepsilon>0,$ consider the truncated Hilbert transforms
\[
H_{\varepsilon}f(x)=\int_{|x-y|>\varepsilon}\frac{f(y)}{x-y}\,\mathrm{d} y,
\]
and the associated maximal Hilbert transform
\[
H^{\ast}f(x)=\sup_{\varepsilon>0}|H_{\varepsilon}f(x)|.
\]
Note that $H_{\varepsilon}f$ is well-defined for any $f\in L^{p}({\mathbb{R}%
})$ and, obviously, $\lim_{\varepsilon\rightarrow0+}H_{\varepsilon}f(x)=Hf(x)$
for $f\in\mathcal{S}({\mathbb{R}})$. Combining Cotlar's inequality
\[
H^{\ast}f(x)\lesssim M(Hf)(x)+Mf(x)
\]
with the Hardy--Littlewood maximal inequality for $M$, one immediately gets
that $H^{\ast}:L^{p}({\mathbb{R}})\rightarrow L^{p}({\mathbb{R}})$ and so,
$\lim_{\varepsilon\rightarrow0+}H_{\varepsilon}f(x)=Hf(x)$ a.e. $x\in
{\mathbb{R}}$ for all $f\in L^{p}({\mathbb{R}})$. According to Theorem
\ref{ThmMaximal} and a simple variant of Remark \ref{RemarkSharpAssertion}, we have
\begin{align*}
\Vert Hf\Vert_{L^{p}({\mathbb{R}})}^{p}  &  = \lim_{\lambda\rightarrow
\infty}\lambda^{p}\mathcal{L}^{2}\bigg(\bigg\{(x,\varepsilon)\in{\mathbb{R}%
}\times(0,\infty):\frac{|H_{\varepsilon}f(x)|}{\varepsilon^{1/p}}%
>\lambda\bigg\}\bigg)\\
&  \leq\sup_{\lambda> 0 }\lambda^{p}\mathcal{L}^{2}\bigg(\bigg\{(x,\varepsilon
)\in{\mathbb{R}}\times(0,\infty):\frac{|H_{\varepsilon}f(x)|}{\varepsilon
^{1/p}}>\lambda\bigg\}\bigg)\\
&  \leq\Vert H^{\ast}f\Vert_{L^{p}({\mathbb{R}})}^{p}%
\end{align*}
and, in particular,
\[
\Vert Hf\Vert_{L^{p}({\mathbb{R}})}^{p} = \lim_{\lambda\rightarrow\infty
}\lambda^{p}\mathcal{L}^{2}\bigg(\bigg\{(x,\varepsilon)\in{\mathbb{R}}%
\times(0,\infty):\frac{|H_{\varepsilon}f(x)|}{\varepsilon^{1/p}}%
>\lambda\bigg\}\bigg)\lesssim\Vert f\Vert_{L^{p}({\mathbb{R}})}^{p}.
\]

The discussion given above can be further extended to deal with
Calder\'{o}n--Zygmund operators on nonhomogeneous spaces (cf. \cite{Nazarov}),
but we will not go into further details here.

\section{Weak-type inequalities of Brezis--Van Schaftingen--Yung type for
Calder\'{o}n--Campanato spaces\label{sec:Ca-Ca-section}}

\subsection{Embedding of Calder\'{o}n--Campanato spaces into Brezis--Van
Schaftingen--Yung spaces}

\label{CC}

In this section we compare the Calder\'{o}n--Campanato spaces $C^{s}%
_{p}({\mathbb{R}}^{N})$ (cf. \eqref{DefCC}, \eqref{DefCC2}) with the
Brezis--Van Schaftingen--Yung spaces $BSY_{p}^{s}({\mathbb{R}}^{N})$ (cf.
\eqref{pppp}). Let us start with some general considerations that will also
serve as a motivation.

In this section we are generally interested in the largest \textquotedblleft
admissible" function space $\mathcal{X}({\mathbb{R}}^{N})$ for which we have
inequalities of the form
\begin{equation}
\Vert f\Vert_{BSY_{p}^{s}({\mathbb{R}}^{N})}\lesssim\Vert f\Vert
_{\mathcal{X}({\mathbb{R}}^{N})},\quad1<p<\infty,\quad s\in(0,1).\label{Frac3}%
\end{equation}
The Gagliardo space $W^{s,p}({\mathbb{R}}^{N})$ (cf. \eqref{DefG}) is
naturally part of the competition. Indeed, since $L^{p}({\mathbb{R}}^{N}%
\times{\mathbb{R}}^{N})\subset L(p,\infty)({\mathbb{R}}^{N}\times{\mathbb{R}%
}^{N}),$ we have
\begin{equation}
\Vert f\Vert_{BSY_{p}^{s}({\mathbb{R}}^{N})}\lesssim\bigg\|\frac
{f(x)-f(y)}{|x-y|^{\frac{N}{p}+s}}\bigg\|_{L^{p}({\mathbb{R}}^{N}%
\times{\mathbb{R}}^{N})}=\Vert f\Vert_{W^{s,p}({\mathbb{R}}^{N})}%
.\label{fracX}%
\end{equation}
In other words, (\ref{Frac3}) is trivially satisfied for the space
$\mathcal{X}({\mathbb{R}}^{N})=W^{s,p}({\mathbb{R}}^{N}).$ But can we do
better? Informally, one reason to believe so is that letting $s\rightarrow1$,
then while the left hand side tends $\Vert\nabla f\Vert_{L^{p}({\mathbb{R}%
}^{N})}$ (cf. \cite{Brezis})$,$ we need mitigating constants in order for
$\Vert f\Vert_{W^{s,p}({\mathbb{R}}^{N})}\rightarrow\Vert\nabla f\Vert
_{L^{p}({\mathbb{R}}^{N})}$ (cf. \cite{Bourgain00})$.$ To sharpen
(\ref{fracX}) we shall analyze and modify the methods of \cite{Brezis}, in
order to be able to incorporate the fractional cases to the analysis. So let
us start by recalling that the proof of
\begin{equation}
\Vert f\Vert_{BSY_{p}^{1}({\mathbb{R}}^{N})}\ =\bigg\|\frac{f(x)-f(y)}%
{|x-y|^{\frac{N}{p}+1}}\bigg\|_{L(p,\infty)({\mathbb{R}}^{N}\times{\mathbb{R}%
}^{N})}\lesssim\Vert\nabla f\Vert_{W^{1,p}({\mathbb{R}}^{N})},\quad
1<p<\infty,\label{Frac4}%
\end{equation}
given in \cite{Brezis} relies on the classical estimate 
\begin{equation}
|f(x)-f(y)|\lesssim|x-y|(M(\nabla f)(x)+M(\nabla f)(y)),\label{Lusin}%
\end{equation}
combined with the $L^{p}$ maximal inequality for the Hardy--Littlewood maximal
function $M$. The inequality \eqref{Lusin} has a long history which goes back to the famous work of John and Nirenberg \cite{JN61} on $\text{BMO}$ and $\text{VMO}$ functions (see also \cite{Liu}
and \cite{Bojarski}). To proceed in analogous way with the fractional case we could
replace (\ref{Lusin}) by  
\begin{equation}
|f(x)-f(y)|\lesssim|x-y|^{s}(M_{1-s}(\nabla f)(x)+M_{1-s}(\nabla f)(y)),\quad
s\in(0,1],\label{Hedberg}%
\end{equation}
where $M_{1-s}f$ is the fractional maximal operator defined by
\[
M_{1-s}f(x)=\sup_{r>0}\frac{r^{1-s}}{\mathcal{L}^{N}(B(x,r))}\int%
_{B(x,r)}|f(y)|\,\mathrm{d}y.
\]
The inequality \eqref{Hedberg} has its roots in the analysis of Campanato \cite{Campanato} (cf. also \cite{Hedberg}). The use of $M_{1-s}$ is consistent with our aims since letting $s=1$ in
\eqref{Hedberg} yields back \eqref{Lusin}. Having at our disposal
\eqref{Hedberg}, we can write
\[
\frac{|f(x)-f(y)|}{|x-y|^{s}|x-y|^{\frac{N}{p}}}\lesssim\frac{1}%
{|x-y|^{\frac{N}{p}}}(M_{1-s}(\nabla f)(x)+M_{1-s}(\nabla f)(y)),
\]
and, therefore by the argument in the Introduction related to the proof of
Theorem \ref{ThmMaximal}(i), we find
\begin{equation}
\bigg\|\frac{f(x)-f(y)}{|x-y|^{\frac{N}{p}+s}}\bigg\|_{L(p,\infty
)({\mathbb{R}}^{N}\times{\mathbb{R}}^{N})}\lesssim\Vert M_{1-s}(\nabla
f)\Vert_{L^{p}({\mathbb{R}}^{N})},\quad s\in(0,1].\label{Frac5}%
\end{equation}
When $s=1$, it follows from the maximal theorem of Hardy--Littlewood that
\begin{equation}
\Vert M(\nabla f)\Vert_{L^{p}({\mathbb{R}}^{N})}\asymp\Vert\nabla
f\Vert_{L^{p}({\mathbb{R}}^{N})},\quad1<p<\infty,\label{Frac6}%
\end{equation}
which implies \eqref{Frac4}. However, the equivalence \eqref{Frac6} is no
longer true when $M$ is replaced by its fractional counterpart $M_{1-s}%
,\,s\in(0,1)$, and we can only expect the one-sided estimate that is provided
by the Hardy--Littlewood--Sobolev inequality. More precisely, for $s\in(0,1],$
we have that for $\frac{1}{q}=\frac{1-s}{N}+\frac{1}{p}$
\begin{equation}
\Vert M_{1-s}(\nabla f)\Vert_{L^{p}({\mathbb{R}}^{N})}\lesssim\Vert\nabla
f\Vert_{L^{q}({\mathbb{R}}^{N})}.\label{Frac7}%
\end{equation}
Thus, combining \eqref{Frac5} and \eqref{Frac7}, we have
\[
\bigg\|\frac{f(x)-f(y)}{|x-y|^{\frac{N}{p}+s}}\bigg\|_{L(p,\infty
)({\mathbb{R}}^{N}\times{\mathbb{R}}^{N})}\lesssim\Vert\nabla f\Vert
_{L^{q}({\mathbb{R}}^{N})}.
\]
However, this estimate is far from being optimal and it is even weaker than
(\ref{fracX}) since $\mathring{W}_{q}^{1}({\mathbb{R}}^{N})\subset
W^{s,p}({\mathbb{R}}^{N}).$

From this point of view to proceed further we need a sharper form of
(\ref{Hedberg}). Fortunately, by Campanato \cite{Campanato} (see also \cite{DeVore}) we have,%
\begin{equation}
|f(x)-f(y)|\lesssim|x-y|^{s}(f_{s}^{\#}(x)+f_{s}^{\#}(y))\quad\text{a.e.}\quad
x,y\in{\mathbb{R}}^{N}, \label{Calderon}%
\end{equation}
where $f_{s}^{\#}$ is given in \eqref{DefCC2}.
One can readily verify that (\ref{Calderon}) indeed sharpens (\ref{Hedberg})
by means of applying Poincar\'{e}'s inequality to obtain
\[
f_{s}^{\#}(x)\lesssim M_{1-s}(\nabla f)(x),\quad0<s\leq1.
\]
Therefore,%
\[
\left\Vert f\right\Vert _{C_{p}^{s}({\mathbb{R}}^{N})}=\left\Vert f_{s}%
^{\#}\right\Vert _{L^{p}({\mathbb{R}}^{N})}\lesssim\left\Vert M_{1-s}(\nabla
f)\right\Vert _{L^{p}({\mathbb{R}}^{N})}\lesssim\Vert\nabla f\Vert
_{L^{q}({\mathbb{R}}^{N})},\quad\frac{1}{q}=\frac{1-s}{N}+\frac{1}{p}.
\]

The preceding discussion was our motivation to study relations between
$C^{s}_{p}({\mathbb{R}}^{N})$ and $BSY^{s}_{p}({\mathbb{R}}^{N})$.


Now, starting from (\ref{Calderon}) instead of \eqref{Hedberg}, the argument
given above, verbatim, yields
\begin{equation}
\label{ppppppp}\bigg\|\frac{f(x)-f(y)}{|x-y|^{\frac{N}{p}+s}}%
\bigg\|_{L(p,\infty)({\mathbb{R}}^{N}\times{\mathbb{R}}^{N})} \lesssim\Vert
f_{s}^{\#} \Vert_{L^{p}({\mathbb{R}}^{N})} =\left\Vert f\right\Vert
_{C_{p}^{s}({\mathbb{R}}^{N})}.
\end{equation}
In other words we have obtained

\begin{theorem}
\label{ThmMarkao} Assume $1 < p < \infty$ and $s \in(0,1]$. Then $C_{p}%
^{s}({\mathbb{R}}^{N})\subset BSY_{p}^{s}({\mathbb{R}}^{N}).$
\end{theorem}

\begin{remark}
The result is not trivial in the sense that for $s\in(0,1)$ we have
\[
W^{s,p}({\mathbb{R}}^{N})\subsetneq C_{p}^{s}({\mathbb{R}}^{N}),
\]
cf. \cite[Section 7]{DeVore}. The interrelations between the spaces
$BSY_{p}^{s}({\mathbb{R}}^{N}),\,W^{s,p}({\mathbb{R}}^{N})$ and $C_{p}%
^{s}({\mathbb{R}}^{N})$ are illustrated in Figure 1 below.

\bigskip

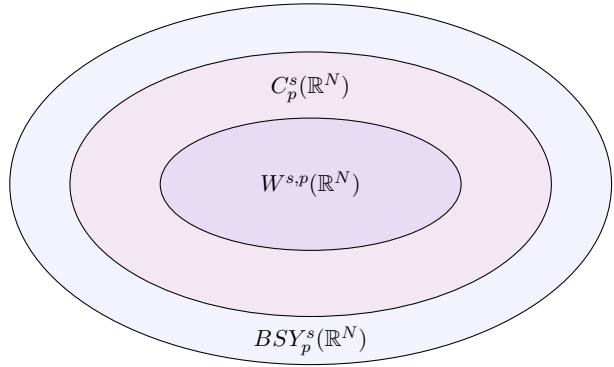
\begin{center}
\begin{tikzpicture}[fill opacity=0.05, ,xscale=0.8,yscale=0.8]
\fill[blue] (0,0) ellipse (5.0 and 3.0);
\draw[fill opacity=0.2] (0,0) ellipse (5.0 and 3.0);
\fill[red] (0,0) ellipse (4.0 and 2.2);
\draw[fill opacity=0.2] (0,0) ellipse (4.0 and 2.2);
\fill[blue] (0,0) ellipse (2.5 and 1.1);
\draw[fill opacity=0.2](0,0) ellipse (2.5 and 1.1);
\node[fill opacity=2,xscale=0.8,yscale=0.8] at (0,1.6) {$C^s_p(\R^N)$};
\node[fill opacity=2,xscale=0.8,yscale=0.8] at (0,0) {$W^{s,p}(\R^N)$};
\node[fill opacity=2,xscale=0.8,yscale=0.8] at (0,-2.6) {$BSY^s_p(\R^N)$};
\end{tikzpicture}
\end{center}

\phantom{qqq}

\bigskip

{\small \textbf{Fig. 1:} Relationships between Brezis--Van Schaftingen--Yung
spaces, Sobolev spaces and Calder\'{o}n--Campanato spaces.}

\bigskip

On the other hand, in the integer case $s=1$, the inequality provided in
Theorem \ref{ThmMarkao} (i.e., \eqref{ppppppp}) coincides with \eqref{Frac4}
since $C_{p}^{1}({\mathbb{R}}^{N})=W^{1,p}({\mathbb{R}}^{N})=BSY_{p}%
^{1}({\mathbb{R}}^{N})$ with equivalence of seminorms (cf. \cite[Theorem
6.2]{DeVore} for the first equivalence and \eqref{I1} for the second one.)
\end{remark}

In the next section we shall show, using Theorem \ref{ThmMaximal}, that the
improvement provided by Theorem \ref{ThmMarkao} is nearly best possible.

\subsection{Calder\'{o}n--Campanato spaces characterized via the
Fefferman--Stein inequality and Theorem \ref{ThmMaximal}\label{ExampleC}}

In this section we give a new characterization of Calder\'{o}n--Campanato
spaces \`{a} la Brezis--Van Schaftingen--Yung.

Let us recall the definition of the local sharp fractional maximal function
(compare with \eqref{DefCC2} above)
\begin{equation}
f_{R, s}^{\#}(x)=\sup_{0<r<R}\frac{1}{r^{s+N}}\int_{B(x,r)}|f(y)-
(f)_{B(x,r)}|\,\mathrm{d}y,\quad x\in{\mathbb{R}}^{N},\quad R>0, \quad s
\in[0,1]. \label{LocalSharpFracMax}%
\end{equation}
In particular, if $s=0$ then we get the restricted versions $f_{R}^{\#}$ of
the Fefferman--Stein maximal function $f^{\#}$.

We start by observing that pointwise,
\[
f_{s}^{\#}(x)=\lim_{R\rightarrow0+}f_{\frac{1}{R},s}^{\#}(x)=\sup
_{R>0}f_{\frac{1}{R},s}^{\#}(x).
\]
Therefore, applying Theorem \ref{ThmMaximal}, we obtain
\begin{equation}
\lim_{\lambda\rightarrow\infty}\lambda^{p}\iint_{E_{\lambda,\gamma/p}%
}R^{\gamma-1}\,\mathrm{d}R\,\mathrm{d}x = \sup_{\lambda> 0}\lambda^{p}%
\iint_{E_{\lambda,\gamma/p}}R^{\gamma-1}\,\mathrm{d}R\,\mathrm{d}x = \frac
{1}{\gamma} \Vert f_{s}^{\#}\Vert_{L^{p}({\mathbb{R}}^{N})}^{p},\quad\gamma>0,
\label{FSx}%
\end{equation}
where
\[
E_{\lambda,\gamma/p}=\{(x,R)\in{\mathbb{R}}^{N}\times(0,\infty):|f_{s,1/R}%
^{\#}(x)|>\lambda\,R^{\gamma/p}\}.
\]
In particular, setting $\gamma=1$,
\begin{align}
\Vert f_{s}^{\#}\Vert_{L^{p}({\mathbb{R}}^{N})}^{p}  &  = \lim_{\lambda
\rightarrow\infty}\lambda^{p}\mathcal{L}^{N+1}\bigg(\bigg\{(x,R)
\in{\mathbb{R}}^{N}\times(0,\infty):\frac{|f_{s,1/R}^{\#}(x)|}{R^{1/p}%
}>\lambda\bigg\}\bigg)\nonumber\\
&  = \sup_{\lambda> 0}\lambda^{p}\mathcal{L}^{N+1}\bigg(\bigg\{(x,R)
\in{\mathbb{R}}^{N}\times(0,\infty):\frac{|f_{s,1/R}^{\#}(x)|}{R^{1/p}%
}>\lambda\bigg\}\bigg). \label{FS}%
\end{align}
Consequently, \eqref{FSx} and \eqref{FS} can be rewritten in terms of $\Vert
f\Vert_{C_{p}^{s}({\mathbb{R}}^{N})}$ (cf. \eqref{DefCC}). For instance,
(\ref{FS}) gives

\begin{theorem}
\label{ThmFSNew} Let $s\in[0,1]$ and $1 <p < \infty$. Then
\begin{align*}
\Vert f\Vert_{C^{s}_{p}(\mathbb{R}^{N})}^{p}  &  = \lim_{\lambda
\rightarrow\infty}\lambda^{p}\mathcal{L}^{N+1}\bigg(\bigg\{(x,R)
\in{\mathbb{R}}^{N}\times(0,\infty):\frac{|f_{s,1/R}^{\#}(x)|}{R^{1/p}
}>\lambda\bigg\}\bigg)\\
&  = \sup_{\lambda> 0}\lambda^{p}\mathcal{L}^{N+1}\bigg(\bigg\{(x,R)
\in{\mathbb{R}}^{N}\times(0,\infty):\frac{|f_{s,1/R}^{\#}(x)|}{R^{1/p}%
}>\lambda\bigg\}\bigg).
\end{align*}

\end{theorem}

Specializing Theorem \ref{ThmFSNew} (or more generally, \eqref{FSx}) with
$s=0$, we apply the Fefferman--Stein inequality (cf. \cite{Fefferman})
\[
\Vert f^{\#}\Vert_{L^{p}({\mathbb{R}}^{N})}\asymp\Vert f\Vert_{L^{p}%
({\mathbb{R}}^{N})}\quad\text{for any}\quad f\in S_{0}({\mathbb{R}}%
^{N})\footnote{$S_{0}({\mathbb{R}}^{N})$ denotes the space of all measurable
functions $f$ on ${\mathbb{R}}^{N}$ such that for any $\lambda>0,\,\mathcal{L}%
^{N}(\{x\in{\mathbb{R}}^{N}:|f(x)|>\lambda\})<\infty.$}%
\]
to give a characterization of Lebesgue norms in terms of restricted
Fefferman--Stein maximal functions. More precisely, we have

\begin{theorem}
Let $\gamma>0, \, 1 < p < \infty,$ and%
\[
E_{\lambda,\gamma/p}=\{(x,R)\in{\mathbb{R}}^{N}\times(0,\infty):|f_{1/R}%
^{\#}(x)|>\lambda\,R^{\gamma/p}\} \quad\text{for each} \quad\lambda> 0.
\]
Then
\[
\lim_{\lambda\rightarrow\infty}\lambda^{p}\iint_{E_{\lambda,\gamma/p}%
}R^{\gamma-1}\,\mathrm{d}R\,\mathrm{d}x = \sup_{\lambda> 0}\lambda^{p}%
\iint_{E_{\lambda,\gamma/p}}R^{\gamma-1}\,\mathrm{d}R\,\mathrm{d}x\asymp\Vert
f\Vert_{L^{p}({\mathbb{R}}^{N})}^{p}.
\]

\end{theorem}

In the same fashion, since $\mathring{W}^{1}_{p}({\mathbb{R}}^{N}) = C^{1}%
_{p}({\mathbb{R}}^{N}), \, 1 < p < \infty,$ (cf. \cite[Theorem 6.2]{DeVore})
we can invoke \eqref{FSx} to obtain an alternative characterization to
\eqref{l1*} in terms of maximal operators.

\begin{theorem}
Let $\gamma>0, \, 1 < p < \infty,$ and%
\[
E_{\lambda,\gamma/p}=\{(x,R)\in{\mathbb{R}}^{N}\times(0,\infty):|f_{1,1/R}%
^{\#}(x)|>\lambda\,R^{\gamma/p}\} \quad\text{for each} \quad\lambda> 0.
\]
Then
\[
\lim_{\lambda\rightarrow\infty}\lambda^{p}\iint_{E_{\lambda,\gamma/p}%
}R^{\gamma-1}\,\mathrm{d}R\,\mathrm{d}x = \sup_{\lambda> 0}\lambda^{p}%
\iint_{E_{\lambda,\gamma/p}}R^{\gamma-1}\,\mathrm{d}R\,\mathrm{d}x\asymp
\Vert\nabla f\Vert_{L^{p}({\mathbb{R}}^{N})}^{p}.
\]
In particular
\begin{align*}
\Vert\nabla f\Vert_{L^{p}(\mathbb{R} ^{N})}^{p}  &  \asymp\lim_{\lambda
\rightarrow\infty} \lambda^{p} \mathcal{L}^{N+1} \bigg(\bigg\{(x,R)
\in{\mathbb{R}}^{N} \times(0,\infty): \frac{|f^{\#}_{1,1/R}(x)|}{R^{1/p}} >
\lambda\bigg\} \bigg)\\
&  = \sup_{\lambda> 0} \lambda^{p} \mathcal{L}^{N+1} \bigg(\bigg\{(x,R)
\in{\mathbb{R}}^{N} \times(0,\infty): \frac{|f^{\#}_{1,1/R}(x)|}{R^{1/p}} >
\lambda\bigg\} \bigg).
\end{align*}

\end{theorem}

\section{Sharp maximal-type operators and Calder\'{o}n--Campanato
spaces\label{sec:max-ca-ca}}

\subsection{Introduction}

As we have seen, through the use of Calder\'{o}n--Campanato spaces, we obtain
sharp fractional embedding results in the spirit of \cite{Brezis}. In this
section we extend the embedding theory of Calder\'{o}n--Campanato spaces on
${\mathbb{R}}^{N}$ to the setting of metric spaces. For this purpose we need
to find a suitable replacement of (\ref{Calderon}). The crux of the matter is
the introduction of generalized sharp maximal operators which we use to obtain
a substitute for (\ref{Calderon}). We believe that the new inequalities as
well as the methods of proof, are of independent interest.

Let $(X,d,m)$ be a metric measure space. Let $\rho:(0,\infty)\rightarrow
(0,\infty)$ and $R\in(0,\text{diam}(X))$. For locally integrable functions
$f:X\rightarrow\mathbb{R}$, we introduce the \emph{sharp maximal-type
operator}
\[
f_{\rho}^{\#}(x)=\sup_{0<r<\text{diam}(X)}\frac{m(B(x,r))}{\rho(r)}%
\int_{B(x,r)}|f(y)-(f)_{B(x,r)}|\,\mathrm{d}m(y),
\]
where $(f)_{B(x,r)}=\frac{1}{m(B(x,r))}\int_{B(x,r)}f(y) \mathrm{d} m(y).$
Furthermore, the restricted maximal operator $f_{R,\rho}^{\#}$ is defined by
\[
f_{R,\rho}^{\#}(x)=\sup_{0<r<R}\frac{m(B(x,r))}{\rho(r)}\int_{B(x,r)}%
|f(y)-(f)_{B(x,r)}|\,\mathrm{d}m(y).
\]
The \emph{Calder\'{o}n--Campanato-type spaces} $C_{p}^{\rho}(X,m)$ are defined
in the usual fashion using the seminorms given by
\[
\Vert f\Vert_{C_{p}^{\rho}(X,m)}=\Vert f_{\rho}^{\#}\Vert_{L^{p}(X,m)}, \quad1
< p < \infty.
\]

Recall that a metric space $X=(X,d)$ endowed with a Borel measure $m$ is said
to be \emph{(Ahlfors) $N$-regular} for some $N\geq0$ (where $N$ is not
necessarily an integer) if there exist constants $c_{0}>0$ and $C_{0}<\infty$
such that
\begin{equation}
c_{0}r^{N}\leq m(B)\leq C_{0}r^{N}, \label{Ahlfors}%
\end{equation}
for every closed ball $B$ in $X$ with radius $r<\text{diam}(X)$. In this
setting, if we let $\rho(r)=r^{2N+s},\,s\in[0,1],$ we obtain the classical
sharp maximal operators $f_{s}^{\#}$ and $f_{R,s}^{\#}$ (cf. \eqref{DefCC2}
and \eqref{LocalSharpFracMax}, respectively, for $X={\mathbb{R}}^{N}$ endowed
with the Lebesgue measure) and the corresponding space $C_{p}^{s}(X,m)$ (cf.
\eqref{DefCC}). On the other hand, setting
\[
\rho_{\beta}(r)=r^{2N}\left\{
\begin{array}
[c]{lcl}%
(\log\frac{e}{r})^{-\beta} & \text{if} & r\in(0,1],\\
1 & \text{if} & r\in(1,\infty),
\end{array}
\right.
\]
we obtain
\begin{equation}
\label{kk}f_{0,\beta}^{\#}(x)=\sup_{0<r<\text{diam}(X)}\frac{r^{N}}%
{\rho_{\beta}(r)}\int_{B(x,r)}|f(y)-(f)_{B(x,r)}|\,\mathrm{d}m(y)
\end{equation}
that corresponds to the logarithmic Calder\'{o}n--Campanato space
$C_{p}^{0,\beta}(X,m)$. Logarithmic-type maximal operators are attracting an
increasing interest in the study of the regularity properties of Lagrangian
flows in Sobolev spaces, cf. \cite{Crippa}, \cite{Brue} and \cite{BrueSemola}.
We will return to this topic in Section \ref{sec:Crippa} below.

\subsection{Comparison of generalized Calder\'{o}n--Campanato spaces and
Brezis--Van Schaftingen--Yung spaces\label{sec:ca-ca}}

In this section we extend the results of the previous sections to the
generalized Calder\'{o}n--Campanato spaces $C_{p}^{\rho}(X,m).$ For this
purpose we need to extend the Brezis--Van Schaftingen--Yung spaces.
Fortunately, the construction of the appropriate spaces is very natural and is
dictated by the underlying maximal inequalities and the mixed norm
inequalities that are available to us.

The main result of this section reads as follows:

\begin{theorem}
\label{ThmCCBSY} Let $(X,d,m)$ be a $N$-regular metric measure space. Let
$\rho:(0,\infty)\rightarrow(0,\infty)$ be continuous and increasing, with
$\lim_{r\rightarrow0+}\rho(r)=0$. Let $\bar{\rho}(r)=\rho(2r), \, r>0.$ Assume
$1 < p < \infty$. Then
\begin{equation}
\left\Vert \frac{f(x)-f(y)}{d(x,y)^{\frac{N}{p}}\Big(\int_{0}^{2d(x,y)}%
\frac{\mathrm{d}\bar{\rho}(\lambda)}{\lambda^{2N}}\Big)}\right\Vert
_{L(p,\infty)(X\times X,m\times m)}\lesssim\Vert f\Vert_{C_{p}^{\bar{\rho}%
}(X,m)}. \label{agregada}%
\end{equation}

\end{theorem}

\begin{remark}
In the model examples given below it holds that $\rho\asymp\bar{\rho}$.
\end{remark}

Specializing Theorem \ref{ThmCCBSY} letting $\rho(r)=r^{2N+s},\,s\in(0,1],$ yields

\begin{corollary}
\label{CorCCBSY2} Suppose that $(X,d,m)$ is $N$-regular, and let $s\in(0,1]$
and $1 < p < \infty$. Then
\begin{equation}
\left\Vert \frac{f(x)-f(y)}{d(x,y)^{\frac{N}{p}+s}}\right\Vert _{L(p,\infty
)(X\times X,m\times m)}\lesssim\Vert f\Vert_{C_{p}^{s}(X,m)}. \label{CorCCBSY}%
\end{equation}

\end{corollary}

\begin{remark}
In view of the previous corollary we see that the functional defined by the
left-hand side of (\ref{agregada}) can be conceived as a suitable
generalization of the $BSY^{s}_{p}$ functional \eqref{pppp}.
\end{remark}

The limiting case $s=0$ in Corollary \ref{CorCCBSY2} reads as follows:

\begin{corollary}
Suppose that $(X,d,m)$ is $N$-regular and let $\beta>1$ and $1 < p < \infty$.
Let
\begin{equation}
\rho_{\beta}(r)=r^{2N}\left\{
\begin{array}
[c]{lcl}%
(\log\frac{e}{r})^{-\beta} & \text{if} & r\in(0,1],\\
1 & \text{if} & r\in(1,\infty),
\end{array}
\right.  \label{shift}%
\end{equation}
and
\begin{equation}
w_{\beta}(r)=\left\{
\begin{array}
[c]{lcl}%
(\log\frac{e}{r})^{-\beta+1} & \text{if} & r\in(0,1],\\
\log r & \text{if} & r\in(1,\infty).
\end{array}
\right.  \label{shift2}%
\end{equation}
Then
\[
\left\Vert \frac{f(x)-f(y)}{d(x,y)^{\frac{N}{p}}w_{\beta}(d(x,y))}\right\Vert
_{L(p,\infty)(X\times X,m\times m)}\lesssim\Vert f\Vert_{C_{p}^{0,\beta}%
(X,m)}.
\]

\end{corollary}

\begin{remark}
The limiting case $s=0$ given in the previous result shows a shift in the
logarithmic smoothness of the involved spaces (see \eqref{shift} and
\eqref{shift2}). This is in sharp contrast with the non-limiting case $s
\in(0,1]$ where both spaces involved in \eqref{CorCCBSY} have smoothness $s$.
\end{remark}

The proof of Theorem \ref{ThmCCBSY} depends on two tools, which are of
independent interest, namely, a new Garsia-type inequality\footnote{Closely
related results can be found in the work by Preston \cite{preston}, although
he does not formulate the results in terms of sharp maximal functions.} which
extends \eqref{Calderon} to functions of generalized smoothness in metric
spaces (cf. Proposition \ref{PropGarsia} below) combined, as usual in this
paper, with mixed-norm inequalities (cf. Proposition \ref{PropMixed} below.)

\begin{proposition}
[Garsia-type inequality]\label{PropGarsia} Let $(X,d,m)$ be a metric measure
space for which the Lebesgue's differentiation theorem holds. Suppose that
$\mu(r):=\inf_{x\in X}m(B(x,r))>0$ for every $r > 0$\footnote{This condition
may be satisfied by non doubling measures on ${\mathbb{R}}^{N}.$ For example,
the measure $\mathrm{d} m(x)=e^{\left\vert x\right\vert } \mathrm{d} x$ on the
line satisfies the condition. On the other hand, for the Gausssian measure
$\gamma_{N}$ on ${\mathbb{R}}^{N}$ we have $\gamma_{N}(B(x,r))\leq\frac
{\omega_{N-1}}{2\pi^{N/2}}r^{N}e^{2r\left\vert x\right\vert }e^{-\left\vert
x\right\vert ^{2}}$ where $\omega_{N-1}$ is the (surface) measure of
$\mathbb{S}^{N-1}$ (cf. \cite[Lemma 1.2, pg. 5]{urbina}). In particular, the
condition is not satisfied in this case.
\par
{}}. Let $\rho:(0,\infty)\rightarrow(0,\infty)$ be continuous and increasing
with $\lim_{r\rightarrow0^{+}}\rho(r)=0$, and set $\bar{\rho}(r)=\rho(2 r), \,
r > 0$. Then
\begin{equation}
|f(x)-f(y)|\leq9 \Big(\int_{0}^{2d(x,y)}\frac{\mathrm{d}\bar{\rho}(\lambda
)}{\mu(\lambda)^{2}}\Big)(f_{2d(x,y),\bar{\rho}}^{\#}(x)+f_{2d(x,y),\bar{\rho
}}^{\#}(y)) \label{G1}%
\end{equation}
for almost every $x,y\in X$.
\end{proposition}

The proof of Proposition \ref{PropGarsia} will be postponed to Section
\ref{Subsection:Garsia}.

\begin{remark}
\begin{enumerate}
\item The reason we attribute \eqref{G1} to Garsia will be made clear in the
process of its proof.

\item The preceding result can be applied to the so-called \emph{Vitali
spaces} (cf. \cite[pg. 6]{Heinonen}). This class of spaces includes not only
metric measure spaces satisfying the doubling condition, but also
${\mathbb{R}}^{N}$ equipped with a Radon measure.
\end{enumerate}
\end{remark}

\begin{example}
Recall that a measurable function $b:(0,\infty)\rightarrow(0,\infty)$ is said
to be \emph{slowly varying} if
\[
\lim_{r\rightarrow\infty}\frac{b(rv)}{b(r)}=1\quad\text{for all}\quad v>0;
\]
see \cite{Bingham}. Special cases of slowly varying functions include powers
of logarithms, iterated logarithms, \textquotedblleft broken\textquotedblright%
\ logarithms defined as
\[
b(r)=\left\{
\begin{array}
[c]{lcl}%
(1-\log r)^{\alpha} & \text{ if } & r\in(0,1],\\
(1+\log r)^{\beta} & \text{ if } & r\in(1,\infty),
\end{array}
\right.
\]
where $\alpha,\beta\in{\mathbb{R}}$ and the family of functions $b(r)=\exp
(|\log r|^{\alpha}),\,\alpha\in(0,1)$.

Let $X={\mathbb{R}}^{N}$ equipped with the Lebesgue measure. Assume
$\rho(r)=r^{2N+s} b(r)$ where $s\in(0,1]$ and $b$ is a slowly varying
function. Applying Proposition \ref{PropGarsia}, we have
\begin{equation}
|f(x)-f(y)|\lesssim|x-y|^{s}b(|x-y|)(f_{2|x-y|,s,b}^{\#}(x)+f_{2|x-y|,s,b}%
^{\#}(y)) \label{G2}%
\end{equation}
for a.e. $x,y\in{\mathbb{R}}^{N}$, where
\[
f_{|x-y|,s,b}^{\#}(x)=\sup_{0<r<|x-y|}\frac{1}{r^{N+s}b(r)}\int_{B(x,r)}%
|f(y)-(f)_{B(x,r)}|\,\mathrm{d}y.
\]
In particular, setting $b\equiv1$ in \eqref{G2} we recover \eqref{Calderon},
i.e.,
\begin{equation}
|f(x)-f(y)|\leq C_{s,N}|x-y|^{s}(f_{s}^{\#}(x)+f_{s}^{\#}(y)) \label{KM}%
\end{equation}
where $C_{s,N}$ is a positive constant depending only on $s$ and $N$. In fact,
a perusal of the proof of \eqref{KM} given in \cite[Theorem 2.5]{DeVore} shows
that $C_{s,N}=c_{N}/s$ where $c_{N}$ depends only on $N$. Since $C_{s,N}%
\rightarrow\infty$ as $s\rightarrow0+$, one expects that a new phenomenon will
appear as a counterpart of \eqref{KM} with $s=0$.

To deal with the limiting case $s=0$ in \eqref{KM}, take $\rho_{\beta
}(r)=r^{2N}(-\log r)^{-\beta},\,\beta>1, \, r \in(0,1)$, in \eqref{G1}. Then
\begin{equation}
|f(x)-f(y)|\lesssim(-\log|x-y|)^{-\beta+1}(f_{2|x-y|,0,\beta}^{\#}%
(x)+f_{2|x-y|,0,\beta}^{\#}(y)) \label{G3}%
\end{equation}
whenever $0<|x-y|<1$. Here, $f_{2|x-y|,0,\beta}^{\#}(x)$ is the restricted
version of \eqref{kk}. This inequality shows an interesting phenomenon,
namely, the optimal blow-up related to the logarithmic H\"{o}lder continuity
$(-\log r)^{-\beta+1}$ is given by logarithmic majorants of the ball averages
of oscillations of order $O((-\log r)^{-\beta})$. This phenomenon is not
observed in the classical setting (cf. \eqref{KM}) where the order of
H\"{o}lder continuity and the sharp maximal function are the same.
\end{example}

\begin{example}
Suppose that $(X,m)$ is doubling, i.e., there exists $C>0$ such that
\[
0<m(B(x,2r))\leq C\,m(B(x,r))<\infty\quad\text{for all}\quad x\in
X\quad\text{and}\quad r>0.
\]
Then it is not hard to show that $\mu(r)\gtrsim r^{\log_{2}C}$. Applying
Proposition \ref{PropGarsia} with $\rho(r)=r^{2\log_{2}C+s}$, we have
\begin{equation}
|f(x)-f(y)|\lesssim d(x,y)^{s}(f_{2d(x,y),s}^{\#}(x)+f_{2d(x,y),s}^{\#}(y)),
\label{DSMetric}%
\end{equation}
(cf. \cite[Lemma 3.6]{Hajlasz} and \cite[Lemma 6.2]{Mastylo}). Furthermore, in
the same fashion as above, one can obtain the analogues of \eqref{G2} and
\eqref{G3} in the metric setting.
\end{example}

The second ingredient that we needed for the proof of Theorem \ref{ThmCCBSY}
is the following mixed-norm inequality.

\begin{proposition}
\label{PropMixed} Let $(X_{1},m_{1})$ and $(X_{2},m_{2})$ be $\sigma$-finite
measure spaces and let $1 \leq p < \infty$. Then
\[
\|F(x_{1},x_{2})\|_{L(p,\infty)(X_{1} \times X_{2}, m_{1} \times m_{2})}
\leq\| \|F_{x_{1}}(x_{2})\|_{L(p,\infty)(X_{2}, m_{2})}\|_{L^{p}(X_{1},m_{1})}%
\]
where $F_{x_{1}}(x_{2}) = F(x_{1},x_{2})$ for $x_{1} \in X_{1}$ and $x_{2} \in
X_{2}$.
\end{proposition}

The proof of Proposition \ref{PropMixed} will be postponed until Section
\ref{Subsection:Mixed}.

Assuming the validity of Propositions \ref{PropGarsia} and \ref{PropMixed}, we
proceed to prove Theorem \ref{ThmCCBSY}.

\begin{proof}
[Proof of Theorem \ref{ThmCCBSY}]According to \eqref{G1} we have
\begin{equation}
\frac{|f(x)-f(y)|}{d(x,y)^{\frac{N}{p}}\Big(\int_{0}^{2d(x,y)}\frac
{\mathrm{d}\bar{\rho}(\lambda)}{\lambda^{2N}}\Big)}\lesssim\frac{f_{\bar{\rho
}}^{\#}(x)+f_{\bar{\rho}}^{\#}(y)}{d(x,y)^{\frac{N}{p}}}. \label{ThmCCBSY1}%
\end{equation}
By the lattice property of $L(p,\infty)$ and symmetry we obtain
\begin{equation}
\left\Vert \frac{f(x)-f(y)}{d(x,y)^{\frac{N}{p}}\Big(\int_{0}^{2d(x,y)}%
\frac{\mathrm{d}\bar{\rho}(\lambda)}{\lambda^{2N}}\Big)}\right\Vert
_{L(p,\infty)(X\times X,m\times m)}\lesssim\left\Vert \frac{f_{\bar{\rho}%
}^{\#}(x)}{d(x,y)^{\frac{N}{p}}}\right\Vert _{L(p,\infty)(X\times X,m\times
m)}. \label{ThmCCBSY2}%
\end{equation}
Invoking Proposition \ref{PropMixed}, we can estimate the right-hand side of
\eqref{ThmCCBSY2} as follows
\[
\left\Vert \frac{f_{\bar{\rho}}^{\#}(x)}{d(x,y)^{\frac{N}{p}}}\right\Vert
_{L(p,\infty)(X\times X,m\times m)}\leq\Vert f_{\bar{\rho}}^{\#}(x)\Vert
d(x,y)^{-\frac{N}{p}}\Vert_{L(p,\infty)(X,m)}\Vert_{L^{p}(X,m)}.
\]
Further, basic computations lead to
\begin{equation}
\Vert d(x,y)^{-\frac{N}{p}}\Vert_{L(p,\infty)(X,m)}\asymp1\quad\text{uniformly
w.r.t}\quad x\in X \label{ThmCCBSY3}%
\end{equation}
(see \eqref{Ahlfors}) and therefore
\[
\left\Vert \frac{f_{\bar{\rho}}^{\#}(x)}{d(x,y)^{\frac{N}{p}}}\right\Vert
_{L(p,\infty)(X\times X,m\times m)}\lesssim\Vert f\Vert_{C_{p}^{\bar{\rho}%
}(X,m)}.
\]
Inserting this estimate into \eqref{ThmCCBSY2} we arrive at
\[
\left\Vert \frac{f(x)-f(y)}{d(x,y)^{\frac{N}{p}}\Big(\int_{0}^{2d(x,y)}%
\frac{\mathrm{d}\bar{\rho}(\lambda)}{\lambda^{2N}}\Big)}\right\Vert
_{L(p,\infty)(X\times X,m\times m)}\lesssim\Vert f\Vert_{C_{p}^{\bar{\rho}%
}(X,m)}.
\]

\end{proof}

\begin{remark}
The estimate \eqref{ThmCCBSY3} illustrates the leitmotiv of the method of
Brezis, Van Schaftingen and Yung \cite{Brezis}, \cite{BrezisSY} whose content may be stated as
``going from $L^{p}$ to $L(p,\infty)$". Indeed, note that $L(p,\infty)$ is the
only space within the scale of the Lorentz spaces $L(p,q), \, 0 < q \leq
\infty,$ for which an estimate like \eqref{ThmCCBSY3} holds. In this regard,
see also \cite[Theorem 1]{BrezisSY}.
\end{remark}

\subsection{Proof of Proposition \ref{PropGarsia}}

\label{Subsection:Garsia}

The proofs of \eqref{Calderon} and \eqref{DSMetric} make an essential use of
the doubling property of the underlying space. More precisely, the strategy is
to use $(f)_{B(x,r_{n})},$ $r_{n}=2^{-n}r,\,n\geq0,$ as a regularization of
$f(x)$. Now, since $(f)_{B(x,r)}\asymp(f)_{B(y,r)}$ whenever $r\asymp d(x,y)$,
the problem can thus be reduced to estimate $|f(x)-(f)_{B(x,r)}|$. Clearly,
these strategies seem to fail when $m$ is not doubling. Furthermore, it is not
obvious how we can choose the sequence of the $r_{n}$'s for general weights
$\rho$. We shall overcome some of these obstructions by adapting some ideas
contained in the elegant proof of the main result of Garsia \cite{Garsia72}
(cf. also Preston \cite{preston}).

\begin{proof}
[Proof of Proposition \ref{PropGarsia}]We claim that for all $r>0,$
\begin{equation}
|f(x)-(f)_{B(x,r)}|\leq4\Big(\int_{0}^{r}\frac{\mathrm{d}\bar{\rho}(\lambda
)}{\mu(\lambda)^{2}}\Big)f_{r,\bar{\rho}}^{\#}(x) \label{LemmaGarsia1}%
\end{equation}
for $r>0$ and $m$-a.e. $x\in X$. Indeed, let $x\in X$ be a Lebesgue point and
construct the sequence $(r_{n})_{n\geq0}$ as follows
\begin{equation}
r_{0}=r,\quad\bar{\rho}(r_{n})=\frac{1}{2}\bar{\rho}(r_{n-1}),\quad n\geq1.
\label{Gagli}%
\end{equation}
By the assumptions on $\rho$, it follows that $\lim_{n\rightarrow\infty}%
\bar{\rho}(r_{n})=\lim_{n\rightarrow\infty}\frac{1}{2^{n}}\bar{\rho}(r)=0,$
consequently $r_{n}\downarrow0$ as $n\rightarrow\infty$. Moreover, it is easy
to see from the definitions that
\begin{equation}
\bar{\rho}(r_{n-1})=4(\bar{\rho}(r_{n})-\bar{\rho}(r_{n+1})).
\label{LemmaGarsia2}%
\end{equation}
For each $n\geq1$, we have%
\begin{align*}
|(f)_{B(x,r_{n})}-(f)_{B(x,r_{n-1})}|  &  \leq\frac{1}{m(B(x,r_{n}))}%
\int_{B(x,r_{n})}|f(y)-(f)_{B(x,r_{n-1})}|\,\mathrm{d}m(y)\\
&  \hspace{-2cm} \leq\frac{1}{m(B(x,r_{n}))}\int_{B(x,r_{n-1})}%
|f(y)-(f)_{B(x,r_{n-1})}|\,\mathrm{d}m(y)\text{ (since }r_{n}\text{
decreases)}\\
&  \hspace{-2cm} =\frac{\bar{\rho}(r_{n-1})}{m(B(x,r_{n}))m(B(x,r_{n-1}%
))}\frac{m(B(x,r_{n-1}))}{\bar{\rho}(r_{n-1})}\int_{B(x,r_{n-1})}%
|f(y)-(f)_{B(x,r_{n-1})}|\,\mathrm{d}m(y)\\
&  \hspace{-2cm} \leq\frac{\bar{\rho}(r_{n-1})}{m(B(x,r_{n}))m(B(x,r_{n-1}%
))}f_{r,\bar{\rho}}^{\#}(x)\\
&  \hspace{-2cm} \leq\frac{\bar{\rho}(r_{n-1})}{\left(  m(B(x,r_{n}))\right)
^{2}}f_{r,\bar{\rho}}^{\#}(x) \leq\frac{\bar{\rho}(r_{n-1})}{\mu(r_{n})^{2}%
}f_{r,\bar{\rho}}^{\#}(x).
\end{align*}
Summing the telescoping series yields
\begin{align}
\left\vert f(x)-(f)_{B(x,r)}\right\vert  &  =\limsup_{n\rightarrow\infty
}\left\vert (f)_{B(x,r_{n})}-(f)_{B(x,r)}\right\vert \nonumber\\
&  \leq\sum_{n=1}^{\infty}|(f)_{B(x,r_{n})}-(f)_{B(x,r_{n-1})}|\nonumber\\
&  \leq f_{r,\bar{\rho}}^{\#}(x)\sum_{n=1}^{\infty}\frac{\bar{\rho}(r_{n-1}%
)}{\mu(r_{n})^{2}}. \label{LemmaGarsia4}%
\end{align}
On the other hand, using the monotonicity of $(r_{n})_{n\geq0}$, the fact that
$\mu$ and $\bar{\rho}$ are increasing functions, and (\ref{LemmaGarsia2}), we
obtain
\begin{align}
\int_{0}^{r}\frac{\mathrm{d}\bar{\rho}(\lambda)}{\mu(\lambda)^{2}}  &
=\sum_{n=1}^{\infty}\int_{r_{n}}^{r_{n-1}}\frac{\mathrm{d}\bar{\rho}(\lambda
)}{\mu(\lambda)^{2}}\geq\sum_{n=2}^{\infty}\frac{1}{\mu(r_{n-1})^{2}}%
\int_{r_{n}}^{r_{n-1}}\mathrm{d}\bar{\rho}(\lambda)\nonumber\\
&  =\sum_{n=2}^{\infty}\frac{1}{\mu(r_{n-1})^{2}}(\bar{\rho}(r_{n-1}%
)-\bar{\rho}(r_{n}))\nonumber\\
&  \geq\frac{1}{4}\sum_{n=2}^{\infty}\frac{\bar{\rho}(r_{n-2})}{\mu
(r_{n-1})^{2}}\nonumber\\
&  =\frac{1}{4}\sum_{n=1}^{\infty}\frac{\bar{\rho}(r_{n-1})}{\mu(r_{n})^{2}}.
\label{f}%
\end{align}
Combining \eqref{LemmaGarsia4} and \eqref{f}, we obtain (\ref{LemmaGarsia1}).

Next we show that
\begin{equation}
|f(y)-(f)_{B(x,r)}|\leq5 \Big(\int_{0}^{2r}\frac{\mathrm{d}\bar{\rho}%
(\lambda)}{\mu(\lambda)^{2}}\Big)f_{2r,\bar{\rho}}^{\#}(y),\quad
\text{whenever}\quad y\in B(x,r). \label{LemmaGarsia3}%
\end{equation}
By the previous argument%
\begin{equation}
|f(y)-(f)_{B(y,2r)}|\leq4\Big(\int_{0}^{2r}\frac{\mathrm{d}\bar{\rho}%
(\lambda)}{\mu(\lambda)^{2}}\Big)f_{2r,\bar{\rho}}^{\#}(y). \label{agreg1}%
\end{equation}
Moreover, if $y\in B(x,r),$ then $B(x,r)\subset B(y,2r).$ Therefore
\begin{align}
\left\vert (f)_{B(y,2r)}-(f)_{B(x,r)}\right\vert  &  \leq\frac{1}%
{m(B(x,r))}\int_{B(x,r)}|f(z)-(f)_{B(y,2r)}| \, \mathrm{d}m(z)\nonumber\\
&  \leq\frac{m(B(y,2r))}{m(B(x,r))}\frac{1}{m(B(y,2r))}\int_{B(y,2r)}%
|f(z)-(f)_{B(y,2r)}| \, \mathrm{d}m(z)\nonumber\\
&  \leq\frac{\rho(2r)}{m(B(x,r))^{2}}f_{2r,\bar{\rho}}^{\#}(y)\nonumber\\
&  \leq\frac{\bar{\rho}(r)}{\mu(r)^{2}}f_{2r,\bar{\rho}}^{\#}(y)\nonumber\\
&  \leq f_{2r,\bar{\rho}}^{\#}(y)\int_{0}^{r}\frac{\mathrm{d}\bar{\rho
}(\lambda)}{\mu(\lambda)^{2}}\text{ (since }\bar{\rho}(0)=0,\text{ and }%
\mu\text{ increases).} \label{agreg2}%
\end{align}
Combining (\ref{agreg1}) and (\ref{agreg2}) with the triangle inequality we
obtain (\ref{LemmaGarsia3}).

Finally, given $x,y,$ and $r>d(x,y),$ we write%
\[
|f(x)-f(y)|\leq|f(x)-(f)_{B(x,r)}|+|(f)_{B(x,r)}-f(y)|,
\]
and collecting the estimates \eqref{LemmaGarsia1} and \eqref{LemmaGarsia3} for
each of the terms on the right-hand side we obtain
\[
|f(x)-f(y)|\leq9\Big(\int_{0}^{2r}\frac{\mathrm{d}\bar{\rho}(\lambda
)}{u(\lambda)^{2}}\Big)(f_{2r,\bar{\rho}}^{\#}(x)+f_{2r,\bar{\rho}}^{\#}(y)).
\]
Taking limits on both sides of the previous inequality as $r\rightarrow
d(x,y)+$ we arrive at the desired estimate \eqref{G1}.
\end{proof}

\begin{remark}
One can find similar constructions to \eqref{Gagli} in Gagliardo's pioneering
work on interpolation theory \cite{Gagliardo61, Gagliardo63}. We should
mention that the stopping time method behind these developments were also
implemented successfully to give a proof the Strong Fundamental Lemma of
Interpolation Theory (cf. \cite{Cwikel}.)
\end{remark}

\subsection{Proof of Proposition \ref{PropMixed}}

\label{Subsection:Mixed} Given $t>0$, we have
\begin{align*}
\left\Vert \left\Vert F_{x_{1}}(x_{2})\right\Vert _{L(p,\infty)(X_{2},m_{2}%
)}\right\Vert _{L^{p}(X_{1},m_{1})}^{p}  &  =\int_{X_{1}}\left\Vert F_{x_{1}%
}(x_{2})\right\Vert _{L(p,\infty)(X_{2},m_{2})}^{p}\,\mathrm{d}m_{1}(x_{1})\\
&  =\int_{X_{1}}\sup_{\lambda>0}\,\{\lambda^{p} m_{2} (\{x_{2} \in X_{2} :
|F_{x_{1}}(x_{2})| > \lambda\}) \}\,\mathrm{d}m_{1}(x_{1})\\
&  \geq t^{p}\int_{X_{1}} m_{2} (\{x_{2} \in X_{2} : |F_{x_{1}}(x_{2})| > t\})
\}\,\mathrm{d}m_{1}(x_{1})\\
&  =t^{p} (m_{1} \times m_{2}) (\{(x_{1},x_{2}) \in X_{1} \times X_{2} :
|F(x_{1},x_{2})| > t\}).
\end{align*}
Taking now the supremum over all $t>0$, we arrive at
\begin{align*}
\left\Vert \left\Vert F_{x_{1}}(x_{2})\right\Vert _{L(p,\infty)(X_{2},m_{2}%
)}\right\Vert _{L^{p}(X_{1},m_{1})}^{p}  &  \geq\sup_{t>0}\,t^{p} (m_{1}
\times m_{2}) (\{(x_{1},x_{2}) \in X_{1} \times X_{2} : |F(x_{1},x_{2})| >
t\})\\
&  =\left\Vert F\right\Vert _{L(p,\infty)(X_{1}\times X_{2},m_{1}\times
m_{2})}^{p}.
\end{align*}
\qed

\section{Logarithmic Brezis--Van Schaftingen--Yung spaces via Crippa--De
Lellis estimates\label{sec:Crippa}}

In the previous sections we have investigated in detail the relationships of
the functionals
\begin{equation}
\left\Vert f\right\Vert _{BSY_{p}^{s}(\mathbb{R}^{N})} = \bigg\|\frac
{f(x)-f(y)}{|x-y|^{\frac{N}{p}+s}}\bigg\|_{L(p,\infty)({\mathbb{R}}^{N}%
\times{\mathbb{R}}^{N})} \label{CL1}%
\end{equation}
with Sobolev and Calder\'{o}n--Campanato seminorms. The goal of this section
is to study the corresponding problems when the classical Gagliardo quotient
$\frac{|f(x)-f(y)|}{|x-y|^{s}}$ in \eqref{CL1} is replaced by $\log
\big(1+\frac{|f(x)-f(y)|}{|x-y|^{s}}\big)$.

Firstly, we shall estimate
\[
\bigg\|\log\bigg(1+\frac{|f(x)-f(y)|}{|x-y|^{s}}\bigg)\frac{1}{|x-y|^{\frac
{N}{p}}}\bigg\|_{L(p,\infty)(\mathbb{R}^{N}\times\mathbb{R}^{N})}%
\]
in terms of the Calder\'{o}n--Campanato functionals $\mathcal{N}%
^{s,p}(\mathbb{R}^{N})$ introduced in \eqref{qqq} and \eqref{qqq56}.

\begin{theorem}
\label{BrezisExp} Let $s \in(0,1]$ and $1 < p < \infty$. Then
\begin{equation}
\label{CL*}\bigg\|\log\bigg(1 + \frac{|f(x) - f(y)|}{|x-y|^{s}} \bigg) \frac
{1}{|x-y|^{\frac{N}{p}}} \bigg\|_{L(p,\infty)(\mathbb{R}^{N} \times
\mathbb{R}^{N})} \lesssim\|f\|_{\mathcal{N}^{s,p}(\mathbb{R}^{N})}.
\end{equation}

\end{theorem}

The proof of Theorem \ref{BrezisExp} is a combination of the mixed-norm
approach developed above (see Proposition \ref{PropMixed}) together with
Crippa--De Lellis estimates (cf. \cite{Crippa}). Broadly speaking, these
estimates assert that for certain functions $f$ there exists a nonnegative $g$
such that
\begin{equation}
|f(x)-f(y)|\leq|x-y|^{s}(\exp{\{g(x)+g(y)\}}-1),\quad x,y\in\mathbb{R}^{N}.
\label{Exp}%
\end{equation}
Typically, they are applied to regular Lagrangian flows associated to
time-dependent vector fields, see e.g. \cite[Proposition 2.3]{Crippa},
\cite[(2.1)]{Brue} or \cite[Theorem 3.11]{BrueSemola}.

\begin{proof}
[Proof of Theorem \ref{BrezisExp}]Let $f\in\mathcal{N}^{s,p}(\mathbb{R}^{N})$.
According to \cite[Proposition 2.2]{Brue}, $f$ satisfies \eqref{Exp} for some
nonnegative $g\in L^{p}(\mathbb{R}^{N})$ and $\inf_{g}\Vert g\Vert
_{L^{p}(\mathbb{R}^{N})}\lesssim\Vert f\Vert_{\mathcal{N}^{s,p}(\mathbb{R}%
^{N})}$. Hence, we have
\[
\log\bigg(1+\frac{|f(x)-f(y)|}{|x-y|^{s}}\bigg)\leq g(x)+g(y)
\]
and by the lattice property of $L(p,\infty)({\mathbb{R}}^{N}\times{\mathbb{R}%
}^{N})$, we arrive at
\begin{equation}
\bigg\|\log\bigg(1+\frac{|f(x)-f(y)|}{|x-y|^{s}}\bigg)\frac{1}{|x-y|^{\frac
{N}{p}}}\bigg\|_{L(p,\infty)(\mathbb{R}^{N}\times{\mathbb{R}}^{N})}%
\leq2\bigg\|\frac{g(x)}{|x-y|^{\frac{N}{p}}}\bigg\|_{L(p,\infty)(\mathbb{R}%
^{N}\times{\mathbb{R}}^{N})}. \label{CL6}%
\end{equation}
The right-hand side can be estimated by using Proposition \ref{PropMixed}.
Indeed, we have
\begin{align*}
\bigg\|\frac{g(x)}{|x-y|^{\frac{N}{p}}}\bigg\|_{L(p,\infty)(\mathbb{R}%
^{N}\times{\mathbb{R}}^{N})}  &  \leq\bigg\|\bigg\|\frac{1}{|x-y|^{\frac{N}%
{p}}}\bigg\|_{L(p,\infty)({\mathbb{R}}^{N})}g(x)\bigg\|_{L^{p}({\mathbb{R}%
}^{N})}\\
&  \asymp\Vert g\Vert_{L^{p}({\mathbb{R}}^{N})}\text{ (by \eqref{ThmCCBSY3}).}%
\end{align*}
Inserting this estimate into \eqref{CL6} we infer that
\[
\bigg\|\log\bigg(1+\frac{|f(x)-f(y)|}{|x-y|^{s}}\bigg)\frac{1}{|x-y|^{\frac
{N}{p}}}\bigg\|_{L(p,\infty)(\mathbb{R}^{N}\times{\mathbb{R}}^{N})}%
\lesssim\Vert g\Vert_{L^{p}({\mathbb{R}}^{N})}.
\]
Taking the infimum over all $g$, we obtain
\[
\bigg\|\log\bigg(1+\frac{|f(x)-f(y)|}{|x-y|^{s}}\bigg)\frac{1}{|x-y|^{\frac
{N}{p}}}\bigg\|_{L(p,\infty)(\mathbb{R}^{N}\times{\mathbb{R}}^{N})}%
\lesssim\Vert f\Vert_{\mathcal{N}^{s,p}(\mathbb{R}^{N})}.
\]

\end{proof}

\begin{remark}
Let $s=1$. Then the previous result admits extensions to a rich class of
metric measure spaces. In particular, let $(X,d,m)$ be a $\text{RCD}^{\ast
}(K,N)$ space, where $m$ is an $n$-Ahlfors regular probability measure for
some $1<n\leq N$. The obtain the corresponding result in this context we only
need to mimic the proof of \eqref{CL*} but now applying the counterpart of
\eqref{Exp}, in the setting of metric spaces, recently obtained in
\cite[Theorem 3.11]{BrueSemola}. We have to leave the details to the reader.
\end{remark}

In the rest of this section we will focus on the converse to Theorem
\ref{BrezisExp} when $s=1$. In this regard, we establish the following

\begin{proposition}
\label{ThmBN} Let $f \in C^{1}({\mathbb{R}}^{N})$. Then
\[
\int_{{\mathbb{R}}^{N}} \log\, (1+ |\nabla f (x)|)^{p} \, \mathrm{d} x
\lesssim\liminf_{\lambda\to\infty} \lambda^{p} \mathcal{L}^{2N}
\bigg( \bigg\{(x,y) \in{\mathbb{R}}^{N} \times{\mathbb{R}}^{N} :\log\bigg(1 +
\frac{|f(x) - f(y)|}{|x-y|} \bigg) \frac{1}{|x-y|^{\frac{N}{p}}} >
\lambda\bigg\} \bigg).
\]

\end{proposition}

\begin{remark}
In view of \eqref{CL*}, Proposition \ref{ThmBN} sharpens the inequality
obtained in \cite[Proposition 2.6]{Brue}, namely,
\[
\int_{{\mathbb{R}}^{N}} \log\, (1+ |\nabla f (x)|)^{p} \, \mathrm{d} x
\lesssim\|f\|_{\mathcal{N}^{1,p}({\mathbb{R}}^{N})}^{p}.
\]

\end{remark}

\begin{proof}
[Proof of Proposition \ref{ThmBN}]Suppose first $N=1$. Since
\[
\lim_{t\rightarrow0+}\log\bigg(1+\frac{|f(x+t)-f(x)|}{t}\bigg)=\log
(1+|f^{\prime}(x)|),
\]
we can invoke Theorem \ref{ThmMaximal}(ii) to derive
\[
\int_{{\mathbb{R}}}\log\, (1+|f^{\prime}(x)|)^{p}\,\mathrm{d}x\leq
\liminf_{\lambda\rightarrow\infty}\lambda^{p}\mathcal{L}^{2}%
\bigg(\bigg\{(x,t)\in{\mathbb{R}}\times(0,\infty):\log\bigg(1+\frac
{|f(x+t)-f(x)|}{t}\bigg)\frac{1}{t^{\frac{1}{p}}}>\lambda\bigg\}\bigg).
\]
This completes the proof in the case $N=1$. The case $N>1$ can be done in a
similar fashion by a simple adaptation of the arguments given in the proof of
Theorem \ref{ThmMaximal}. We shall leave the somewhat tedious details to the reader.
\end{proof}

\section{Brezis--Van Schaftingen--Yung inequalities via Caffarelli--Silvestre
extensions\label{sec:Caff-Silv}}

In \cite{Brezis} the local properties of $\nabla f$ play an essential role,
which is to be contrasted with the non-local operator $(-\Delta)^{s}%
,\,s\in(0,1).$ Here we overcome the localization issues by means of employing
the celebrated Caffarelli--Silvestre extension theorem \cite{Caffarelli}. We
will not repeat here the formulation of this theorem, but we refer the reader
to Section \ref{sec:introduction}.


Next we establish Brezis--Van Schaftingen--Yung type inequalities in terms of
$(-\Delta)^{s}$ and the Riesz potential space $H^{2s,p}({\mathbb{R}}^{N})$
(cf. \eqref{IntroRiesz}.)
%

\begin{theorem}
[Brezis--Van Schaftingen--Yung inequalities for the fractional Laplacian]%
\label{ThmCSBSY} Let $s\in(0,1)$ and $1 < p < \infty$. Let $P_{(-\Delta)^{s}%
}[f]$ be the Caffarelli--Silvestre extension of $f \in\mathcal{S}({\mathbb{R}%
}^{N})$ (cf. \eqref{CSExplicit}). Then
\begin{equation}
\label{ThmCSBSYppp}\Vert f\Vert_{H^{2s,p}(\mathbb{R}^{N})}\leq\mu
_{s}\bigg\|\frac{P_{(-\Delta)^{s}} [f](x,t) - f(x)}{t^{2s+\frac{1}{p}}%
}\bigg\|_{L(p,\infty)(\mathbb{R}_{+}^{N+1})}%
\end{equation}
where $\mu_{s}$ is the constant appearing in the formula \eqref{combinada}.
\end{theorem}

\begin{remark}
Specializing the previous result with $s=1/2$ and applying the Riesz's
theorem, we infer that
\begin{equation}
\Vert\nabla f\Vert_{L^{p}({\mathbb{R}}^{N})}\asymp\Vert(-\Delta)^{1/2}%
f\Vert_{L^{p}(\mathbb{R}^{N})}\lesssim\bigg\|\frac{P [f](x,t) - f(x)}%
{t^{1+\frac{1}{p}}}\bigg\|_{L(p,\infty)(\mathbb{R}_{+}^{N+1})} \label{rec1}%
\end{equation}
where $P[f]$ is the Poisson extension of $f$ (cf. \eqref{ooo}). Informally,
this estimate recovers (cf. \eqref{I1})
\[
\Vert\nabla f\Vert_{L^{p}({\mathbb{R}}^{N})}\lesssim\left\Vert \frac
{f(x)-f(y)}{\left\vert x-y\right\vert ^{1+ \frac{N}{p}}}\right\Vert
_{L(p,\infty)({\mathbb{R}}^{N}\times{\mathbb{R}}^{N})}%
\]
modulo the \textquotedblleft change of variables" $f(x)\in{\mathbb{R}}%
^{N}\leftrightarrow P[f](x,t)\in{\mathbb{R}}^{N}\times(0,\infty)$.
\end{remark}

\begin{proof}
[Proof of Theorem \ref{ThmCSBSY}]According to \eqref{combinada} and
\eqref{CSExplicit}, we have
\begin{equation}
(-\Delta)^{s}f(x)=-\mu_{s}\lim_{t\rightarrow0+}\frac{u_{(-\Delta)^{s}%
}(x,t)-u_{(-\Delta)^{s}}(x,0)}{t^{2s}}=-\mu_{s}\lim_{t\rightarrow0+}%
\frac{P_{(-\Delta)^{s}} [f](x,t) - f(x)}{t^{2s}}. \label{ProofCS1}%
\end{equation}
Let us introduce the sequence of operators given by
\[
T_{t}f(x)=\frac{P_{(-\Delta)^{s}} [f](x,t) -f(x)}{t^{2s}},\quad x\in
{\mathbb{R}}^{N},\quad t>0.
\]
Then we may rewrite \eqref{ProofCS1} as
\[
\mu_{s}^{-1}(-\Delta)^{s}f(x)=-\lim_{t\rightarrow0+}T_{t}f(x).
\]
At this point we are in position to apply Theorem \ref{ThmMaximal}. Indeed, it
follows from Theorem \ref{ThmMaximal}(ii) (with $\gamma=1$) that
\begin{align*}
\mu_{s}^{-p}\Vert(-\Delta)^{s}f\Vert_{L^{p}({\mathbb{R}}^{N})}^{p}  &
\leq\liminf_{\lambda\rightarrow\infty}\lambda^{p}\mathcal{L}^{N+1}%
\bigg(\bigg\{(x,t)\in{\mathbb{R}}_{+}^{N+1}:\frac{|T_{t}f(x)|}{t^{1/p}%
}>\lambda\bigg\}\bigg)\\
&  \leq\sup_{\lambda>0}\lambda^{p}\mathcal{L}^{N+1}\bigg(\bigg\{(x,t)\in
{\mathbb{R}}_{+}^{N+1}:\frac{|T_{t}f(x)|}{t^{1/p}}>\lambda\bigg\}\bigg)\\
&  =\bigg\|\frac{P_{(-\Delta)^{s}} [f](x,t) -f(x)}{t^{2s+1/p}}%
\bigg\|_{L(p,\infty)({\mathbb{R}}_{+}^{N+1})}^{p}.
\end{align*}

\end{proof}

\begin{remark}
In fact, the proof of Theorem \ref{ThmCSBSY} shows a stronger inequality than
\eqref{ThmCSBSYppp}. Namely
\[
\mu_{s}^{-p}\Vert f\Vert_{H^{2s,p}({\mathbb{R}}^{N})}^{p} \leq\liminf
_{\lambda\rightarrow\infty}\lambda^{p}\mathcal{L}^{N+1}\bigg(\bigg\{(x,t)\in
{\mathbb{R}}_{+}^{N+1}:\frac{|P_{(-\Delta)^{s}} [f](x,t) -f(x)|}{t^{2s+1/p}%
}>\lambda\bigg\}\bigg).
\]

\end{remark}

The previous theorem can be examined within the more general framework of
non-local hypoelliptic operators in the sense of H\"{o}rmander
\cite{Hormander}. To be more precise, consider the class of degenerate
equations (see \eqref{Hor})
\begin{equation}
\mathscr{K}u=\mathscr{A}u-u_{t}=\text{tr}(Q\nabla^{2}u)+\langle Bx,\nabla
u\rangle-u_{t}=0\quad\text{in}\quad{\mathbb{R}}_{+}^{N+1} \label{Hormander2}%
\end{equation}
(cf. Example \ref{ExampleHormander}.)

The generalization of the Caffarelli--Silvestre theorem \eqref{11111},
\eqref{CSExplicit} and \eqref{combinada} for the fractional Laplacian
$(-\Delta)^{s}$ to fractional powers of the H\"{o}rmander equation
$(-\mathscr{K})^{s}$ and its diffusive part $(-\mathscr{A})^{s}$ has been
recently obtained by Garofalo and Tralli \cite{Garofalo20}. Here, the
operators $(-\mathscr{A})^{s}$ and $(-\mathscr{K})^{s}$ are defined on
$\mathcal{S}({\mathbb{R}}^{N})$ by means of the corresponding H\"{o}rmander
semigroups; and we refer to \cite{Garofalo20} for precise definitions. Then
the Caffarelli--Silvestre extension theorem for $(-\mathscr{A})^{s}%
,\,s\in(0,1)$, reads as follows \cite[Theorem 5.5]{Garofalo20}: If
$f\in\mathcal{S}({\mathbb{R}}^{N})$ then there exists a function
$u_{(-\mathscr{A})^{s}}\equiv u\in C^{\infty}({\mathbb{R}}_{+}^{N+1})$ such
that
\begin{equation}
\left\{
\begin{array}
[c]{ll}%
u(x,0)=f(x) & \text{on}\quad\mathbb{R}^{N},\\
t^{1-2s}(\mathscr{A}u+\frac{1-2s}{t}u_{t}+u_{tt})=0 & \text{in}\quad
\mathbb{R}_{+}^{N+1}.
\end{array}
\right.  \label{CSGarofaloPDE}%
\end{equation}
Moreover, we also have in $L^{\infty}({\mathbb{R}}^{N})$
\begin{equation}
(-\mathscr{A})^{s}f(x)=-\mu_{s}\lim_{t\rightarrow0+}\frac{u(x,t)-u(x,0)}%
{t^{2s}}. \label{CSGarofalo}%
\end{equation}

As far as the extension theorem for $(-\mathscr{K})^{s}$, it was shown in
\cite[Theorem 4.1]{Garofalo20} that if $f\in\mathcal{S}({\mathbb{R}}^{N+1})$
then there exists a function $u_{(-\mathscr{K})^{s}}\equiv u\in C^{\infty
}({\mathbb{R}}_{+}^{N+2})$ such that
\begin{equation}
\left\{
\begin{array}
[c]{ll}%
u(x,0)=f(x) & \text{on}\quad\mathbb{R}^{N+1},\\
t^{1-2s}(\mathscr{A}u+\frac{1-2s}{t}u_{t}+u_{tt}-u_{t})=0 & \text{in}%
\quad\mathbb{R}_{+}^{N+2}.
\end{array}
\right.  \label{CSGarofaloPDE2}%
\end{equation}
Moreover, we also have that, in $L^{\infty}({\mathbb{R}}^{N}),$
\begin{equation}
(-\mathscr{K})^{s}f(x)=-\mu_{s}\lim_{t\rightarrow0+}\frac{u(x,t)-u(x,0)}%
{t^{2s}}. \label{CSGarofalo2}%
\end{equation}

The solutions $u_{(-\mathscr{A})^{s}}$ and $u_{(-\mathscr{K})^{s}}$ of
\eqref{CSGarofaloPDE} and \eqref{CSGarofaloPDE2}, respectively, can be written
explicitly in terms of $f$ via the corresponding Poisson kernels
$P_{(-\mathscr{A})^{s}}$ and $P_{(-\mathscr{K})^{s}}$; cf. \cite{Garofalo20}
for the precise definitions. Consequently, one can rewrite \eqref{CSGarofalo}
and \eqref{CSGarofalo2} as
\begin{equation}
(-\mathscr{A})^{s}f(x)=-\mu_{s}\lim_{t\rightarrow0+}\frac
{P_{(-\mathscr{A})^{s}}[f](x,t)-f(x)}{t^{2s}} \label{CSGarofalo4444}%
\end{equation}
and
\begin{equation}
(-\mathscr{K})^{s}f(x)=-\mu_{s}\lim_{t\rightarrow0+}\frac
{P_{(-\mathscr{K})^{s}}[f](x,t)-f(x)}{t^{2s}}. \label{CSGarofalo24444}%
\end{equation}

Having the extension theorems \eqref{CSGarofaloPDE} and \eqref{CSGarofalo4444}
(respectively, \eqref{CSGarofaloPDE2} and \eqref{CSGarofalo24444}) for
$(-\mathscr{A})^{s}$ (respectively, $(-\mathscr{K})^{s}$) in hand, we are now
able to apply the method of proof of Theorem \ref{ThmCSBSY} to establish the
following Brezis--Van Schaftingen--Yung-type inequalities.

\begin{theorem}
[Brezis--Van Schaftingen--Yung inequalities for the diffusive part of
H\"{o}rmander equations]\label{ThmCSBSYGarofalo} Let $f\in\mathcal{S}%
({\mathbb{R}}^{N})$, and let $s\in(0,1), 1 < p < \infty$. Then
\[
\Vert f\Vert_{H_{\mathscr{A}}^{2s,p}(\mathbb{R}^{N})}\leq\mu_{s}%
\bigg\|\frac{P_{(-\mathscr{A})^{s}}[f](x,t)-f(x)} {t^{2s+\frac{1}{p}}%
}\bigg\|_{L(p,\infty)(\mathbb{R}_{+}^{N+1})}%
\]
where
\[
\Vert f\Vert_{H_{\mathscr{A}}^{2s,p}(\mathbb{R}^{N})}:=\Vert(-\mathscr{A})^{s}%
f\Vert_{L^{p}({\mathbb{R}}^{N})}.
\]

\end{theorem}

\begin{theorem}
[Brezis--Van Schaftingen--Yung inequalities for H\"{o}rmander equations]%
\label{ThmCSBSYGarofalo2} Let $f\in\mathcal{S}({\mathbb{R}}^{N})$, and let
$s\in(0,1), 1 < p < \infty$. Then
\[
\Vert f\Vert_{H_{\mathscr{K}}^{2s,p}(\mathbb{R}^{N+1})}\leq\mu_{s}%
\bigg\|\frac{P_{(-\mathscr{K})^{s}}[f](x,t)-f(x)} {t^{2s+\frac{1}{p}}%
}\bigg\|_{L(p,\infty)(\mathbb{R}_{+}^{N+2})}
\]
where
\[
\Vert f\Vert_{H_{\mathscr{K}}^{2s,p}(\mathbb{R}^{N+1})}:=\Vert
(-\mathscr{K})^{s}f\Vert_{L^{p}({\mathbb{R}}^{N+1})}.
\]

\end{theorem}

\begin{remark}
Setting $Q=I_{N}$ and $B=0_{N}$ in \eqref{Hormander2} one has
$\mathscr{A}=\Delta$ and thus Theorem \ref{ThmCSBSYGarofalo} covers Theorem
\ref{ThmCSBSY}. On the other hand, $\mathscr{K}=\Delta-u_{t}$ and thus Theorem
\ref{ThmCSBSYGarofalo2} may be regarded as the caloric counterpart of Theorem
\ref{ThmCSBSY}. However, the applicability of Theorems \ref{ThmCSBSYGarofalo}
and \ref{ThmCSBSYGarofalo2} goes well beyond the Laplace and heat equations,
providing a wide variety of Brezis--Van Schaftingen--Yung inequalities related
to the Ornstein--Uhlenbeck equation, Kolmogorov equation, Kramers equation,
... cf. \cite{Garofalo20, Garofalo} for further details.
\end{remark}

Our approach is not limited only to the class of PDE's given by
\eqref{Hormander2}, but also can be applied to basic Schr\"{o}dinger operators
in ${\mathbb{R}}^{N}$. Indeed, consider the harmonic oscillator
\[
\mathcal{H}=-\Delta+|x|^{2}%
\]
and its fractional powers $\mathcal{H}^{s},\,s\in(0,1)$. For smooth functions
$f$, the operator $\mathcal{H}^{s}f$ can be defined by a pointwise formula
involving the corresponding heat semigroup (cf. \cite[Remark 4.3]{Stinga}.)

The extension theorem associated with $\mathcal{H}^{s}$ was shown by Stinga
and Torrea \cite[Theorem 4.2]{Stinga}. It claims that if $f$ is a
$C^{2}({\mathbb{R}}^{N})$ function with a certain decay at infinity, then
there exists a function $u_{\mathcal{H}^{s}}\equiv u$ on ${\mathbb{R}}%
_{+}^{N+1},$ such that
\[
\left\{
\begin{array}
[c]{ll}%
u(x,0)=f(x) & \text{on}\quad\mathbb{R}^{N},\\
\mathcal{H}_{x}u+\frac{1-2s}{t}u_{t}+u_{tt}=0 & \text{in}\quad\mathbb{R}%
_{+}^{N+1}.
\end{array}
\right.
\]
Moreover, if $P_{\mathcal{H}^{s}}$ denotes the corresponding Poisson kernel
then, up to a multiplicative constant depending only on $s$,
\begin{equation}
\mathcal{H}^{s}f(x)=-\lim_{t\rightarrow0+}\frac{P_{\mathcal{H}^{s}%
}[f](x,t)-f(x)}{t^{2s}} \label{CSStinga}%
\end{equation}
where the convergence is understood in the classical sense.

\begin{theorem}
[Brezis--Van Schaftingen--Yung inequalities for the harmonic oscillator]%
\label{ThmCSBSYStinga} Let $f\in C^{2}({\mathbb{R}}^{N})$, and let $s\in(0,1),
1 < p < \infty$. Then there exists a constant $C,$ depending only on $s,$ such
that
\[
\Vert f\Vert_{H_{\mathcal{H}}^{2s,p}(\mathbb{R}^{N})}\leq C \, \bigg\|\frac
{P_{\mathcal{H}^{s}}[f](x,t)-f(x)}{t^{2s+\frac{1}{p}}}\bigg\|_{L(p,\infty
)(\mathbb{R}_{+}^{N+1})},
\]
where
\[
\Vert f\Vert_{H_{\mathcal{H}}^{2s,p}(\mathbb{R}^{N})}:=\Vert(-\mathcal{H}%
)^{s}f\Vert_{L^{p}({\mathbb{R}}^{N})}.
\]

\end{theorem}

The proof of Theorem \ref{ThmCSBSYStinga} follows line by line the proof of
Theorem \ref{ThmCSBSY}, but now relying on the formula \eqref{CSStinga}. The
rigorous details can be safely left to the reader.

The set of techniques presented here allow us to address similar problems on
other manifolds. For instance, one may consider fractional powers related to
the Laplace--Beltrami operator $\Delta_{\mathbb{S}^{N-1}}$ on the unit sphere
$\mathbb{S}^{N-1}$. As in the Euclidean setting, the fractional powers
$(-\Delta_{\mathbb{S}^{N-1}})^{s}$ can be introduced using spherical
harmonics. In this setting, the Caffarelli--Silvestre extension theorem was
recently obtained in \cite[Theorem 7.2]{DeNapoli}. Accordingly, applying the
methodology described in the proof of Theorem \ref{ThmCSBSY}, we can show the following

\begin{theorem}
[Brezis--Van Schaftingen--Yung inequality on the sphere]Let $f\in C^{\infty
}(\mathbb{S}^{N-1})$, and let $s\in(0,1), 1 < p < \infty$. Then
\[
\Vert(-\Delta_{\mathbb{S}^{N-1}})^{s}f\Vert_{L^{p}(\mathbb{S}^{N-1}%
,\mathrm{d}\sigma^{N-1})}\lesssim\bigg\|\frac{P_{(-\Delta_{\mathbb{S}^{N-1}%
})^{s}}[f](x,t)-f(x)}{t^{2s+\frac{1}{p}}}\bigg\|_{L(p,\infty)(\mathbb{S}%
^{N-1}\times(0,\infty),\mathrm{d}\sigma^{N-1}\times\mathcal{L})}%
\]
where $P_{(-\Delta_{\mathbb{S}^{N-1}})^{s}}[f]$ solves the boundary value
problem
\[
\left\{
\begin{array}
[c]{ll}%
u(x,0)=f(x) & \text{on}\quad\mathbb{S}^{N-1},\\
\Delta_{\mathbb{S}^{N-1}}u+\frac{1-2s}{t}u_{t}+u_{tt}=0 & \text{in}%
\quad\mathbb{S}^{N-1}\times(0,\infty).
\end{array}
\right.
\]

\end{theorem}

We can also deal with fractional operators on Carnot groups. Recall that a
simply connected Lie group $\mathbb{G}$ is a \emph{Carnot group of step two}
if its Lie algebra $\mathfrak{g}$ admits a stratification $\mathfrak{g}%
=V_{1}\oplus V_{2}$ which is $2$-nilpotent, that is, $[V_{1},V_{1}]=V_{2}$ and
$[V_{1},V_{2}]=\{0\}$. Suppose that $\mathfrak{g}$ is endowed with a scalar
product. Let $m=\text{dim}(V_{1})$ and $k=\text{dim}(V_{2})$. If
$\{e_{1},\ldots,e_{m}\}$ is a orthonormal basis on $V_{1},$ then the
associated horizontal Laplacian $\Delta_{H}$ is defined by
\[
\Delta_{H}f=\sum_{j=1}^{m}X_{j}^{2}f,
\]
where $X_{j}$ is the vector field on $\mathbb{G}$ given by Lie's rule
\[
X_{j}f(g)=\frac{d}{ds}f(g\circ\exp se_{j})\big\lvert_{s=0}.
\]
In the special case $\mathbb{G}=\mathbb{H}^{N}$, the Heisenberg group, the
underlying manifold is ${\mathbb{R}}^{2N+1}$, and the Laplacian $\Delta_{H}$
can be written in real coordinates $(x,y,t)$ as
\[
\Delta_{H}=\Delta_{x,y}+\frac{|x|^{2}+|y|^{2}}{4}\partial_{t}^{2}+\partial
_{t}\Big(\sum_{j=1}^{N}(x_{j}\partial_{y_{j}}-y_{j}\partial_{x_{j}})\Big);
\]
see \cite{Garofalo16}.

Following \cite{Folland} and \cite{Garofalo18}, one can introduce in a natural
way the fractional powers $(\partial_{t} - \Delta_{H})^{s}$ and $(-\Delta
_{H})^{s}, \, s \in(0,1),$ by pointwise representations in terms of
(evolutive) heat semigroups. Moreover, these fractional operators can also be
represented by the Dirichlet-to-Neumann relations provided by
Caffarelli--Silvestre extensions. Specifically, one has that
\begin{equation}
\label{CSGarofaloPDE2+1}- C_{s} \lim_{y \to0+} \frac{u(g,t,y) - u(g,t,0)}{y^{2
s}} = (\partial_{t} - \Delta_{H})^{s} f (g,t), \quad u \in C^{\infty
}(\mathbb{G} \times{\mathbb{R}}),
\end{equation}
where
\begin{equation}
\label{CSGarofaloPDE2+}\left\{
\begin{array}
[c]{ll}%
u(g,t,0) = f(g,t) & \text{on} \quad\mathbb{G} \times{\mathbb{R}},\\
\Delta_{H} u + \frac{1- 2 s}{y} u_{y} + u_{y y} - u_{t} = 0 & \text{in}
\quad\mathbb{G} \times{\mathbb{R}} \times(0,\infty),
\end{array}
\right.
\end{equation}
and
\begin{equation}
\label{CSGarofaloPDE2+2}- C_{s} \lim_{y \to0+} \frac{u(g,y) - u(g,0)}{y^{2 s}}
= (- \Delta_{H})^{s} f (g), \quad u \in C^{\infty}(\mathbb{G}),
\end{equation}
where
\begin{equation}
\label{CSGarofaloPDE2++}\left\{
\begin{array}
[c]{ll}%
u(g,0) = f(g) & \text{on} \quad\mathbb{G},\\
\Delta_{H} u + \frac{1- 2 s}{y} u_{y} + u_{y y} = 0 & \text{in} \quad
\mathbb{G} \times(0,\infty).
\end{array}
\right.
\end{equation}
Here, the limits \eqref{CSGarofaloPDE2+1} and \eqref{CSGarofaloPDE2+2} are not
only pointwise, but also hold in $L^{p}$ for any $1 \leq p \leq\infty$.

\begin{theorem}
[Brezis--Van Schaftingen--Yung inequalities on Carnot groups]\label{ThmCarnot}
Let $s\in(0,1)$.

\begin{enumerate}
\item Let $f \in C^{\infty}(\mathbb{G} \times{\mathbb{R}})$. Then
\[
\| (\partial_{t}-\Delta_{H})^{s} f\|_{L^{p}(\mathbb{G} \times{\mathbb{R}},
\mathcal{L}^{m+k+1})} \lesssim\bigg\|\frac{u(g,t,y)-f(g,t)}{y^{2 s + \frac
{1}{p}}} \bigg\|_{L(p,\infty)(\mathbb{G} \times{\mathbb{R}} \times(0,\infty),
\mathcal{L}^{m+k+2})}%
\]
where $u$ solves \eqref{CSGarofaloPDE2+}.

\item Let $f\in C^{\infty}(\mathbb{G})$. Then
\[
\Vert(-\Delta_{H})^{s}f\Vert_{L^{p}(\mathbb{G},\mathcal{L}^{m+k})}%
\lesssim\bigg\|\frac{u(g,y)-f(g)}{y^{2s+\frac{1}{p}}}\bigg\|_{L(p,\infty
)(\mathbb{G}\times(0,\infty),\mathcal{L}^{m+k+1})}%
\]
where $u$ solves \eqref{CSGarofaloPDE2++}.
\end{enumerate}
\end{theorem}

\section{Application: Identifying constant functions via Harnack's
inequalities\label{sec:further}}

The aim of this final section is to establish new criteria which enable us to
identify constant functions. Before we present our contributions to this
topic, we gather some known results which are relevant for our purposes.

The Bourgain--Brezis--Mironescu theorem \cite{Bourgain00} tells us that if
$u:{\mathbb{R}}^{N}\rightarrow{\mathbb{R}}$ is a measurable function and
$1\leq p<\infty$ satisfying
\begin{equation}
\bigg\|\frac{u(x)-u(y)}{|x-y|^{1+\frac{N}{p}}}\bigg\|_{L^{p}({\mathbb{R}}%
^{N}\times{\mathbb{R}}^{N})}^{p}=\iint\limits_{{\mathbb{R}}^{N}\times
{\mathbb{R}}^{N}}\frac{|u(x)-u(y)|^{p}}{|x-y|^{p+N}}\,\mathrm{d}%
x\,\mathrm{d}y<\infty\label{BBMConstant}%
\end{equation}
then $u$ is a constant function; see also \cite{Brezis02}, \cite{davila} and
\cite{Milman}.

Very recently, Brezis, Van Schaftingen and Yung established \cite[Proposition
6.3]{Brezis} a weak-$L^{p}$ variant of \eqref{BBMConstant}. Specifically, they
showed that if $u : {\mathbb{R}}^{N} \to{\mathbb{R}}$ is measurable, $1 < p <
\infty$ and
\begin{equation}
\label{BBMConstant2}\lim_{\lambda\to\infty} \lambda^{p} \mathcal{L}^{2 N}
\bigg( \bigg\{(x,y) \in{\mathbb{R}}^{N} \times{\mathbb{R}}^{N} : \frac{|u(x)
-u(y)|}{|x-y|^{1+\frac{N}{p}}} \geq\lambda\bigg\} \bigg) = 0
\end{equation}
then $u$ must be a constant.

Combining our previous results together with PDE's techniques, namely,
Harnack's inequalities for fractional operators, we will be able to establish
(scale-invariant) fractional counterparts of \eqref{BBMConstant2} in terms of
Caffarelli--Silvestre extensions.

\begin{corollary}
[Constancy via fractional Laplacian/harmonic oscillator]\label{CorollaryBBM}
Let $f$ be a $\mathcal{S}({\mathbb{R}}^{N})$ nonnegative function, $R>0$ and
$s\in(0,1), 1 < p < \infty$. For $\mathfrak{R}\in\{-\Delta,\mathcal{H}\}$, we
let $P_{\mathfrak{R}^{s}}[f] \equiv P[f]$ be the
Poisson--Caffarelli--Silvestre extension of $f$. Suppose that
\begin{equation}
\lim_{\lambda\rightarrow\infty}\lambda^{p}\mathcal{L}^{N+1}%
\bigg(\bigg\{(x,t)\in B(0,R)\times(0,\infty):\frac{|P[f](x,t)-f(x)|}%
{t^{2s+\frac{1}{p}}}\geq\lambda\bigg\}\bigg)=0. \label{BBMConstant3}%
\end{equation}
Then there exists $C>0$, which depends only on $s,p$ and $N$, such that
\[
\sup_{B(0,R/2)}f\leq C\inf_{B(0,R/2)}f.
\]

\end{corollary}

\begin{remark}
The condition \eqref{BBMConstant3} with $s=1/2$ corresponds to
\eqref{BBMConstant2} (modulo the ``change of variables" $f(x) \in{\mathbb{R}%
}^{N} \leftrightarrow P[f](x,t) \in{\mathbb{R}}^{N} \times(0,\infty)$.)
\end{remark}

\begin{proof}
[Proof of Corollary \ref{CorollaryBBM}]Consider the sequence of operators
defined by
\[
T_{t} f (x) := \frac{P[f](x,t)-f(x)}{t^{2 s}}, \quad x \in B(0,R), \quad t >
0.
\]
According to Theorem \ref{ThmMaximal}(ii), \eqref{ProofCS1} and
\eqref{CSStinga}, we have
\begin{align}
\| \mathfrak{R}^{s} f\|_{L^{p}(B(0,R))}^{p}  &  \lesssim\nonumber\\
&  \hspace{-2cm} \lim_{\lambda\to\infty} \lambda^{p}\mathcal{L}^{N+1}
\bigg( \bigg\{(x,t) \in B(0,R) \times(0,\infty) : \frac{|P[f](x,t)-f(x)|}{t^{2
s +\frac{1}{p}}} \geq\lambda\bigg\} \bigg). \label{BBMConstant4}%
\end{align}
Since $(-\Delta)^{s} f \in C^{\infty}({\mathbb{R}}^{N})$ (respectively,
$\mathcal{H}^{s} f \in\mathcal{S}({\mathbb{R}}^{N})$), it follows from
\eqref{BBMConstant3} and \eqref{BBMConstant4} that $\mathfrak{R}^{s} f = 0$ in
$B(0,R)$. According to the Harnack inequalities for $\mathfrak{R}=-\Delta$
(cf. \cite{Landkof} or \cite{Caffarelli}) and for $\mathfrak{R}=\mathcal{H}$
(cf. \cite{Stinga}), we conclude the desired result.
\end{proof}

The above proof can be easily adapted to deal with fractional heat operators
$(\partial_{t} - \Delta)^{s}$ via the formula \eqref{CSGarofalo24444} (with
$Q= I_{N}$ and $B=0_{N}$) and the corresponding Harnack inequality given in
\cite[Theorem 5.2]{Banerjee}.

\begin{corollary}
[Constancy via fractional heat operator]\label{CorollaryBBMHeat} Let $f$ be a
$\mathcal{S}({\mathbb{R}}^{N+1})$ nonnegative function, $R>0$ and $s\in(0,1),
1 < p < \infty$. Suppose that
\[
\lim_{\lambda\rightarrow\infty}\lambda^{p}\mathcal{L}^{N+2}%
\bigg(\bigg\{(x,t)\in B(0,R)\times(0,\infty):\frac{|P_{(\partial_{t} -
\Delta)^{s}}[f](x,t)-f(x)|}{t^{2s+\frac{1}{p}}}\geq\lambda\bigg\}\bigg)=0.
\]
Then there exists $C>0$, which is independent of $R$, such that
\[
\sup_{B(0,R/2)}f\leq C\inf_{B(0,R/2)}f.
\]

\end{corollary}

Our approach also allows us to obtain conditions implying constancy for
functions on the sphere via the corresponding Harnack's inequality given in
\cite[Theorem 7.4]{DeNapoli}.

\begin{corollary}
[Constancy of functions on the sphere]\label{CorollaryBBMSphere} Let $f$ be a
$C^{\infty}(\mathbb{S}^{N-1})$ nonnegative function and $s\in(0,1),1<p<\infty
$. Let $\Omega^{\prime}\subset\subset\Omega\subset\mathbb{S}^{N-1}$ be open
sets. Suppose that
\[
\lim_{\lambda\rightarrow\infty}\lambda^{p}(\mathrm{d}\sigma^{N-1}%
\times\mathcal{L})\bigg(\bigg\{(x,t)\in\Omega\times(0,\infty):\frac
{|P_{(-\Delta_{\mathbb{S}^{N-1}})^{s}}[f](x,t)-u(x,0)|}{t^{2s+\frac{1}{p}}%
}\geq\lambda\bigg\}\bigg)=0.
\]
Then there exists $C>0$, which depends on $\Omega^{\prime},\Omega,N$ and $s$,
such that \ \
\[
\sup_{\Omega^{\prime}}f\leq C\inf_{\Omega^{\prime}}f.
\]

\end{corollary}

\begin{remark}
Combining in the above fashion Theorem \ref{ThmCarnot} with the Harnack's
inequality for fractional sub-Laplacians on $\mathbb{G}$ (cf. \cite{Ferrari}),
we may establish a scale-invariant criterion to identify constant functions on
Carnot groups. Further details are left to the reader.
\end{remark}

\section{Appendix: Further results and applications\label{appendix}}

\subsection{Connections and comparisons with the work of Bourgain--Nguyen}\label{SectionBN}

Following comments and suggestions by the referees we prove a variant of
Theorem \ref{ThmMaximal} that allows to show how the results of
Bourgain--Nguyen (cf. (\ref{nb1}), (\ref{nb2}), (\ref{nb3})) fit in our
abstract framework. As a result we are able to clarify the connections between
the Brezis--Van Schaftingen--Yung and the Bourgain--Nguyen approaches. The key
result in this direction is the following

\begin{theorem}
\label{ThmMaximal2}Let $(X,m)$ be a $\sigma$-finite measure space, and let
$\{T_{t}:t>0\}$ be a one-parameter family of (not necessarily linear)
operators on $L^{p}(X,m),\,1\leq p<\infty$. Let
$\gamma\neq0,$ and define
\[
L=\left\{
\begin{array}
[c]{cl}%
\infty & \text{if}\quad\gamma>0,\\
0 & \text{if}\quad\gamma<0.
\end{array}
\right.
\]

\begin{enumerate}
[\upshape(i)] 
\item Assume 
\begin{equation}\label{AssumptionMaximal}
T^{\ast}f:=\sup_{t>0}|T_{t}f|\in L^{p}(X,m).
\end{equation}
Then
\begin{equation}
\sup_{\lambda>0}\lambda^{p}\underset{t^{-\gamma/p}|T_{t}f(x)|>\lambda
}{\int_{X}\int_{0}^{\infty}}t^{\gamma-1}\,\mathrm{d}t\,\mathrm{d}m(x)\leq
\frac{1}{|\gamma|}\Vert T^{\ast}f\Vert_{L^{p}(X,m)}^{p}. \label{desired}%
\end{equation}
If, in addition, 
\begin{equation}\label{FinitenessLimit}
g(x)=\lim_{t\rightarrow0+}T_{t}f(x)\text{ exists and it is finite}%
\quad\text{$m$-a.e.}\quad x\in X,
\end{equation}
then
\begin{equation}
\frac{1}{|\gamma|}\,\|g\|_{L^{p}(X,m)}^{p} = \lim_{\lambda
\rightarrow L}\lambda^{p}\underset{t^{-\gamma/p}|T_{t}f(x)|>\lambda}{\int%
_{X}\int_{0}^{\infty}}t^{\gamma-1}\,\mathrm{d}t\,\mathrm{d}m(x).
\label{Statement1.2}%
\end{equation}

\item Assume that \eqref{FinitenessLimit} holds. Then
\begin{equation}
\frac{1}{|\gamma|}\, \|g \|_{L^{p}(X,m)}^{p}\leq\liminf_{\lambda
\rightarrow L}\lambda^{p}\underset{t^{-\gamma/p}|T_{t}f(x)|>\lambda}{\int%
_{X}\int_{0}^{\infty}}t^{\gamma-1}\,\mathrm{d}t\,\mathrm{d}m(x).
\label{Statement1}%
\end{equation}

\end{enumerate}
\end{theorem}

\begin{proof}
For each $\lambda>0$, we set
\[
E(f,\lambda)=\{(x,t)\in X\times(0,\infty):t^{-\gamma/p}|T_{t}f(x)|>\lambda\}.
\]
Clearly
\begin{equation}\label{MaximalInclusion}
E(f,\lambda)\subseteq\widetilde{E}(f,\lambda)=\{(x,t)\in X\times
(0,\infty):t^{-\gamma/p}T^{\ast}f(x)>\lambda\},
\end{equation}
consequently
\begin{equation}\label{TrivialEstim}
\underset{E(f,\lambda)}{\iint}t^{\gamma-1}\,\mathrm{d}t\,\mathrm{d}%
m(x)\leq\underset{\widetilde{E}(f,\lambda)}{\iint}t^{\gamma-1}\,\mathrm{d}%
t\,\mathrm{d}m(x).
\end{equation}
We claim that
\begin{equation}
\underset{\widetilde{E}(f,\lambda)}{\iint}t^{\gamma-1}\,\mathrm{d}%
t\,\mathrm{d}m(x)=\frac{1}{|\gamma|\lambda^{p}}\int_{X}(T^{\ast}%
f(x))^{p}\,\mathrm{d}m(x). \label{ProofUnifying1}%
\end{equation}
To prove \eqref{ProofUnifying1} we use Fubini's theorem, which naturally leads
us to consider two cases. Namely, if $\gamma>0,$ then we have
\[
\underset{\widetilde{E}(f,\lambda)}{\iint}t^{\gamma-1}\,\mathrm{d}%
t\,\mathrm{d}m(x)=\int_{X}\int_{0}^{(\frac{T^{\ast}f(x)}{\lambda})^{p/\gamma}%
}t^{\gamma-1}\,\mathrm{d}t\,\mathrm{d}m(x)=\frac{1}{\gamma\lambda^{p}}\int%
_{X}(T^{\ast}f(x))^{p}\,\mathrm{d}m(x),
\]
likewise, if $\gamma<0,$ then%
\[
\underset{\widetilde{E}(f,\lambda)}{\iint}t^{\gamma-1}\,\mathrm{d}%
t\,\mathrm{d}m(x)=\int_{X}\int_{(\frac{\lambda}{T^{\ast}f(x)})^{-p/\gamma}%
}^{\infty}t^{\gamma-1}\,\mathrm{d}t\,\mathrm{d}m(x)=\frac{1}{-\gamma
\lambda^{p}}\int_{X}(T^{\ast}f(x))^{p}\,\mathrm{d}m(x).
\]
Thus (\ref{ProofUnifying1}) holds and (\ref{desired}) readily follows from \eqref{TrivialEstim}.

A simple change of variables yields that, for $\lambda>0$,
\begin{equation}
\lambda^{p}\underset{t^{-\gamma/p}|T_{t}f(x)|>\lambda}{\int_{X}\int%
_{0}^{\infty}}t^{\gamma-1}\,\mathrm{d}t\,\mathrm{d}m(x)=\underset{t^{-\gamma
/p}|T_{t\lambda^{-p/\gamma}}f(x)|>1}{\int_{X}\int_{0}^{\infty}}t^{\gamma
-1}\,\mathrm{d}t\,\mathrm{d}m(x). \label{NewEquation}%
\end{equation}
According to \eqref{FinitenessLimit}, it is easy to see that
\begin{equation}\label{AnalogLimInf}
\chi_{\{(x,t):t^{-\gamma/p}|g(x)| > 1\}}  \leq \liminf_{\lambda\rightarrow L}\chi_{\{(x,t):t^{-\gamma/p}|T_{t\lambda
^{-p/\gamma}}f(x)|>1\}}.
\end{equation}
Therefore, combining \eqref{NewEquation} and Fatou's lemma yields,
\[
\liminf_{\lambda\rightarrow L}\lambda^{p}\underset{t^{-\gamma/p}%
|T_{t}f(x)|>\lambda}{\int_{X}\int_{0}^{\infty}}t^{\gamma-1}\,\mathrm{d}%
t\,\mathrm{d}m(x)\geq\int_{X}\int_{0}^{\infty}t^{\gamma-1}\chi
_{\{(x,t):t^{-\gamma/p}|g(x)| > 1\}}\,\mathrm{d}t\,\mathrm{d}m(x).
\]
The integral on the right-hand side can be computed \textit{mutatis mutandis}
by formally replacing $T^{\ast}f$ by $g$ in the proof of \eqref{ProofUnifying1} (with $\lambda = 1$) to obtain
\begin{equation}
\int_{X}\int_{0}^{\infty}t^{\gamma-1}\chi_{\{(x,t):t^{-\gamma/p}|g(x)| > 
1\}}\,\mathrm{d}t\,\mathrm{d}m(x)=\frac{1}{|\gamma|}\int_{X}|g(x)|^{p}%
\,\mathrm{d}m(x) \label{NewEquation2}%
\end{equation}
and \eqref{Statement1} follows.

Finally we deal with \eqref{Statement1.2} under \eqref{AssumptionMaximal} and \eqref{FinitenessLimit}. According to \eqref{Statement1}, it remains to show that
\begin{equation}\label{Limsupg}
	\limsup_{\lambda
\rightarrow L}\lambda^{p}\underset{t^{-\gamma/p}|T_{t}f(x)|>\lambda}{\int%
_{X}\int_{0}^{\infty}}t^{\gamma-1}\,\mathrm{d}t\,\mathrm{d}m(x) \leq \frac{1}{|\gamma|}\, \|g \|_{L^{p}(X,m)}^{p}.
\end{equation}
Indeed, the analog of \eqref{AnalogLimInf} when $\liminf_{\lambda \to L}$ is replaced by $\limsup_{\lambda \to L}$ is given by
\begin{equation}\label{LpEstimMaximalg}
 \limsup_{\lambda\rightarrow L}\chi_{\{(x,t):t^{-\gamma/p}|T_{t\lambda
^{-p/\gamma}}f(x)|>1\}} \leq \chi_{\{(x,t):t^{-\gamma/p}|g(x)| \geq 1\}}.
\end{equation}
Further, since (cf. \eqref{MaximalInclusion})
$$
	 \sup_{\lambda > 0}\chi_{\{(x,t):t^{-\gamma/p}|T_{t\lambda
^{-p/\gamma}}f(x)|>1\}} \leq \chi_{\{(x, t) : t^{-\gamma/p} T^*f(x) > 1\}}
$$
where \eqref{AssumptionMaximal} implies (cf. \eqref{ProofUnifying1} with $\lambda=1$)
$$
	\int_X \int_0^\infty t^{\gamma-1} \chi_{\{(x, t) : t^{-\gamma/p} T^*f(x) > 1\}} \, \mathrm{d} t \, \mathrm{d} m(x)= \frac{1}{|\gamma|}\int_{X}(T^{\ast}%
f(x))^{p}\,\mathrm{d}m(x) < \infty,
$$ 
one can apply reverse Fatou's lemma and \eqref{LpEstimMaximalg} to estimate
\begin{align*}
	\limsup_{\lambda \to L} \underset{t^{-\gamma
/p}|T_{t\lambda^{-p/\gamma}}f(x)|>1}{\int_{X}\int_{0}^{\infty}}t^{\gamma
-1}\,\mathrm{d}t\,\mathrm{d}m(x) & \leq \int_{X}\int_{0}^{\infty}t^{\gamma-1}\chi_{\{(x,t):t^{-\gamma/p}|g(x)| \geq 
1\}}\,\mathrm{d}t\,\mathrm{d}m(x) \\
& =\frac{1}{|\gamma|}\int_{X}|g(x)|^{p}%
\,\mathrm{d}m(x) 
\end{align*}
where the last step follows from Fubini's theorem as in the proof of \eqref{ProofUnifying1} (cf. also  \eqref{NewEquation2}). The proof of \eqref{Limsupg} is complete. 
\end{proof}

As a byproduct of Theorem \ref{ThmMaximal2}, all the applications given in Section \ref{SubsectionAppl} admit now counterparts in the spirit of Bourgain--Nguyen. In particular, following the methodology described in
Section \ref{ExampleBSY}, one can establish the following weighted version of \eqref{nb1} and \eqref{nb2}. 

\begin{corollary}
Let $f\in W^{1, p}(\mathbb{R}^{N}),\,p\in(1,\infty)$, and $\gamma<0$. Then
\begin{equation*}
\sup_{\delta>0}\underset{|x-y|^{-\gamma/p-1}|f(x)-f(y)|>\delta
}{\int_{\mathbb{R}^{N}}\int_{\mathbb{R}^{N}}}\frac{\delta^{p}}%
{|x-y|^{N-\gamma}}\,\mathrm{d}x\,\mathrm{d}y \leq -\frac{C(p, N)}{\gamma} \Vert\nabla f\Vert_{L^{p}(\mathbb{R}^{N})}^{p}
\end{equation*}
and
\begin{equation*}
	\lim_{\delta
\rightarrow0}\underset{|x-y|^{-\gamma/p-1}|f(x)-f(y)|>\delta}{\int%
_{\mathbb{R}^{N}}\int_{\mathbb{R}^{N}}}\frac{\delta^{p}}{|x-y|^{N-\gamma}%
}\,\mathrm{d}x\,\mathrm{d}y = - \frac{k(p, N)}{\gamma} \Vert\nabla f\Vert_{L^{p}(\mathbb{R}^{N})}^{p}
\end{equation*}
where $k(p, N)$ is given by \eqref{ConstantkpN}. 
\end{corollary}

In particular, setting $\gamma=-p$
in the previous result one recovers \eqref{nb1} and \eqref{nb2}.

\subsection{The limiting $L^{p}$ inequality for the Brezis--Van
Schaftingen--Yung spaces and the Bourgain--Nguyen functionals}

As we have pointed before (cf. footnote 5) our methodology allows us to
readily establish
\[
\bigg\|\frac{f(x)-f(y)}{|x-y|^{\frac{N}{p}}}\bigg\|_{L(p,\infty)(\mathbb{R}%
^{N}\times\mathbb{R}^{N})}\lesssim\Vert f\Vert_{L^{p}(\mathbb{R}^{N})}%
,\qquad1\leq p<\infty.
\]
Since the preprint of our paper appeared on line, Gu and Yung \cite[Theorem
1.1]{gupo} have shown that the converse estimate is also true, and therefore
established that%
\begin{equation}
\bigg\|\frac{f(x)-f(y)}{|x-y|^{\frac{N}{p}}}\bigg\|_{L(p,\infty)(\mathbb{R}%
^{N}\times\mathbb{R}^{N})}\asymp\Vert f\Vert_{L^{p}(\mathbb{R}^{N})}.
\label{maz1}%
\end{equation}
Moreover in \cite{gupo} it is also proven that
\[
\lim_{\lambda\rightarrow0}\lambda^{p}\mathcal{L}^{2N}\bigg(\bigg\{(x,y)\in
\mathbb{R}^{N}\times\mathbb{R}^{N}:\frac{f(x)-f(y)}{|x-y|^{\frac{N}{p}}%
}>\lambda\bigg\}\bigg)=2\kappa_{N}\,\Vert f\Vert_{L^{p}(\mathbb{R}^{N})}^{p},
\]
where $\kappa_{N}$ is the volume of the unit ball in $\mathbb{R}^{N}$. These
results can be considered counterparts to the limiting theorems for Gagliardo
seminorms due to Maz'ya--Shaposhnikova \cite{Maz}%

\[
\lim_{s \to 0+}s\int_{\mathbb{R}^{N}}\int_{\mathbb{R}^{N}}\frac{\left\vert
f(x)-f(y)\right\vert ^{p}}{\left\vert x-y\right\vert ^{N+sp}} \, \mathrm{d}x \, \mathrm{d}y= \frac{2 N}{p} \kappa_N  \int_{\mathbb{R}^{N}}\left\vert
f(x)\right\vert ^{p} \, \mathrm{d}x.
\]

Our aim in this section is to establish the corresponding version of
(\ref{maz1}) for the functionals $I_{\delta}$ (cf. \eqref{Idelta}) studied by Bourgain--Nguyen.
Namely, we prove the following

\begin{theorem}
Let $f\in L^{p}(\mathbb{R}^{N}),1\leq p<\infty.$ Then%
\begin{equation}
\underset{|x-y||f(x)-f(y)|>\delta}{\int_{\mathbb{R}^{N}}\int%
_{\mathbb{R}^{N}}}\frac{\delta^{p}}{|x-y|^{N+p}} \, \mathrm{d}x \, \mathrm{d}y\leq\frac{2^{p+1}%
\kappa_{N}}{p}\int_{\mathbb{R}^{N}}|f(x)|^{p}\,\mathrm{d}x, \qquad \text{for all} \qquad \delta>0, \label{GY1}%
\end{equation}
and%
\begin{equation}
\liminf_{\delta\rightarrow0}\underset{|x-y||f(x)-f(y)|>\delta
}{\int_{\mathbb{R}^{N}}\int_{\mathbb{R}^{N}}}\frac{\delta^{p}}{|x-y|^{N+p}%
} \, \mathrm{d}x \, \mathrm{d}y\geq\frac{4\kappa_{N}}{p}\int_{\mathbb{R}^{N}}|f(x)|^{p}\,\mathrm{d}x.
\label{GY2}%
\end{equation}

In particular,
\begin{align}
\lim_{\delta\rightarrow0}\underset{|x-y||f(x)-f(y)|>\delta}{\int%
_{\mathbb{R}^{N}}\int_{\mathbb{R}^{N}}}\frac{\delta^{p}}{|x-y|^{N+p}} \, \mathrm{d}x \, \mathrm{d}y  &
\simeq\sup_{\delta>0}\underset{|x-y||f(x)-f(y)|>\delta}{\int%
_{\mathbb{R}^{N}}\int_{\mathbb{R}^{N}}}\frac{\delta^{p}}{|x-y|^{N+p}} \, \mathrm{d}x \, \mathrm{d}y \nonumber\\
&  \simeq\left\Vert f\right\Vert _{L^{p}(\mathbb{R}^{N}).} \label{dep}%
\end{align}

\end{theorem}

\begin{proof}
We start by showing \eqref{GY1}. Using the triangle inequality
\[
|f(x)-f(y)|\leq|f(x)|+|f(y)|
\]
yields
\begin{align*}
\underset{|x-y||f(x)-f(y)|>\delta}{\int_{\mathbb{R}^{N}}\int_{\mathbb{R}^{N}}%
}\frac{\delta^{p}}{|x-y|^{N+p}}\,\mathrm{d}x\,\mathrm{d}y  &  \leq
2\int_{\mathbb{R}^{N}}\int_{|x-y|>\frac{\delta}{2|f(x)|}}\frac{\delta^{p}%
}{|x-y|^{N+p}}\,\mathrm{d}y\,\mathrm{d}x\\
&  =\frac{2^{p+1}\kappa_{N}}{p}\int_{\mathbb{R}^{N}}|f(x)|^{p}\,\mathrm{d}x.
\end{align*}
Consequently,
\[
\sup_{\delta>0}\underset{|x-y||f(x)-f(y)|>\delta}{\int_{\mathbb{R}^{N}}%
\int_{\mathbb{R}^{N}}}\frac{\delta^{p}}{|x-y|^{N+p}}\,\mathrm{d}%
x\,\mathrm{d}y\leq\frac{2^{p+1}\kappa_{N}}{p}\int_{\mathbb{R}^{N}}%
|f(x)|^{p}\,\mathrm{d}x.
\]

To prove \eqref{GY2}, we first assume that $f$ has compact support. Therefore,
there exists $R>0$ such that
\[
\text{supp }f\subseteq B_{R},
\]
where $B_{R}$ denotes the Euclidean ball centered at the origin and radius
$R$. For each $\delta>0$, set
\[
E_{\delta}:=\{(x,y):|x-y||f(x)-f(y)|>\delta\},
\]
and
\[
H_{\delta}:=E_{\delta}\cap\{(x,y):|y|>|x|\}.
\]
Note that $(x,y)\in H_{\delta}$ implies $x\in B_{R}$. Given $x\in B_{R}$, let
\[
H_{\delta,x}:=\{y:|y|>|x|,\quad|x-y||f(x)-f(y)|>\delta\}.
\]
Therefore we can write
\begin{equation}
\underset{H_{\delta}}{\iint}\frac{\delta^{p}}{|x-y|^{N+p}}\,\mathrm{d}%
x\,\mathrm{d}y=\int_{B_{R}}\int_{H_{\delta,x}}\frac{\delta^{p}}{|x-y|^{N+p}%
}\,\mathrm{d}y\,\mathrm{d}x. \label{Yung1}%
\end{equation}
Clearly
\[
H_{\delta,x,R}:=\{y:|y|>R,\quad|x-y||f(x)|>\delta\}\subseteq H_{\delta,x}.
\]
A simple change of variables yields,
\begin{align*}
\int_{H_{\delta,x}}\frac{\delta^{p}}{|x-y|^{N+p}}\,\mathrm{d}y  &  \geq
\int_{H_{\delta,x,R}}\frac{\delta^{p}}{|x-y|^{N+p}}\,\mathrm{d}y\\
&  =\int_{|y-x|>\frac{\delta}{|f(x)|}}\frac{\delta^{p}}{|x-y|^{N+p}%
}\,\mathrm{d}y\\
&  \hspace{1cm}-\int_{|y|<R,\,|x-y|>\frac{\delta}{|f(x)|}}\frac{\delta^{p}%
}{|x-y|^{N+p}}\,\mathrm{d}y\\
&  \geq\kappa_{N}\int_{\frac{\delta}{|f(x)|}}^{\infty}\frac{\delta^{p}%
}{t^{p+1}}\,\mathrm{d}t-\kappa_{N}\int_{\frac{\delta}{|f(x)|}}^{2R}%
\frac{\delta^{p}}{t^{p+1}}\,\mathrm{d}t\\
&  =2\kappa_{N}\frac{|f(x)|^{p}}{p}-\kappa_{N}\delta^{p}\frac{1}{p2^{p}R^{p}}.
\end{align*}
Hence, by \eqref{Yung1},%
\[
\underset{H_{\delta}}{\iint}\frac{\delta^{p}}{|x-y|^{N+p}}\,\mathrm{d}%
x\,\mathrm{d}y\geq\frac{2\kappa_{N}}{p}\int_{B_{R}}|f(x)|^{p}\,\mathrm{d}%
x-\kappa_{N}^{2}\delta^{p}\frac{R^{N-p}}{p2^{p}}.
\]
Taking limits as $\delta\rightarrow0$ we find
\begin{equation*}
\liminf_{\delta\rightarrow0}\underset{H_{\delta}}{\iint}\frac{\delta^{p}%
}{|x-y|^{N+p}}\,\mathrm{d}x\,\mathrm{d}y\geq\frac{2\kappa_{N}}{p}\int_{B_{R}%
}|f(x)|^{p}\,\mathrm{d}x. 
\end{equation*}
By symmetry,
\[
\liminf_{\delta\rightarrow0}\underset{E_{\delta}}{\iint}\frac{\delta^{p}%
}{|x-y|^{N+p}}\,\mathrm{d}x\,\mathrm{d}y\geq\frac{4\kappa_{N}}{p}\int_{B_{R}%
}|f(x)|^{p}\,\mathrm{d}x.
\]

More generally, if $f\in L^{p}(\mathbb{R}^{N})$ in \eqref{GY2} we can adapt
the approximation method used in \cite{gupo} when dealing with (\ref{maz1}).
For the sake of completeness we provide full details. Given $R>0,$ let
$f_{R}=f\chi_{B(0,R)}$ and $g_{R}=f-f_{R}$. For $\lambda,\delta>0$, we set
\[
A=\{(x,y):|f_{R}(x)-f_{R}(y)||x-y|>\delta(1+\lambda)\}
\]
and
\[
B=\{(x,y):|g_{R}(x)-g_{R}(y)||x-y|>\delta\lambda\}.
\]
Using the triangle inequality it follows that,
\[
A\backslash B\subseteq E_{\delta},
\]
and therefore
\begin{align}
\iint_{E_{\delta}}\frac{\delta^{p}}{|x-y|^{N+p}}\,\mathrm{d}x\,\mathrm{d}y  &
\geq\iint_{A}\frac{\delta^{p}}{|x-y|^{N+p}}\,\mathrm{d}x\,\mathrm{d}%
y-\iint_{A\cap B}\frac{\delta^{p}}{|x-y|^{N+p}}\,\mathrm{d}x\,\mathrm{d}%
y\nonumber\\
&  \geq\iint_{A}\frac{\delta^{p}}{|x-y|^{N+p}}\,\mathrm{d}x\,\mathrm{d}%
y-\iint_{B}\frac{\delta^{p}}{|x-y|^{N+p}}\,\mathrm{d}x\,\mathrm{d}y.
\label{hehehe}%
\end{align}
According to \eqref{GY1}, we can estimate
\[
\iint_{B}\frac{\delta^{p}}{|x-y|^{N+p}}\,\mathrm{d}x\,\mathrm{d}y\leq
\frac{2^{p+1}\kappa_{N}}{p\lambda^{p}}\Vert g_{R}\Vert_{L^{p}(\mathbb{R}^{N}%
)}^{p}.
\]
Inserting this information into \eqref{hehehe} we find
\[
\iint_{E_{\delta}}\frac{\delta^{p}}{|x-y|^{N+p}}\,\mathrm{d}x\,\mathrm{d}%
y\geq\iint_{A}\frac{\delta^{p}}{|x-y|^{N+p}}\,\mathrm{d}x\,\mathrm{d}%
y-\frac{2^{p+1}\kappa_{N}}{p\lambda^{p}}\Vert g_{R}\Vert_{L^{p}(\mathbb{R}%
^{N})}^{p}.
\]
Taking into account that $f_{R}$ has compact support, we now take limits, as
$\delta\rightarrow0,$ on both sides of the previous inequality$,$ and find
\begin{align*}
\liminf_{\delta\rightarrow0}\iint_{E_{\delta}}\frac{\delta^{p}}{|x-y|^{N+p}%
}\,\mathrm{d}x\,\mathrm{d}y  &  \geq\frac{1}{(1+\lambda)^{p}}\liminf
_{\delta\rightarrow0}\iint_{A}\frac{(\delta(1+\lambda))^{p}}{|x-y|^{N+p}%
}\,\mathrm{d}x\,\mathrm{d}y-\frac{2^{p+1}\kappa_{N}}{p\lambda^{p}}\Vert
g_{R}\Vert_{L^{p}(\mathbb{R}^{N})}^{p}\\
&  \geq\frac{1}{(1+\lambda)^{p}}\frac{4\kappa_{N}}{p}\Vert f_{R}\Vert
_{L^{p}(\mathbb{R}^{N})}^{p}-\frac{2^{p+1}\kappa_{N}}{p\lambda^{p}}\Vert
g_{R}\Vert_{L^{p}(\mathbb{R}^{N})}^{p}.
\end{align*}
Since
\[
\lim_{R\rightarrow\infty}\Vert f_{R}\Vert_{L^{p}(\mathbb{R}^{N})}^{p}=\Vert
f\Vert_{L^{p}(\mathbb{R}^{N})}^{p}\quad\text{and}\quad\lim_{R\rightarrow
\infty}\Vert g_{R}\Vert_{L^{p}(\mathbb{R}^{N})}^{p}=0
\]
it follows that
\[
\liminf_{\delta\rightarrow0}\iint_{E_{\delta}}\frac{\delta^{p}}{|x-y|^{N+p}%
}\,\mathrm{d}x\,\mathrm{d}y\geq\frac{1}{(1+\lambda)^{p}}\frac{4\kappa_{N}}%
{p}\Vert f\Vert_{L^{p}(\mathbb{R}^{N})}^{p}%
\]
for all $\lambda>0$. Passing to the limit as $\lambda\rightarrow0$ we conclude
that
\[
\liminf_{\delta\rightarrow0}\iint_{E_{\delta}}\frac{\delta^{p}}{|x-y|^{N+p}%
}\,\mathrm{d}x\,\mathrm{d}y\geq\frac{4\kappa_{N}}{p}\Vert f\Vert
_{L^{p}(\mathbb{R}^{N})}^{p}.
\]

\end{proof}

\begin{remark}
In the setting of Bourgain--Nguyen functionals it could be of interest to
consider the *interpolation* functionals%
\[
I_{\delta,s}(f)=\underset{|x-y|^{1-s}|f(x)-f(y)|>\delta}{\int_{\mathbb{R}%
^{N}}\int_{\mathbb{R}^{N}}}\frac{\delta^p}{|x-y|^{N+p}}\,\mathrm{d}x\,\mathrm{d}y, \qquad s \in [0, 1], 
\]
and let%
\[
\left\Vert f\right\Vert _{BN_{p}^{s} (\mathbb{R}^N)}^{p}=\sup_{\delta>0} I_{\delta
,s}(f).
\]
Here the case $s=1,$ corresponds to the recovery of the gradient (cf. \eqref{nb1} and \eqref{nb2}) while the case $s=0$ corresponds to the recovery of $L^{p}$
norms (cf. (\ref{dep})). Furthermore, a simple adaptation of the method of proof of Theorem \ref{ThmMarkao} enables us to show that
$$
	C^s_p(\mathbb{R}^N) \subset BN^s_p(\mathbb{R}^N), \qquad s \in (0, 1], \qquad p \in (1, \infty). 
$$
\end{remark}

\end{document}